\documentclass[reqno, a4paper, 10pt,oneside]{amsart}
\makeatletter
\def\specialsection{\@startsection{section}{1}%
	\z@{\linespacing\@plus\linespacing}{.5\linespacing}%
	{\normalfont}}
\def\section{\@startsection{section}{1}%
	\z@{.7\linespacing\@plus\linespacing}{.5\linespacing}%
	{\normalfont\scshape\bfseries}}
\makeatother

\usepackage[top=1in, bottom=1.4in, inner=1in, outer=1in, includehead, showframe=false]{geometry}
\usepackage[algo2e, linesnumbered, noline, noend, ruled]{algorithm2e}
\usepackage[foot]{amsaddr}
\usepackage{amsfonts}
\usepackage{amsmath}
\usepackage{amssymb}
\usepackage{amsthm}
\usepackage[english]{babel}
\usepackage{booktabs}
\usepackage{caption}
\usepackage{braket}
\usepackage{float}
\usepackage[T1]{fontenc}
\usepackage{graphicx}
\usepackage[utf8]{inputenc}
\usepackage{mathtools}
\usepackage{nicefrac}
\usepackage[defaultlines=2,all]{nowidow}
\usepackage[per-mode=symbol, group-digits=integer]{siunitx}
\usepackage[table,dvipsnames]{xcolor}
\usepackage{enumitem}
\usepackage{subcaption}
\usepackage{stmaryrd}
\usepackage{textcomp}
\usepackage{todonotes}

\usepackage{cite}
\usepackage{tikz}
\usepackage[allcolors=cyan,colorlinks]{hyperref}

\usetikzlibrary{patterns}

\tikzset{%
	dimen/.style={<->,>=latex,thin,every rectangle node/.style={fill=white,midway,font=\sffamily}},
}

\definecolor{myblue}{RGB}{0, 107, 164}
\definecolor{myorange}{RGB}{255, 128, 14}
\definecolor{grayone}{RGB}{171, 171, 171}

\SetCommentSty{mycommfont}
\SetKwComment{Comment}{\# }{}

\theoremstyle{plain}
\newtheorem{Theorem}{Theorem}[section]

\theoremstyle{definition}
\newtheorem{Definition}[Theorem]{Definition}
\newtheorem{Remark}[Theorem]{Remark}

\newtheorem*{Example*}{Example}
\newtheorem*{Remark*}{Remark}

\setenumerate{label=(\arabic*), ref=(\arabic*)}

\newcommand{\qe}[1]{``#1''}
\newcommand{\tr}[1]{\mathrm{tr} #1 }

\newcommand{\R}{\mathbb{R}}

\newcommand{\transposed}{^\top}

\newcommand{\norm}[2]{\left\lvert \left\lvert #1 \right\rvert \right\rvert_{#2}}
\newcommand{\abs}[1]{\left\lvert #1 \right\rvert}
\newcommand{\dmeas}[1]{\ \mathrm{d}#1}

\newcommand{\integral}[1]{\int_{#1}}

\makeatletter
\DeclareOldFontCommand{\rm}{\normalfont\rmfamily}{\mathrm}
\DeclareOldFontCommand{\sf}{\normalfont\sffamily}{\mathsf}
\DeclareOldFontCommand{\tt}{\normalfont\ttfamily}{\mathtt}
\DeclareOldFontCommand{\bf}{\normalfont\bfseries}{\mathbf}
\DeclareOldFontCommand{\it}{\normalfont\itshape}{\mathit}
\DeclareOldFontCommand{\sl}{\normalfont\slshape}{\@nomath\sl}
\DeclareOldFontCommand{\sc}{\normalfont\scshape}{\@nomath\sc}
\makeatother

\numberwithin{equation}{section}

\newcommand*\closure[1]{%
	\hbox{%
		\vbox{%
			\hrule height 0.5pt 
			\kern0.5ex
			\hbox{%
				\ensuremath{#1}%
			}%
		}%
	}%
} 

\newcommand{\eps}{\varepsilon}

\newcommand{\Dsf}{{\mathsf{D}}}

\providecommand{\Cd}{{\mathcal D}}

\providecommand{\Cp}{{\mathcal P}}

\providecommand*{\abs}[1]{\left|{#1}\right|} 

\begin{document}
{\noindent\footnotesize This is a post-peer-review, pre-copyedit version of an article published in Structural and Multidisciplinary Optimization. The final version is available online at \url{https://doi.org/10.1007/s00158-023-03653-2}.
}	

\title[Quasi-Newton Methods for Topology Optimization Using a Level-Set Method]{Quasi-Newton Methods for Topology Optimization Using a Level-Set Method}
\author{Sebastian Blauth$^{*,1}$}
\address{$^*$ Corresponding Author}
\address{$^1$ Fraunhofer ITWM, Kaiserslautern, Germany}
\email{\href{mailto:sebastian.blauth@itwm.fraunhofer.de}{sebastian.blauth@itwm.fraunhofer.de}}
\author{Kevin Sturm$^2$}
\address{$^2$ Institute of Analysis and Scientific Computing, TU Wien, Vienna, Austria}
\email{\href{mailto:kevin.sturm@tuwien.ac.at}{kevin.sturm@tuwien.ac.at}}

\begin{abstract}
	The ability to efficiently solve topology optimization problems is of great importance for many practical applications. Hence, there is a demand for efficient solution algorithms. In this paper, we propose novel quasi-Newton methods for solving PDE-constrained topology optimization problems. Our approach is based on and extends the popular solution algorithm of Amstutz and Andr\"a (A new algorithm for topology optimization using a level-set method, Journal of Computational Physics, 216, 2006). To do so, we introduce a new perspective on the commonly used evolution equation for the level-set method, which allows us to derive our quasi-Newton methods for topology optimization. We investigate the performance of the proposed methods numerically for the following examples: Inverse topology optimization problems constrained by linear and semilinear elliptic Poisson problems, compliance minimization in linear elasticity, and the optimization of fluids in Navier-Stokes flow, where we compare them to current state-of-the-art methods. Our results show that the proposed solution algorithms significantly outperform the other considered methods: They require substantially less iterations to find a optimizer while demanding only slightly more resources per iteration. This shows that our proposed methods are highly attractive solution methods in the field of topology optimization.
	
	
\end{abstract}

\maketitle

\section*{Acknowledgments and Funding}

No funding was received.


\vspace{-0.5cm}

\section{Introduction}
\label{sec:introduction}

Topology optimization is concerned with the optimization of a domain by altering its geometrical features. Whereas in shape optimization only the boundary of a domain is variable, topology optimization considers the addition or removal of material to the geometry. Originally introduced in the context of solid mechanics, topology optimization has been considered for many practical applications, e.g., compliance minimization in elasticity \cite{a_ESKOSC_1994a, Allaire2005Structural, Amstutz2006new}, design optimization in the context of fluid mechanics \cite{Borrvall2003Topology, N.Sa2016Topological}, electrical machines \cite{Gangl2012Topology}, as well as the solution of inverse problems \cite{a_CALANO_2015a, Hintermueller2008Electrical, a_LAHIFRSC_2013a, Beretta2017reconstruction} and the modeling and simulation of fracture evolution \cite{Xavier2018simplified, Xavier2017Topological}.

The goal of topology optimization is to optimize a cost functional $J$ depending on the set $\Omega$, which plays the role of the design variable. To achieve this goal, usually all possible designs are assumed to belong to a fixed domain $\Dsf\subset \R^d$, $d\in \mathbb{N}_{>0}$, where $\mathbb{N}$ denotes the set of positive integers, which is referred to as hold-all domain. The understanding how a functional varies under perturbations of $\Omega$ is crucial for the development of numerical methods. Of particular importance are perturbations obtained by removing or adding a small inclusion $\omega_\eps(x_0)$ at $x_0\in \Omega$ or $x_0\in \Dsf\setminus \overline\Omega$ of size $\eps$. The first order variation of the functional under this perturbation is called the \qe{topological derivative}. This leads to the definition of the perturbed domains 
\begin{align}\label{eq:omega_epss}
\Omega_\eps & := \Omega_\eps(x_0,\omega) :=  \begin{cases}
\Omega\cup \omega_\eps(x_0) &   \text{ for }x_0\in \Dsf\setminus \overline\Omega, \\
\Omega\setminus \overline{\omega_\eps(x_0)} &  \text{ for } x_0 \in \Omega,
\end{cases}  
\end{align}
and $\omega_\eps(x_0):= x_0 + \eps \omega$ with $\omega\subset \R^d$ being a simply connected domain with $0\in \omega$. Then, the 
topological derivative can be defined by 
\begin{equation}\label{def:topo_derivative_trans}
DJ(\Omega)(x_0,\omega) = \lim_{\eps\searrow0} \frac{J(\Omega_\eps(x_0,\omega)) - J(\Omega)}{f(\eps)},
\end{equation}
where $f$ satisfies $\lim_{\eps \searrow 0} f(\eps) = 0$. In fact, the topological derivative can also depend on the inclusion $\omega$, which we omit in our notation for simplicity and just write $DJ(\Omega)$. For many shape functionals in practice the function $f$ is given by the volume of $\omega_\eps$, however, there are problems where this is not the case, e.g., when Dirichlet boundary conditions are imposed on the inclusion boundary $\partial \omega_\eps$ (cf.~\cite{Amstutz2022introduction}).

The topological derivative was first formally introduced as the bubble method in \cite{a_ESKOSC_1994a} and later mathematically justified in \cite{a_SOZO_1999a,a_GAGUMA_2001a}. Topological derivatives have been established for a variety of PDE constrained functionals and we refer to the monographs \cite{b_NOSO_2013a,b_NOSO_2020a} for further information. For the computation of the topological derivative there exist several methods, for instance a direct approach \cite{b_NOSO_2013a}, where first the expansion of the state variable is computed and afterwards the topological expansion of the cost function is derived via Taylor's formula. Lagrangian approaches provide another way to compute topological derivatives, e.g., a Lagrangian method using an averaged adjoint equation is presented in \cite{a_ST_2020a}, a method using a perturbed adjoint equation in conjunction with an unperturbed state equation is found in \cite{Amstutz2006new}, and another Lagrangian method using an unperturbed adjoint equation in \cite{a_GAST_2020a}. We refer to \cite{Baumann2022Adjoint} for a comparison of Lagrangian methods and their advantages and disadvantages.  

From a numerical perspective, one can either use the topological derivative directly to find the location of an optimal design (see, e.g., \cite{Hintermueller2008Electrical,b_NOSO_2020a}) or one can use a level-set approach \cite{Amstutz2006new} in an iterative fashion to find the optimal topology. However, it should be noted that only the topological derivative is capable to create new holes while the level-set approach only allows to close already present holes. Additionally, Newton-type algorithms to find circular inclusions have been introduced using higher order topological expansions (cf.~\cite{a_LAHIFRSC_2013a,a_CALANO_2015a} or \cite[Chapter 10]{b_NOSOZO_2019a}). This leads to rapidly converging algorithms, however, it requires a combinatorial search, which makes it numerically expensive and also the study of nonlinear problems remains an open problem.

In the related field of shape optimization, the development, modification, and analysis of efficient solution algorithms has received lots of attention in recent years, e.g., in \cite{Blauth2021Nonlinear,Blauth2022Shape} and \cite{Schulz2016Efficient}, where nonlinear conjugate gradient and quasi-Newton methods for shape optimization have been proposed, in \cite{Blauth2022Space}, where a space mapping technique for shape optimization is presented, and in \cite{Deckelnick2022novel}, where a $W^{1,\infty}$ approach for shape optimization is introduced. 

In this paper, we follow these developments and propose novel quasi-Newton methods for solving PDE constrained topology optimization problems based on the topological derivative. Our approach is based on and extends the popular level-set approach of Amstutz and Andr\"a \cite{Amstutz2006new}. In particular, we first provide a new perspective on the evolution equation for the level-set function given in \cite{Amstutz2006new} to derive a gradient descent-type algorithm. This algorithm is the foundation for the quasi-Newton methods we propose afterwards. In particular, we derive a limited-memory BFGS method for topology optimization. We investigate the novel methods numerically for several problem classes: Inverse topology optimization problems constrained by linear and semilinear Poisson problems, compliance minimization in linear elasticity, and the optimization of fluids in Navier-Stokes flow. We compare our methods behavior to the widely popular method of Amstutz and Andr\"a \cite{Amstutz2006new} as well as a simpler convex combination method from \cite{Gangl2021Asymptotic}. The results show that the novel quasi-Newton methods usually have a significantly better convergence behavior than the other solution algorithms. Particularly, they require substantially less iterations to find an optimizer of the problems while demanding only slightly more numerical resources per iteration. This makes our proposed quasi-Newton methods highly attractive for solving topology optimization problems. 

This paper is structured as follows. In Section~\ref{sec:preliminaries}, we briefly recall basic results from topology optimization as well as the level-set method for topology optimization and the solution algorithms from \cite{Amstutz2006new} and \cite{Gangl2021Asymptotic}, which are widely-used in the literature. In Section~\ref{sec:quasi_newton}, we present a new perspective on the evolution equation for the level-set function. This allows us to derive a gradient descent algorithm for topology optimization and, afterwards, propose novel quasi-Newton methods for topology optimization. 
Finally, we investigate the methods numerically and compare our quasi-Newton methods to the current state-of-the-art methods in Section~\ref{sec:numerics}. 

\section{Preliminaries}
\label{sec:preliminaries}

\subsection{Topological Sensitivity Analysis}
\label{ssec:topological_sensitivity_analysis}

In this section we recall basics on the topological derivative. We start by recalling the definition of the topological derivative (cf.~\cite{Amstutz2022introduction}).

\begin{Definition}
    Let $\omega$ be a simply connected, bounded, and open subset of $\R^d$, $d\in \mathbb{N}_{>0}$, with $0\in \omega$, let $\Dsf$ be a bounded hold-all domain, and let $\mathcal{P}(\Dsf)$ be the power set of $\Dsf$, i.e., $\mathcal{P}(\Dsf) = \Set{\Omega\subset \mathbb{R}^d:\; \Omega \subset \Dsf}$. Let $J\colon \mathcal{P}(\Dsf) \to \R$ be a shape functional. We say that $J$ has a topological derivative at $\Omega \in \mathcal{P}(\Dsf)$ and at the point $x_0\in \Dsf$ w.r.t.~$\omega$ if there exists some function $f\colon \R_{>0} \to \R_{>0}$, where $\R_{>0}$ denotes positive real numbers, with $\lim_{\eps \searrow 0} f(\eps) = 0$ so that the following limit exists
	\begin{equation*}
		DJ(\Omega)(x_0) := \lim_{\eps \searrow 0} \frac{J(\Omega_\eps) - J(\Omega)}{f(\eps)},
	\end{equation*}
	where the perturbed domain $\Omega_\eps$ is defined by
	\begin{equation*}
	\Omega_\eps = \Omega_\eps(x_0, \omega) = \begin{cases}
	\Omega \setminus \closure{\omega_\eps(x_0)} \quad &\text{for } x_0 \in \Omega,\\
	\Omega \cup \omega_\eps(x_0) \quad &\text{for } x_0 \in \Dsf \setminus \closure{\Omega},
	\end{cases}
	\end{equation*}
	and $\omega_\eps(x_0) := x_0 + \eps \omega$.
\end{Definition}

We follow the approach of \cite{Amstutz2006new} and represent a set $\Omega \subset \Dsf$ with the help of a continuous level set function $\psi \colon \Dsf \to \R$ as follows
\begin{equation*}
\begin{aligned}
\psi(x) < 0 \quad &\Leftrightarrow \quad x \in \Omega, \\
\psi(x) > 0 \quad &\Leftrightarrow \quad x \in \Dsf \setminus \closure{\Omega}, \\
\psi(x) = 0 \quad &\Leftrightarrow \quad x \in \partial \Omega \setminus \partial \Dsf.
\end{aligned}
\end{equation*}
We write $\Omega_\psi := \Omega$ for a domain $\Omega \subset \Dsf$ which is represented by the level-set function $\psi$.

\begin{Definition} 
	Let $J\colon \Cp(\Dsf)\to \R$ be a shape functional, let $\Omega \in \mathcal{P}(\Dsf)$ be an open set, and let $\Gamma = \partial \Omega \setminus \partial \Dsf$. Assume that the topological derivative $DJ(\Omega)(x)$ exists for all $x\in \Dsf\setminus \Gamma$. Then, we define the generalized topological derivative by 
	\begin{equation}
	\Cd J(\Omega)(x) := \begin{cases}
	-DJ(\Omega)(x) & \text{ for } x\in \Omega, \\
	DJ(\Omega)(x) & \text{ for } x\in \Dsf\setminus\overline\Omega.
	\end{cases}
	\end{equation}
\end{Definition}

The idea behind the generalized topological derivative originates from the following observation. Let $\psi\colon \Dsf\to \R$ be a level set function representing $\Omega$. If there is a constant $c>0$ such that
\begin{equation}\label{eq:opt_levelset}
\Cd J(\Omega_\psi)(x) = c\psi(x) \quad \text{ for } x\in \Dsf,
\end{equation}
then 
\begin{equation}
\label{eq:fonc}
DJ(\Omega)(x)\ge 0 \quad \text{ for all } x\in \Dsf\setminus\partial \Omega,
\end{equation}
which is the necessary condition for $\Omega$ to be optimal (cf.~\cite{Amstutz2006new}). This observation is the starting point of the popular solution algorithm from \cite{Amstutz2006new}, which we present in the next section.

\subsection{Topology Optimization Algorithms Using a Level-Set Method}
\label{ssec:optimization_algorithms}

Based on the discussion in the previous section, a solution algorithm for topology optimization has been derived in \cite{Amstutz2006new}, which we briefly recall in the following. For the derivation, we introduce a fictitious time $t$ and consider a family of domains $\Omega(t)$ represented by a level-set function $\psi\colon [0,T] \times \Dsf \to \R$. To derive a solution algorithm, the idea is to establish an evolution equation for the level-set function $\psi$ which ensures that $\Omega(t)$ converges to a minimizer for $t\to \infty$ and to discretize this evolution equation. A natural idea would be to evolve the level-set function according to the generalized topological derivative $\mathcal{D}J(\Omega(t))$, so that
\begin{equation*}
\frac{\partial \psi(t, \cdot)}{\partial t} = \mathcal{D}J(\Omega(t)), \quad t\ge 0.
\end{equation*}
For the sake of better readability, we drop the dependence on the fictitious time $t$ in the domain $\Omega(t)$ throughout the rest of this paper.

In \cite{Amstutz2006new}, the authors note that, in general, using this formulation does not guarantee convergence since the generalized topological derivative $\mathcal{D}J(\Omega)$ does not vanish for an optimal geometry (cf.~\eqref{eq:fonc}). Instead, the main idea of \cite{Amstutz2006new} is to use a modified version of this equation, where the topological derivative is projected to the orthogonal complement of $\psi$ in $L^2(\Dsf)$, which is written as
\begin{equation}
\label{eq:projected_amstutz}
\frac{\partial \psi}{\partial t} = P_{\psi^\perp}(\mathcal{D}J(\Omega)),
\end{equation}
where the operator $P_{\psi^\perp}$ is defined as
\begin{equation*}
P_{\psi^\perp}(a) = a - \frac{(a, \psi)}{\norm{\psi}{L^2(\Dsf)}^2} \psi.
\end{equation*}
Here, $(a, b) = (a,b)_{L^2(\Dsf)}$ denotes the $L^2(\Dsf)$ scalar product between $a, b\in L^2(\Dsf)$. Now, if the right-hand side of \eqref{eq:projected_amstutz} vanishes, i.e., if $P_{\psi^\perp}(\mathcal{D}J(\Omega)) = 0$, then it holds that there exists some $\alpha \in \mathbb{R}$ so that $\mathcal{D}J(\Omega) = \alpha \psi$. If $\alpha > 0$, then the necessary optimality conditions \eqref{eq:fonc} for the optimization problem are satisfied and the geometry $\Omega$ described by the corresponding level-set function $\psi$ is a local minimizer.
For this reason, we also consider the norm of the projected topological derivative, i.e.,
\begin{equation}
	\label{eq:norm_projected_derivative}
	\norm{P_\psi^\perp(\mathcal{D}J(\Omega))}{L^2(\Dsf)}
\end{equation}
as a second convergence criterion for our numerical experiments in Section~\ref{sec:numerics}.

\begin{algorithm2e}[!t]
	\KwIn{Initial level-set function $\psi_0$ and corresponding geometry $\Omega_0$, initial step size $\lambda_0 = 1$, stopping tolerance $\tau > 0$, maximum number of iterations $k_{\mathrm{max}} \in \mathbb{N}$ \label{line:input}}
	\For{$k=0,1,2,\dots, k_\mathrm{max}$}{
		Compute the generalized topological derivative $\mathcal{D}J(\Omega_k)$\\
		Compute $\theta_k = \arccos\left( \frac{(\psi_k, \mathcal{D}J(\Omega_k))}{\norm{\psi_k}{L^2(D)} \norm{\mathcal{D}J(\Omega_k)}{L^2(D)}} \right)$\\
		\If{$\theta_k < \tau$}{
			Stop with approximate minimizer $\Omega_k$\\
		}
		Set $\lambda_k = \min\left( 1, 1.5 \lambda_{k-1} \right)$\\
		\While{$J\left(\frac{1}{\sin(\theta_k)} \sin\left( (1- \lambda_k) \theta_k \right) \psi_k + \sin\left( \lambda_k \theta_k \right) \frac{\mathcal{D}J(\Omega_k)}{\norm{\mathcal{D}J(\Omega_k)}{L^2(D)}}\right) > J(\psi_k)$ \label{line:line_search}}{
			Decrease the step size $\lambda_k = \frac{\lambda_k}{2}$\\
		}
		Update the level-set function: $\psi_{k+1} = \frac{1}{\sin(\theta_k)} \sin\left( (1- \lambda_k) \theta_k \right) \psi_k + \sin\left( \lambda_k \theta_k \right) \frac{\mathcal{D}J(\Omega_k)}{\norm{\mathcal{D}J(\Omega_k)}{L^2(D)}}$
	}
	\caption{Euler's method on the sphere algorithm for topology optimization.}
	\label{algo:euler}
\end{algorithm2e}

For the numerical solution, it is proposed in \cite{Amstutz2006new} to discretize \eqref{eq:projected_amstutz} with Euler's scheme on the sphere using a step size $\lambda$, leading to
\begin{equation}
	\label{eq:update_euler}
	\psi_{k+1} = \frac{1}{\sin(\theta_k)} \sin\left( (1- \lambda_k) \theta_k \right) \psi_k + \sin\left( \lambda_k \theta_k \right) \frac{\mathcal{D}J(\Omega_k)}{\norm{\mathcal{D}J(\Omega_k)}{L^2(\Dsf)}}, 
\end{equation}
where $\theta_k$ is the angle between $\psi_k$ and $DJ(\Omega_k)$, i.e. 
\begin{equation*}
\theta_k = \arccos\left( \frac{(\psi_k, \mathcal{D}J(\Omega_k))}{\norm{\psi_k}{L^2(\Dsf)} \norm{\mathcal{D}J(\Omega_k)}{L^2(D)}} \right).
\end{equation*}
Note that the iteration is terminated if the angle between generalized topological derivative and level-set function is sufficiently small, so that $P_{\psi^\perp}(\mathcal{D}J(\Omega)) \approx 0$ and the necessary optimality conditions are satisfied approximately.

The resulting numerical algorithm, which is analyzed in \cite{Amstutz2011Analysis}, can be seen in Algorithm~\ref{algo:euler}. By slight misuse of notation, we write $J(\psi)$ instead of $J(\Omega)$, where $\psi$ is the level-set function representing the domain $\Omega$, in line~\ref{line:line_search} of Algorithm~\ref{algo:euler} for the sake of better readability.

\begin{algorithm2e}[!t]
	\KwIn{Initial level-set function $\psi_0$ and corresponding geometry $\Omega_0$, initial step size $\lambda_0 = 1$, stopping tolerance $\tau > 0$, maximum number of iterations $k_{\mathrm{max}} \in \mathbb{N}$}
	\For{$k=0,1,2,\dots, k_\mathrm{max}$}{
		Compute the generalized topological derivative $\mathcal{D}J(\Omega_k)$\\
		Compute $\theta_k = \arccos\left( \frac{(\psi_k, \mathcal{D}J(\Omega_k))}{\norm{\psi_k}{L^2(D)} \norm{\mathcal{D}J(\Omega_k)}{L^2(\Dsf)}} \right)$\\
		\If{$\theta_k < \tau$}{
			Stop with approximate minimizer $\Omega_k$\\
		}
		Set $\lambda_k = \min\left( 1, 2 \lambda_{k-1} \right)$\\
		\While{$J\left(\lambda_k \frac{\mathcal{D}J(\Omega_k)}{\norm{\mathcal{D}J(\Omega_k)}{L^2(\Dsf)}} + (1 - \lambda_k) \psi_k\right) > J(\psi_k)$ \label{line:line_search_convex}}{
			Decrease the step size $\lambda_k = \frac{\lambda_k}{2}$\\
		}
		Update the level-set function: $\psi_{k+1} = \lambda_k \frac{\mathcal{D}J(\Omega_k)}{\norm{\mathcal{D}J(\Omega_k)}{L^2(\Dsf)}} + (1 - \lambda_k) \psi_k$
	}
	\caption{Convex combination algorithm for topology optimization.}
	\label{algo:convex_combination}
\end{algorithm2e}

The main idea of Algorithm~\ref{algo:euler} is to iteratively use a convex combination of generalized topological derivative and level-set function in order to reach a minimum of the optimization problem, where the weights for the convex combination are chosen according to Euler's method on the sphere (cf.~\cite{Amstutz2006new}). A simpler idea was used in \cite{Gangl2021Asymptotic}, where the authors consider the following convex combination to evolve the level-set function
\begin{equation}
\label{eq:update_convex}
\psi_{k+1} = \lambda \frac{\mathcal{D}J(\Omega_k)}{\norm{\mathcal{D}J(\Omega_k)}{L^2(\Dsf)}} + (1 - \lambda) \psi_k,
\end{equation}
where the parameter $\lambda$ plays the role of a step size. For a fixed step size $\lambda > 0$, it is easy to see that if the method becomes stationary, i.e., if $\psi_{k+1} = \psi_k$, then we have that $\psi_k = \nicefrac{\mathcal{D}J(\Omega_k)}{\norm{\mathcal{D}J(\Omega_k)}{L^2(\Dsf)}}$ and, in particular, the necessary optimality conditions \eqref{eq:opt_levelset} are satisfied. This idea gives rise to the Algorithm~\ref{algo:convex_combination}. 

Note that particularly the solution method presented in Algorithm~\ref{algo:euler} is widely popular and represents the state-of-the-art algorithm for solving topology optimization problems with a level-set function.

\section{Quasi-Newton Methods for Topology Optimization}
\label{sec:quasi_newton}

In this section, we present novel quasi-Newton methods for topology optimization. To do so, we first take a different perspective on equation~\eqref{eq:projected_amstutz} and formulate a new gradient descent method for topology optimization. This enables us to define (limited memory) BFGS methods for topology optimization afterwards.

\subsection{A Novel Perspective on the Level-Set Evolution Equation}

Before we can introduce quasi-Newton methods for topology optimization, we note that the algorithmic frameworks presented in Algorithms~\ref{algo:euler} and~\ref{algo:convex_combination} are not suitable for defining such methods. The reason for this is that the update rules for the level-set function given in \eqref{eq:update_euler} and \eqref{eq:update_convex} do consider a convex combination of the level-set function and the generalized topological derivative, which is different to the classical form of descent methods, where the update of the design variables (the level-set function) is performed by subtracting the gradient (topological derivative), scaled by an appropriate step size, from the current iterate. 

To remedy this problem, we start by considering the continuous equation~\eqref{eq:projected_amstutz} for evolving the level-set function from \cite{Amstutz2006new}. Instead of using the elaborate approach of \cite{Amstutz2006new}, we discretize \eqref{eq:projected_amstutz} by an explicit Euler method, which yields the discretized equation
\begin{equation}
\label{eq:gradient_flow_discrete}
\frac{\psi_{k+1} - \psi_k}{\Delta t} = P_{\psi_k^\perp}(\mathcal{D}J(\Omega_k)) \qquad \Leftrightarrow \qquad \psi_{k+1} = \psi_k + \Delta t P_{\psi_k^\perp}(\mathcal{D}J(\Omega_k))
\end{equation}
where $k$ denotes the current time step. For discretizing equation \eqref{eq:projected_amstutz}, one usually would consider small time steps $\Delta t$ converging to $0$, which would yield a so-called gradient flow method. However, we change our viewpoint and interpret \eqref{eq:gradient_flow_discrete} as gradient descent method, where the time step $\Delta t$ now plays the role of a step size. The benefit of this interpretation is that we can (potentially) make use of the convergence behavior of the gradient descent method and use large step sizes when appropriate, reducing the computational cost of our algorithm. Therefore, we can now interpret $g_k = -P_{\psi_k^\perp}(\mathcal{D}J(\Omega_k))$ as the \qe{gradient} associated to our topology optimization problem. The resulting optimization algorithm is presented in Algorithm~\ref{algo:gradient_descent}. Note, that the main benefit of this method is that it follows the \qe{standard} form of a gradient descent method, which makes it amenable to define quasi-Newton methods, which we do in the next section. Additionally, our numerical experiments in Section~\ref{sec:numerics} show, that Algorithm~\ref{algo:gradient_descent} can also yield faster convergence compared to Algorithms~\ref{algo:euler} and~\ref{algo:convex_combination}.

\begin{algorithm2e}[!t]
	\KwIn{Initial level-set function $\psi_0$, initial step size $\lambda_0 = 1$, stopping tolerance $\tau > 0$, maximum number of iterations $k_{\mathrm{max}} \in \mathbb{N}$}
	\For{$k=0,1,2,\dots, k_\mathrm{max}$}{
		Compute the generalized topological derivative $\mathcal{D}J(\Omega_k)$\\
		Set $g_k = -P_{\psi_k^\perp}(\mathcal{D}J(\Omega_k)) = - \left( \mathcal{D}J(\Omega_k) - \frac{\left( \mathcal{D}J(\Omega_k), \psi_k \right)}{\norm{\psi_k}{L^2(\Dsf)}^2} \psi_k \right)$\\
		Compute $\theta_k = \arccos\left( \frac{(\psi_k, \mathcal{D}J(\Omega_k))}{\norm{\psi_k}{L^2(\Dsf)} \norm{\mathcal{D}J(\Omega_k)}{L^2(\Dsf)}} \right)$\\
		\If{$\theta_k < \tau$}{
			Stop with approximate minimizer $\Omega_k$\\
		}
		Set $\lambda_k = 2 \lambda_{k-1}$\\
		\While{$J\left(\psi_k - \lambda_k g_k\right) > J(\psi_k)$}{
			Decrease the step size $\lambda_k = \frac{\lambda_k}{2}$\\
		}
		Update the level-set function: $\psi_{k+1} = \psi_k - \lambda_k g_k$
	}
	\caption{Gradient descent method for topology optimization.}
	\label{algo:gradient_descent}
\end{algorithm2e}

\subsection{A Limited Memory BFGS Method for Topology Optimization}

With the gradient descent method described in the previous section, we now focus our attention to quasi-Newton methods for topology optimization, which can now be derived analogously to the finite-dimensional case (see, e.g., \cite{Nocedal2006Numerical, Kelley1999Iterative}). To do so, we introduce the functions $s_k = \psi_{k+1} - \psi_k$ and $y_k = g_{k+1} - g_k$. The quasi-Newton methods rely on the so-called secant equation, which in our setting can be written as
\begin{equation*}
B_{k+1} s_k = y_k,
\end{equation*}
where $B_{k+1}$ is an isomorphism from $L^2(\Dsf)$ to $L^2(\Dsf)$ which can be seen as approximation of the Hessian, and we denote its inverse by $H_{k+1}$. In the following, we will describe a BFGS method for topology optimization based on $H_k$, the inverse of the Hessian approximation, which makes it easier to derive the limited-memory version of the method which we have implemented in the software package cashocs \cite{Blauth2021cashocs, Blauth2023Version}. In particular, the search direction for the BFGS method is given by
\begin{equation}
\label{eq:search_direction_bfgs}
p_k = - H_k g_k
\end{equation}
and the update formula for $H_k$ is given by
\begin{equation}
\label{eq:inverse_hessian_update}
H_{k+1} = \left( \mathrm{Id}_{L^2(\Dsf)} - \frac{s_k \otimes y_k}{(y_k, s_k)_{L^2(\Dsf)}} \right) H_k \left( \mathrm{Id}_{L^2(\Dsf)} - \frac{y_k \otimes s_k}{(s_k, y_k)_{L^2(\Dsf)}} \right) + \frac{s_k \otimes s_k}{(y_k, s_k)_{L^2(\Dsf)}},
\end{equation}
where $\otimes$ denotes the outer product of $L^2(\Dsf)$, i.e., $(a \otimes b) c = (b, c)_{L^2(\Dsf)} a$. 


To avoid storing large dense matrices as discretizations of the operator $H_{k}$, we employ a limited memory BFGS method for our numerical implementation which is shown in Algorithm~\ref{algo:lbfgs}. Particularly, the limited memory BFGS method only requires us to additionally store $s_k$ and $y_k$ to compute the application of $H_{k}$ to some right-hand side, which is shown in Algorithm~\ref{algo:double_loop} (cf.~\cite{Nocedal2006Numerical} for a description of the L-BFGS method in a finite-dimensional setting).

\begin{Remark}
	Our numerical experiments presented in Section~\ref{sec:numerics} showed that the search direction computed with \eqref{eq:search_direction_bfgs} may sometimes not yield descent of the cost functional. Therefore, we employ a restarted version of the BFGS method, which replaces the search direction with $p_k = - g_k$ if the line search procedure in line~\ref{algo_bfgs:line_search} of Algorithm~\ref{algo:lbfgs} does not converge.
\end{Remark}

\begin{Remark}
	The idea presented above could also be applied analogously to derive nonlinear conjugate gradient (NCG) methods for topology optimization. As a thorough investigation of the several popular NCG methods is beyond the scope of this paper, we do not consider these methods in the following (cf.~\cite{Blauth2021Nonlinear} for a discussion of NCG methods in the context of shape optimization). However, the NCG methods are already implemented and ready for use in our software package cashocs \cite{Blauth2021cashocs, Blauth2023Version}. An investigation of such NCG methods for topology optimization is planned for future research.
\end{Remark}

\begin{algorithm2e}[!t]
	\KwIn{Initial level-set function $\psi_0$, initial step size $\lambda_0 = 1$, initial inverse Hessian approximation $H_0 = \mathrm{Id}_{L^2(D)}$, stopping tolerance $\tau > 0$, maximum number of iterations $k_{\mathrm{max}} \in \mathbb{N}$, memory size $m\geq 0$}
	\For{$k=0,1,2,\dots, k_\mathrm{max}$}{
		Compute the generalized topological derivative $\mathcal{D}J(\Omega_k)$\\
		Set $g_k = -P_{\psi_k^\perp}(\mathcal{D}J(\Omega_k)) = - \left( \mathcal{D}J(\Omega_k) - \frac{\left( \mathcal{D}J(\Omega_k), \psi_k \right)}{\norm{\psi_k}{L^2(\Dsf)}^2} \psi_k \right)$\\
		Compute the search direction $p_k = -H_k g_k$ via Algorithm~\ref{algo:double_loop}\\
		Compute $\theta_k = \arccos\left( \frac{(\psi_k, \mathcal{D}J(\Omega_k))}{\norm{\psi_k}{L^2(\Dsf)} \norm{\mathcal{D}J(\Omega_k)}{L^2(\Dsf)}} \right)$\\
		\If{$\theta_k < \tau$}{
			Stop with approximate minimizer $\Omega_k$\\
		}
		Set $\lambda_k = 1$\\
		\While{$J\left(\psi_k + \lambda_k p_k\right) > J(\psi_k)$ \label{algo_bfgs:line_search}}{
			Decrease the step size $\lambda_k = \frac{\lambda_k}{2}$\\
		}
		Update the level-set function: $\psi_{k+1} = \psi_k + \lambda_k p_k$\\
		Store $y_k = g_{k+1} - g_k$, $s_k = \lambda_k p_k$, and $\rho_k = \frac{1}{(y_k, s_k)}$
	}
	\caption{Limited-memory BFGS method for topology optimization.}
	\label{algo:lbfgs}
\end{algorithm2e}

\begin{algorithm2e}[!t]
	\KwIn{Right-hand side $b$ and stored elements $y_k, s_k, \rho_k$ for $i=k-m,\dots, k-1$}
	$q = b$ \\
	\For{$i=k-1, k-2, \dots, k-m$}{
		Store $\alpha_i = \rho_i (s_i, q)$\\
		$q = q - \alpha_i y_i$ \\
	}
	$H_0^k =  \frac{(s_{k-1}, y_{k-1})}{(y_{k-1}, y_{k-1})} \text{Id}_{L^2(\Dsf)}$ \\
	$r = H_0^k q$ \\
	\For{$i=k-m, k-m+1, \dots, k-1$}{
		$\beta_i = \rho_i (y_i, r)$\\
		$r = r + (\alpha_i - \beta_i) s_i $ \\
	}
	\KwRet{$H_k b = r$}
	\caption{Limited-memory double loop for computing $H_k b$.}
	\label{algo:double_loop}
\end{algorithm2e}

\section{Numerical Investigation of Quasi-Newton Methods for Topology Optimization}
\label{sec:numerics}

In this section, we consider the practical performance of the proposed quasi-Newton methods from Section~\ref{sec:quasi_newton}. We consider four problem classes: Inverse topology optimization problems constrained by linear and semilinear Poisson problems, compliance minimization in linear elasticity, and the topological design optimization in Navier-Stokes flow. We solve each of these problems with the four solution algorithms presented in this paper, where we use the following notation for the sake of simplicity. Algorithm~\ref{algo:euler}, originally introduced in \cite{Amstutz2006new}, is called the sphere combination method, Algorithm~\ref{algo:convex_combination} is called convex combination method, Algorithm~\ref{algo:gradient_descent} is called gradient descent method and, finally, Algorithm~\ref{algo:lbfgs} is called (limited memory) BFGS method. We remark that, for all numerical test cases considered in this paper, we choose a memory size of $m=5$ for the limited memory BFGS methods.

Note that we have implemented all of the optimization algorithms for topology optimization considered in Section~\ref{sec:preliminaries} and the novel gradient descent and BFGS methods from Section~\ref{sec:quasi_newton} in our open-source software package cashocs \cite{Blauth2021cashocs, Blauth2023Version}, which is a software for solving arbitrary PDE constrained shape optimization and optimal control problems. Our software is based on the finite element software FEniCS \cite{Alnes2015FEniCS,Logg2012Automated} and derives the necessary adjoint systems for the optimization with the help of automatic differentiation. Therefore, for the discretization of the PDE constraints, the finite element method is naturally employed. Moreover, the source code for our numerical experiments is available freely on GitHub \cite{Blauth2023Software}.

\subsection{Linear Poisson Problem}
\label{ssec:numerics_model}

In this section, we investigate a topology optimization problem constrained by a linear Poisson problem. Let $\Dsf \subset \R^d$ with $d\in \mathbb{N}_{>0}$ be an open and bounded domain with boundary $\partial \Dsf$ and let $\Omega \subset \Dsf$ be an open subset. We denote by $\Gamma = \partial\Omega \setminus \partial \Dsf$ the interior boundary of $\Omega$ in $\Dsf$ and by $\Omega^c = \Dsf \setminus \closure{\Omega}$ the complement of $\Omega$ in $\Dsf$. We consider the following problem
\begin{equation}
\label{eq:state_system_poisson}
\begin{aligned}
&\min_{\Omega} J(\Omega, u) = \frac{1}{2} \int_{\Dsf} \left( u - u_\mathrm{des} \right)^2 \dmeas{x} \\
&\text{s.t.} \quad \begin{alignedat}[t]{2}
-\Delta u + \alpha_\Omega u &= f_\Omega \quad &&\text{ in } \Dsf,\\
u &= 0 \quad &&\text{ on } \partial \Dsf,
\end{alignedat}
\end{aligned}
\end{equation}
where $\alpha_\Omega(x) = \chi_\Omega(x) \alpha_\mathrm{in} + \chi_{\Omega^c}(x) \alpha_\mathrm{out}$ and $f_\Omega(x) = \chi_\Omega(x) f_\mathrm{in} + \chi_{\Omega^c}(x) f_\mathrm{out}$ with $\alpha_\mathrm{in}, \alpha_\mathrm{out} > 0$ as well as $f_\mathrm{in}, f_\mathrm{out} \in \mathbb{R}$ are constant in $\Omega$ and $\Omega^c$.
In our setting, $u_\mathrm{des}$ is given as the solution of the PDE constraint on a desired domain $\Omega_\mathrm{des}$. Hence, the above problem can be interpreted as an inverse problem of identifying the unknown domain $\Omega_\mathrm{des}$ using the measurement $u_\mathrm{des}$. The generalized topological derivative for this problem is derived, e.g., in \cite{Amstutz2022introduction} and is given by
\begin{equation*}
\begin{aligned}
\mathcal{D}J(\Omega)(x) = (\alpha_\mathrm{out} - \alpha_\mathrm{in}) u(x) p(x) - (f_\mathrm{out} - f_\mathrm{in}) p(x) \quad \text{ for all } x\in \Dsf \setminus \Gamma,
\end{aligned}
\end{equation*}
where $u$ solves the state equation \eqref{eq:state_system_poisson} and $p$ solves the following adjoint equation
\begin{equation*}
\begin{alignedat}{2}
-\Delta p + \alpha_\Omega p &= - (u-u_\mathrm{des}) \quad &&\text{ in } \Dsf,\\
p &= 0 \quad &&\text{ on } \partial \Dsf.
\end{alignedat}
\end{equation*}

\begin{Remark}
	Usually, the term \qe{inverse problem} is used to denote a identification problem which aims to reconstruct some unknown domain $\Omega_\mathrm{des} \subset \Dsf$ using measurements obtained on (parts of) the boundary of $\Dsf$, see, e.g., \cite{a_CALANO_2015a, Hintermueller2008Electrical, a_LAHIFRSC_2013a, Beretta2017reconstruction}. However, for our model problems \eqref{eq:state_system_poisson} and \eqref{eq:semilinear_poisson}, we use (artificial) measurements in the entire domain $\Dsf$. Still, we refer to these problems as inverse problems for the sake of simplicity.
\end{Remark}

\begin{figure}[!b]
	\centering
	\begin{subfigure}{0.333\textwidth}
		\centering
		\includegraphics[width=\textwidth]{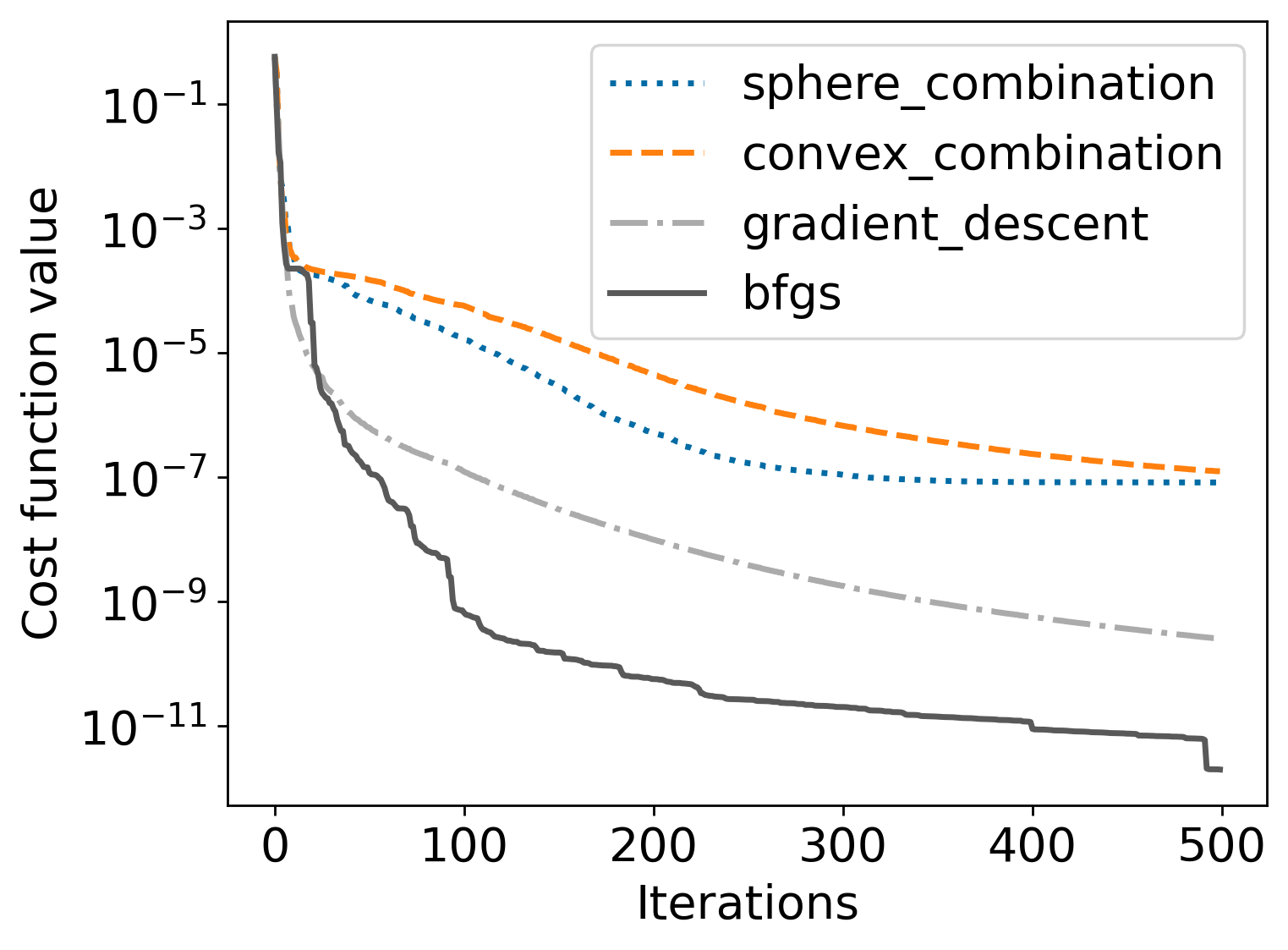}
		\caption{Cost functional.}
	\end{subfigure}%
	\begin{subfigure}{0.333\textwidth}
		\centering
		\includegraphics[width=\textwidth]{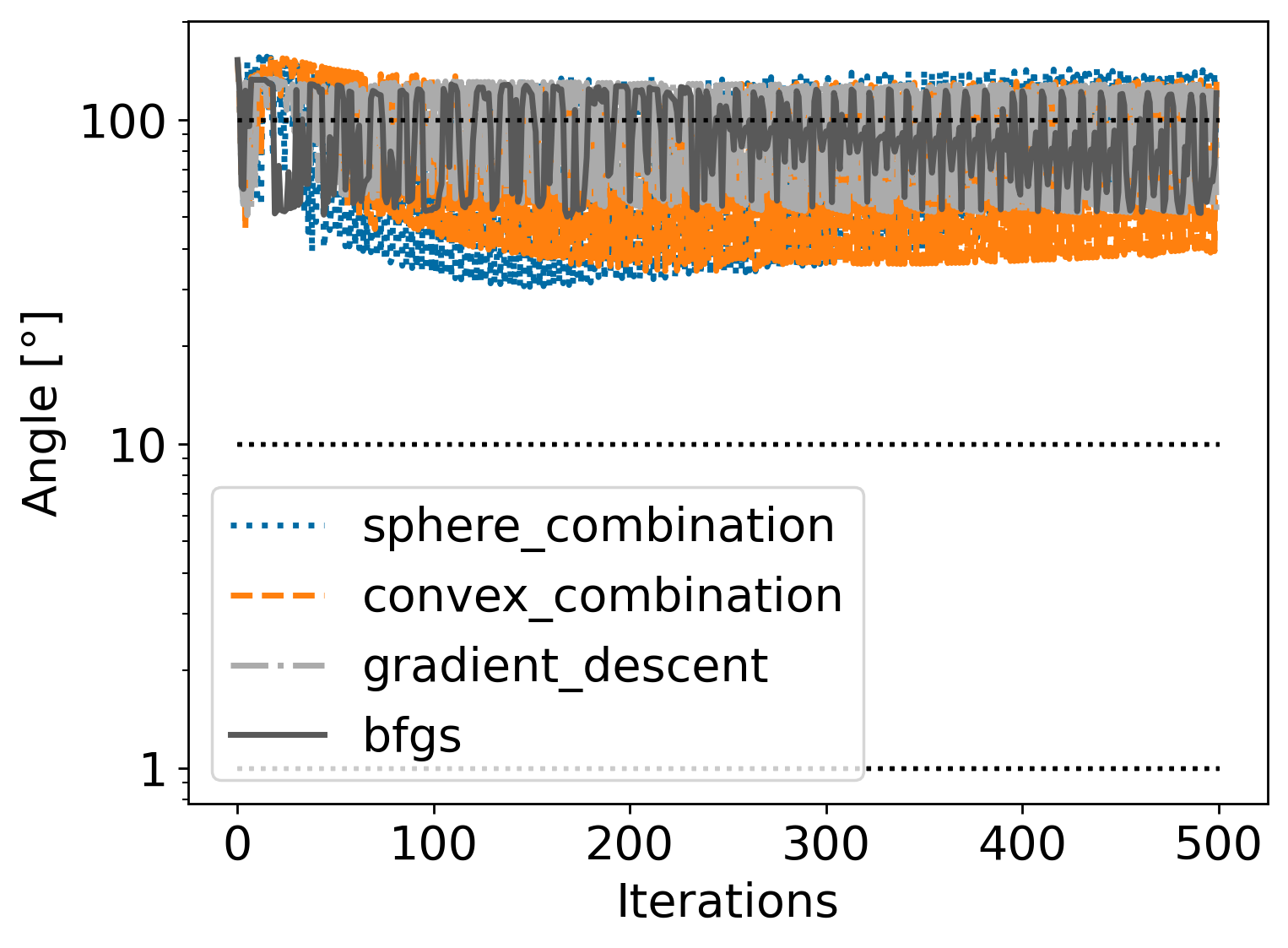}
		\caption{Angle.}
	\end{subfigure}%
	\begin{subfigure}{0.333\textwidth}
		\centering
		\includegraphics[width=\textwidth]{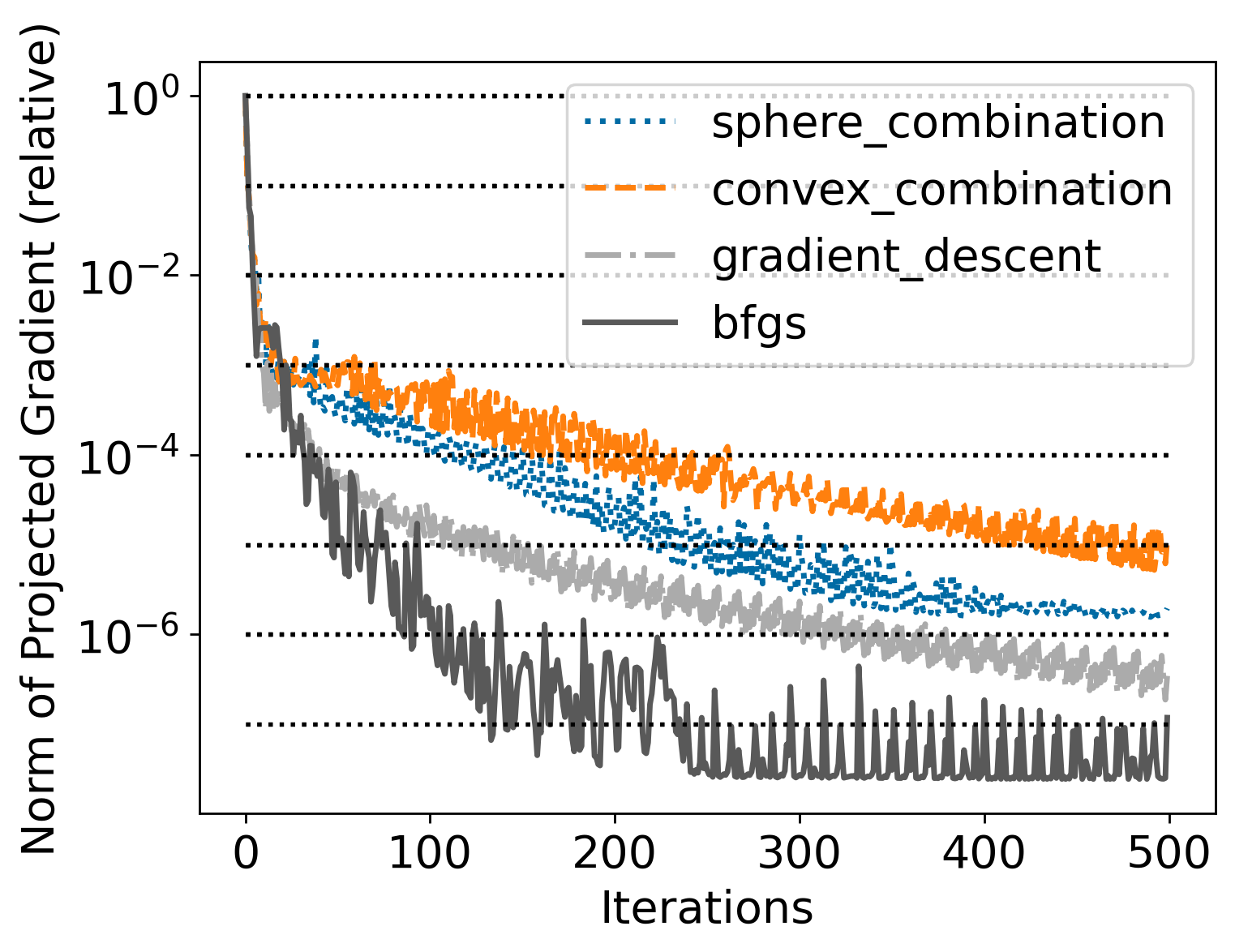}
		\caption{Norm of the projected topological derivative.}
	\end{subfigure}
	\caption{Evolution of the optimization for the linear Poisson problem \eqref{eq:state_system_poisson}.}
	\label{fig:poisson_nonprincipal}
\end{figure}

In the following, we solve this problem with the four optimization algorithms described in Sections~\ref{sec:preliminaries} and~\ref{sec:quasi_newton}, where we consider a hold-all domain of $\Dsf = (-2, 2)^2$. We discretize the geometry using four different mesh sizes in order to investigate the dependence of the algorithms on the discretization. The considered meshes are generated by creating a uniform quadrangular grid with $32 \times 32$, $48\times 48$, $64\times 64$, and $96\times 96$ squares, where each square is in turn subdivided into four triangles, so that the meshes consist of \num{2113}~nodes and \num{4096}~triangles ($32\times 32$), \num{4705}~nodes and \num{9216}~triangles ($48\times 48$), \num{8321}~nodes and \num{16384}~triangles ($64\times 64$), and \num{18625}~nodes and \num{36864}~triangles ($96\times 96$), respectively. Moreover, we discretize both the state and adjoint variables with linear Lagrange finite elements. For the parameters, we use $\alpha_\mathrm{in} = f_\mathrm{in} = 10$ and $\alpha_\mathrm{out} = f_\mathrm{out} = 1$. The sought design $\Omega_\mathrm{des}$ is chosen as the one corresponding to a clover shape, which can be seen in the right-most column of Table~\ref{tab:comparison_poisson}. 

The results of the optimization can be seen in Figure~\ref{fig:poisson_nonprincipal}, where we show the evolution of the cost functional, the angle criterion, and the norm of the projected gradient (cf.~\eqref{eq:norm_projected_derivative}) over the optimization for the finest discretization. Here, we observe that our proposed gradient descent and BFGS methods perform significantly better than the established methods. We observe that the convergence behavior, based on the cost functional and the norm of the projected gradient, is best for the BFGS method, followed by the gradient descent and sphere combination methods, and that the convex combination algorithm performs worst. Particularly, the BFGS method reaches stationarity in the cost functional and norm of the projected gradient after only about 125 iterations, whereas all remaining methods continue to decrease the respective measures until the final iteration. 

\begin{table}[!t]
	\newcommand{\sizebox}{0.215\textwidth}
	\newcolumntype{C}{>{\centering\arraybackslash}m{\sizebox}}
	\setlength{\tabcolsep}{0pt}
	\caption{Evolution of the geometries for the linear Poisson problem \eqref{eq:state_system_poisson}.}
	\label{tab:comparison_poisson}
	\begin{tabular}{m{0.1\textwidth} @{\hskip 0.5em} C C C @{\hskip 1em} C}
		\toprule
		Algorithm & 50 iterations & 100 iterations & 500 iterations & Reference \\ 
		\midrule
		Sphere combination & \includegraphics[width=\sizebox]{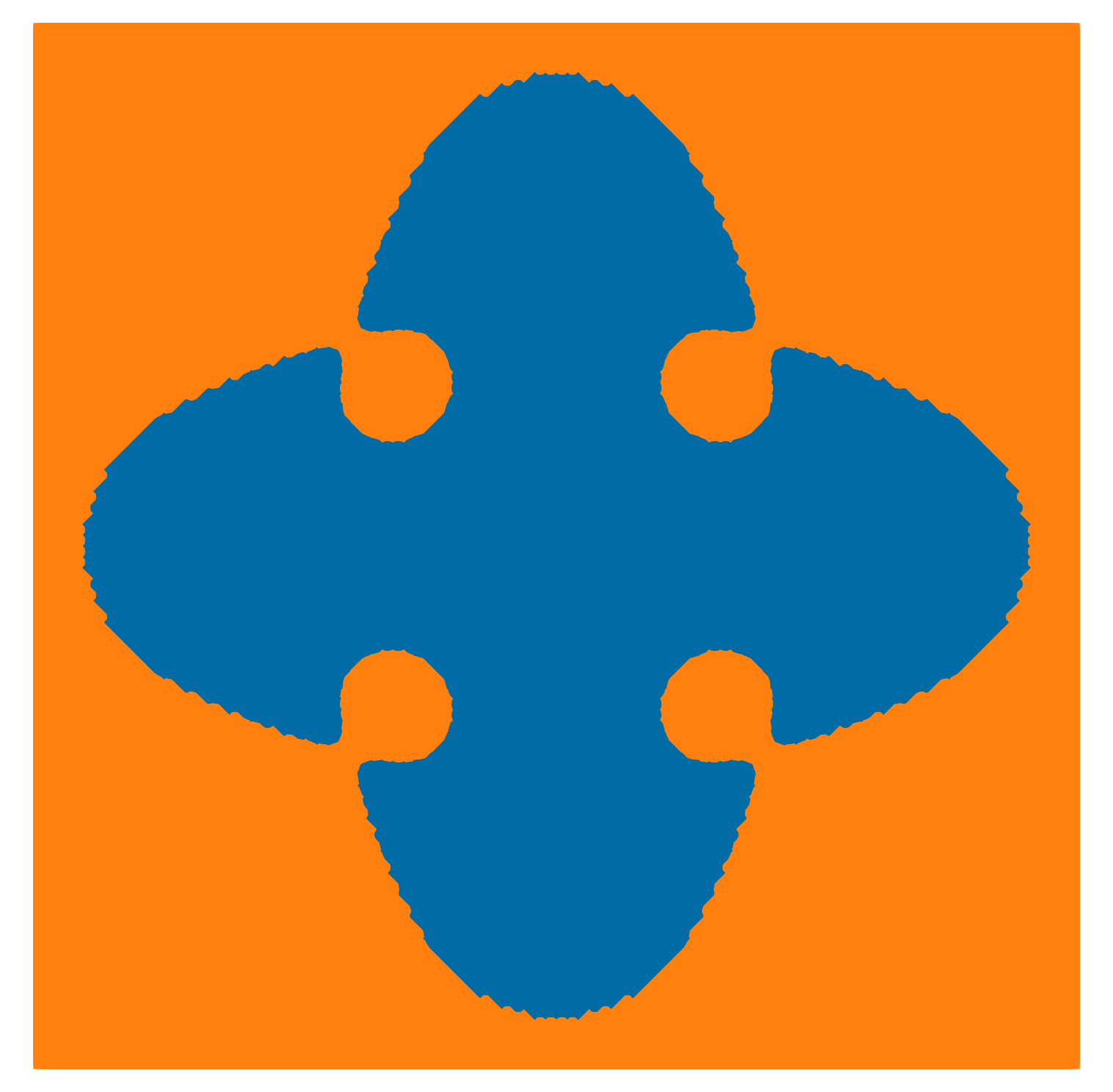} & \includegraphics[width=\sizebox]{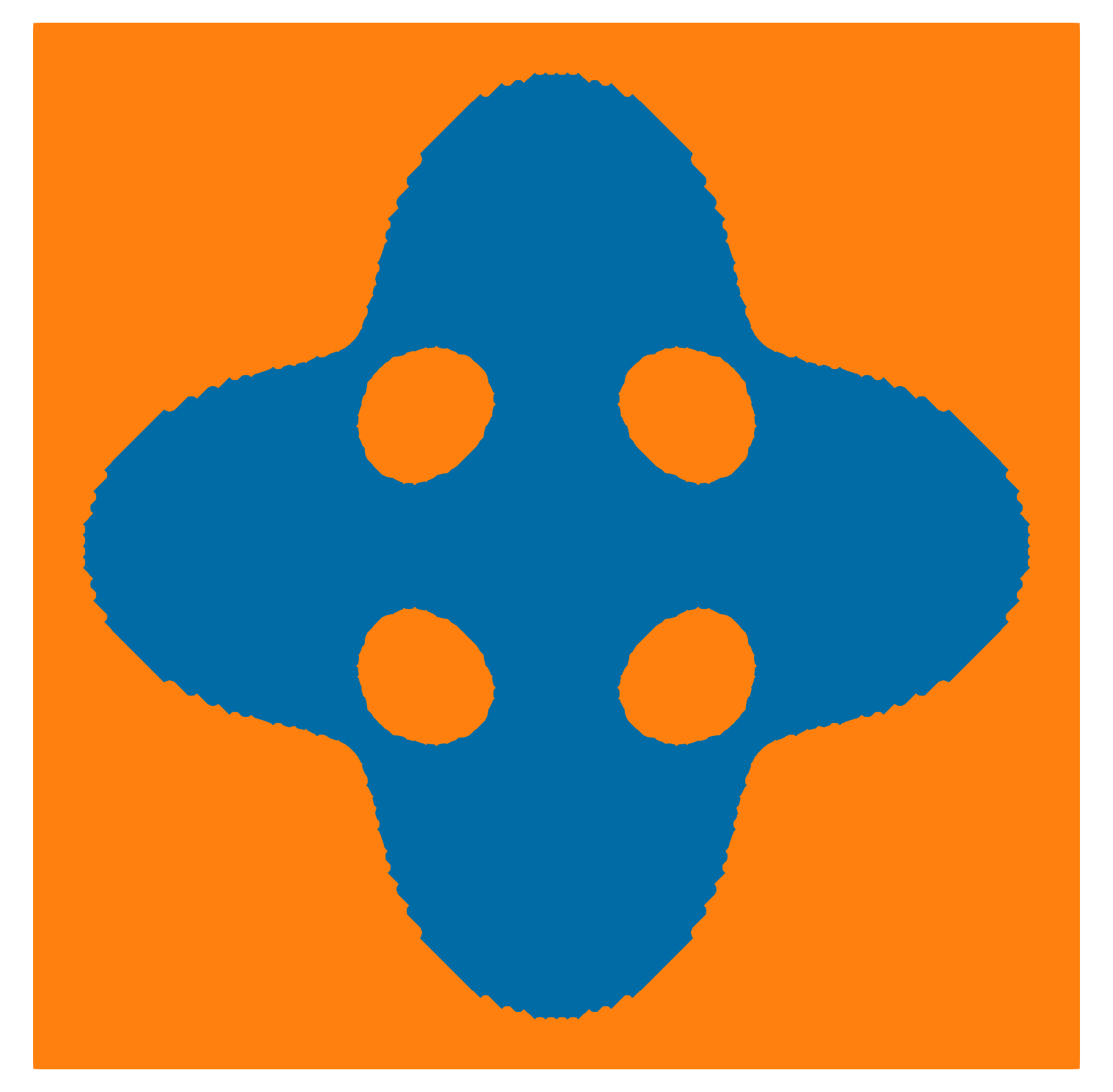} & \includegraphics[width=\sizebox]{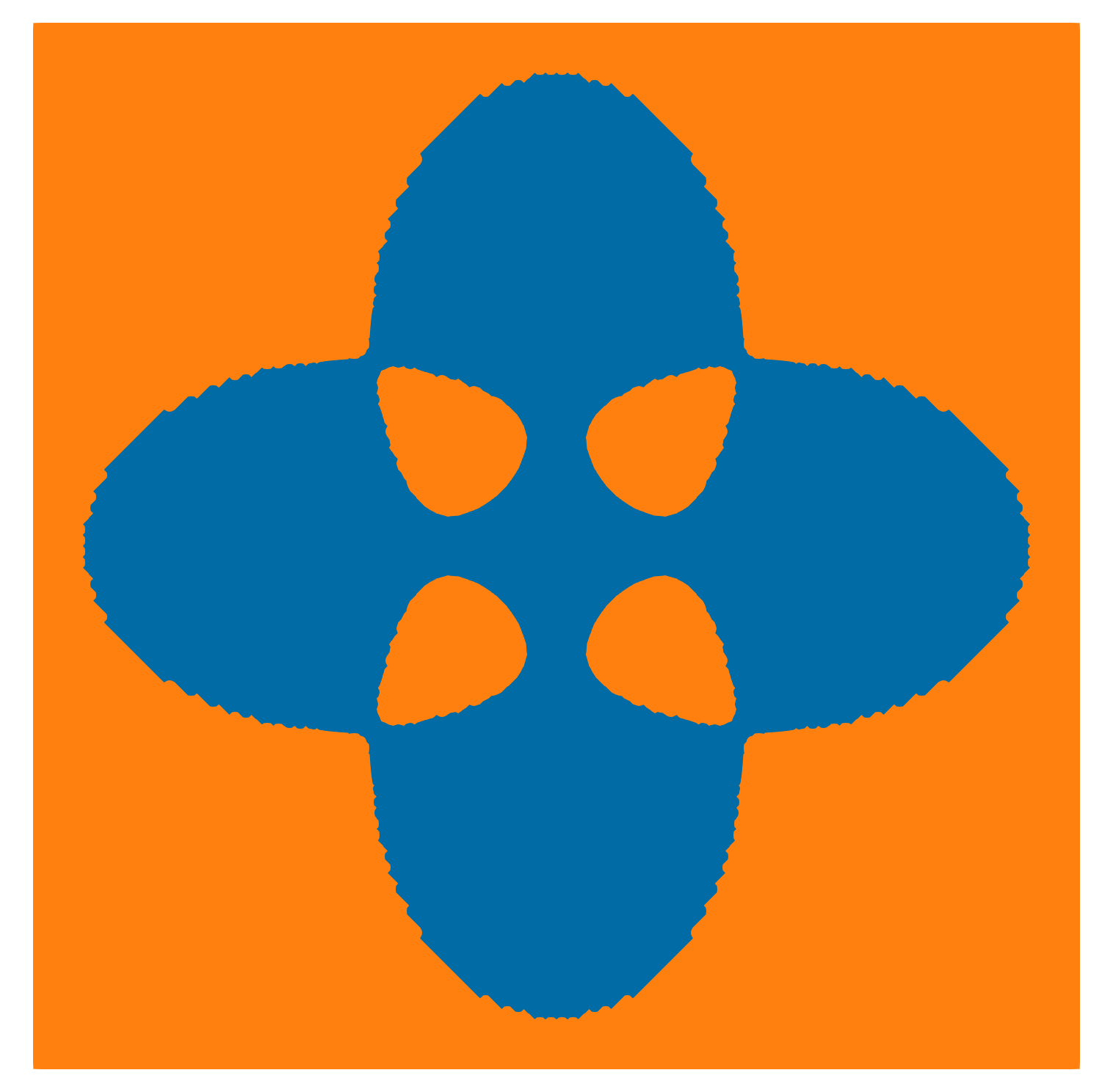} & \includegraphics[width=\sizebox]{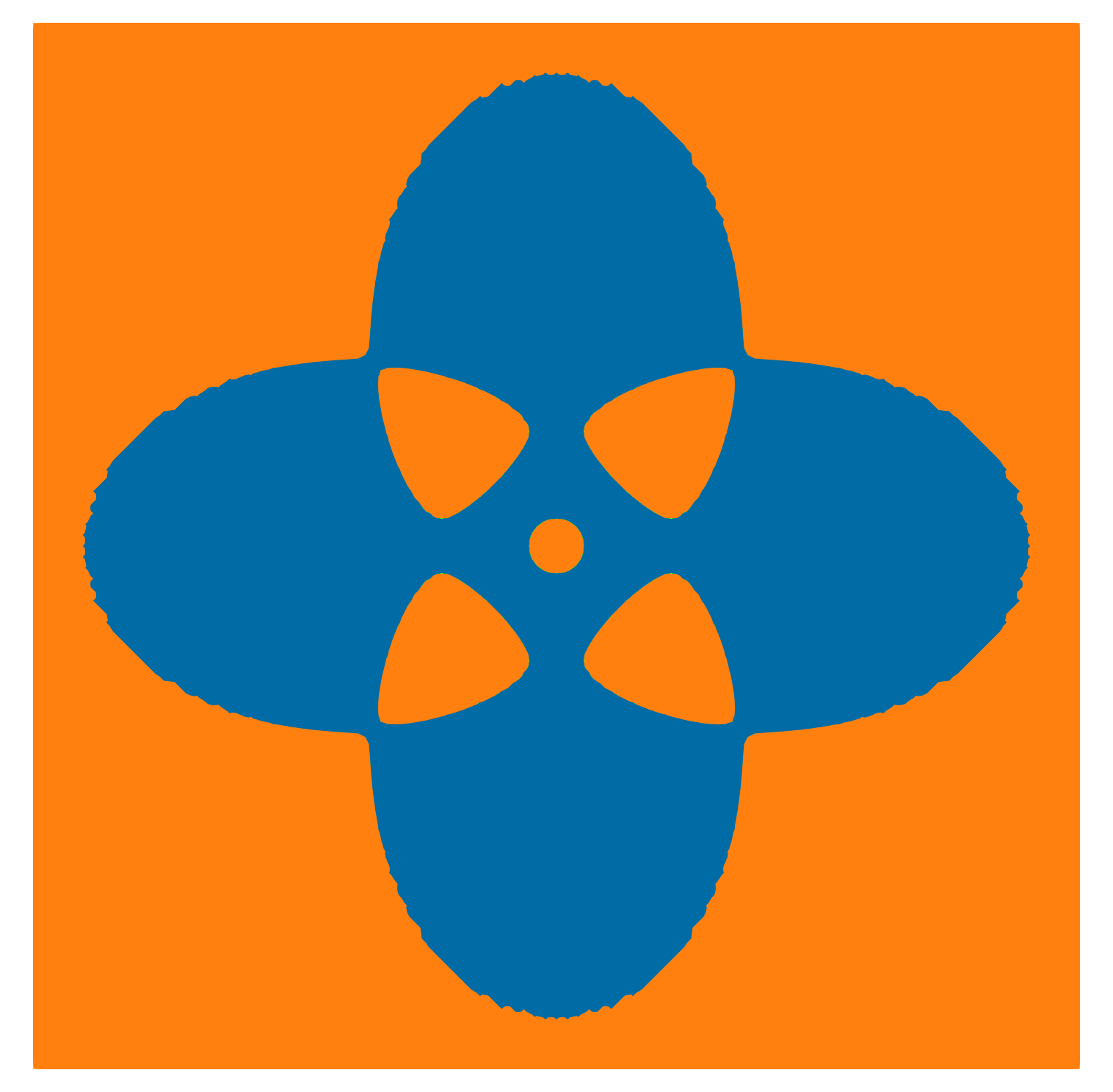} \\ 
		Convex combination & \includegraphics[width=\sizebox]{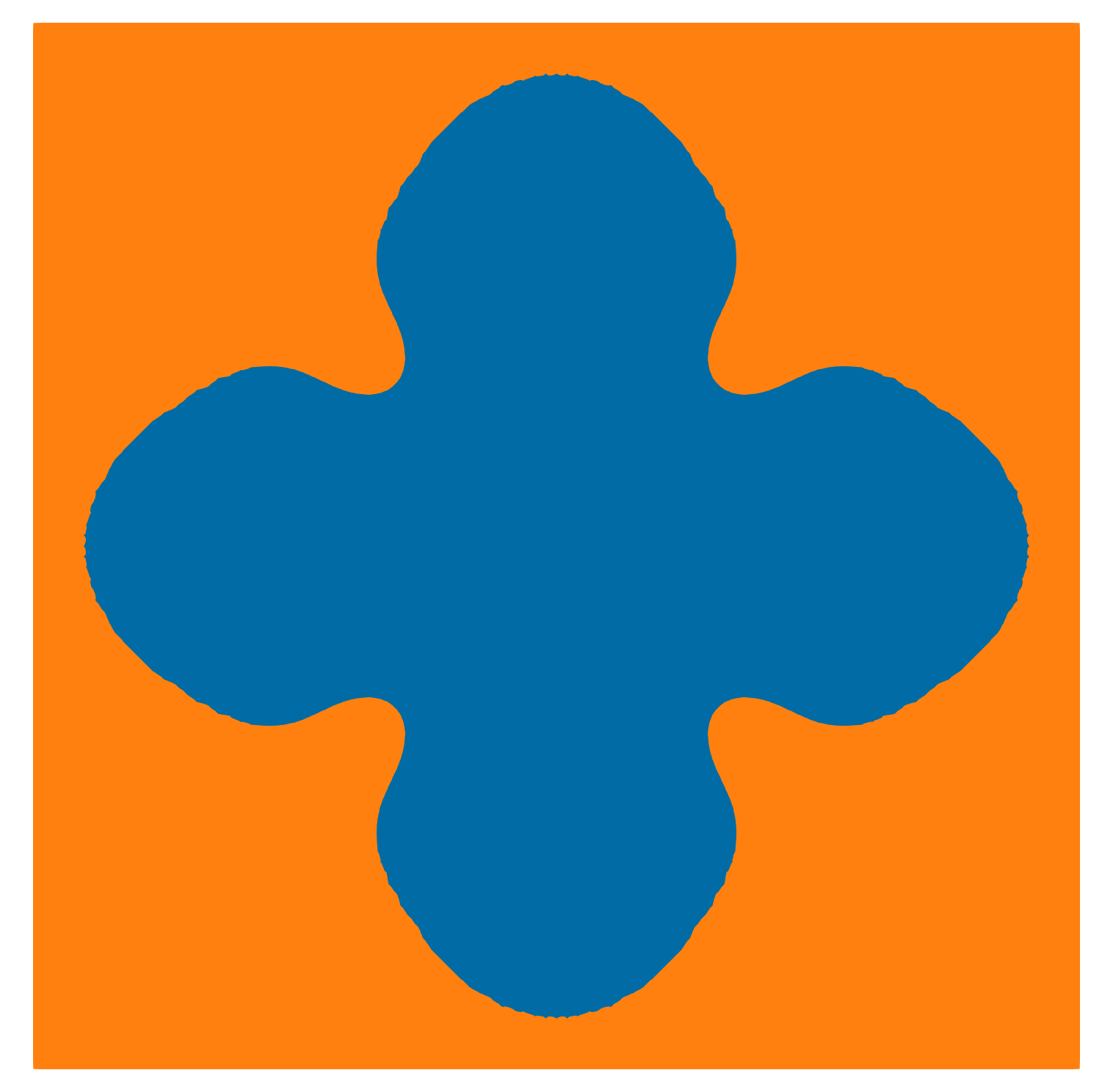} & \includegraphics[width=\sizebox]{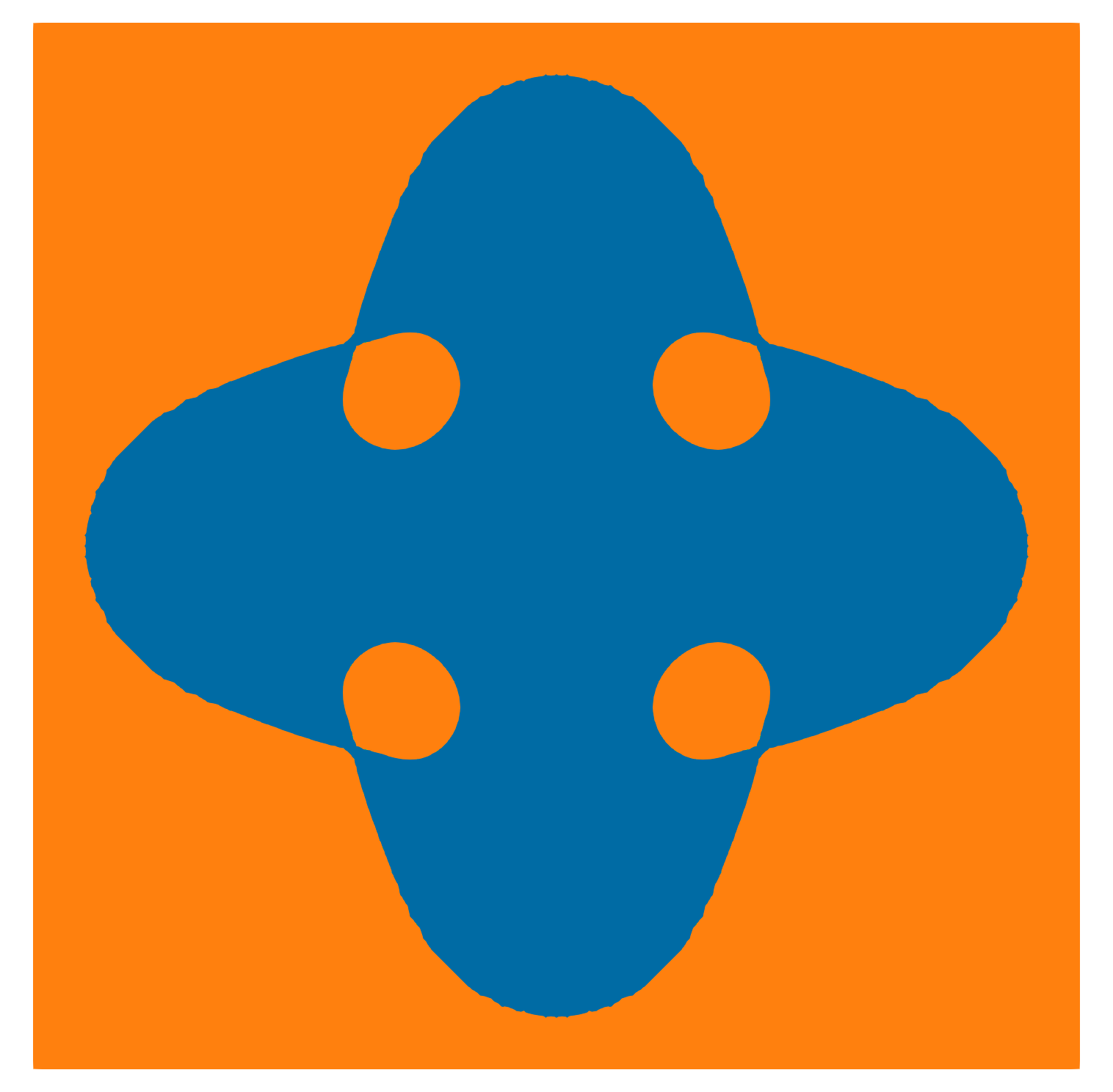} & \includegraphics[width=\sizebox]{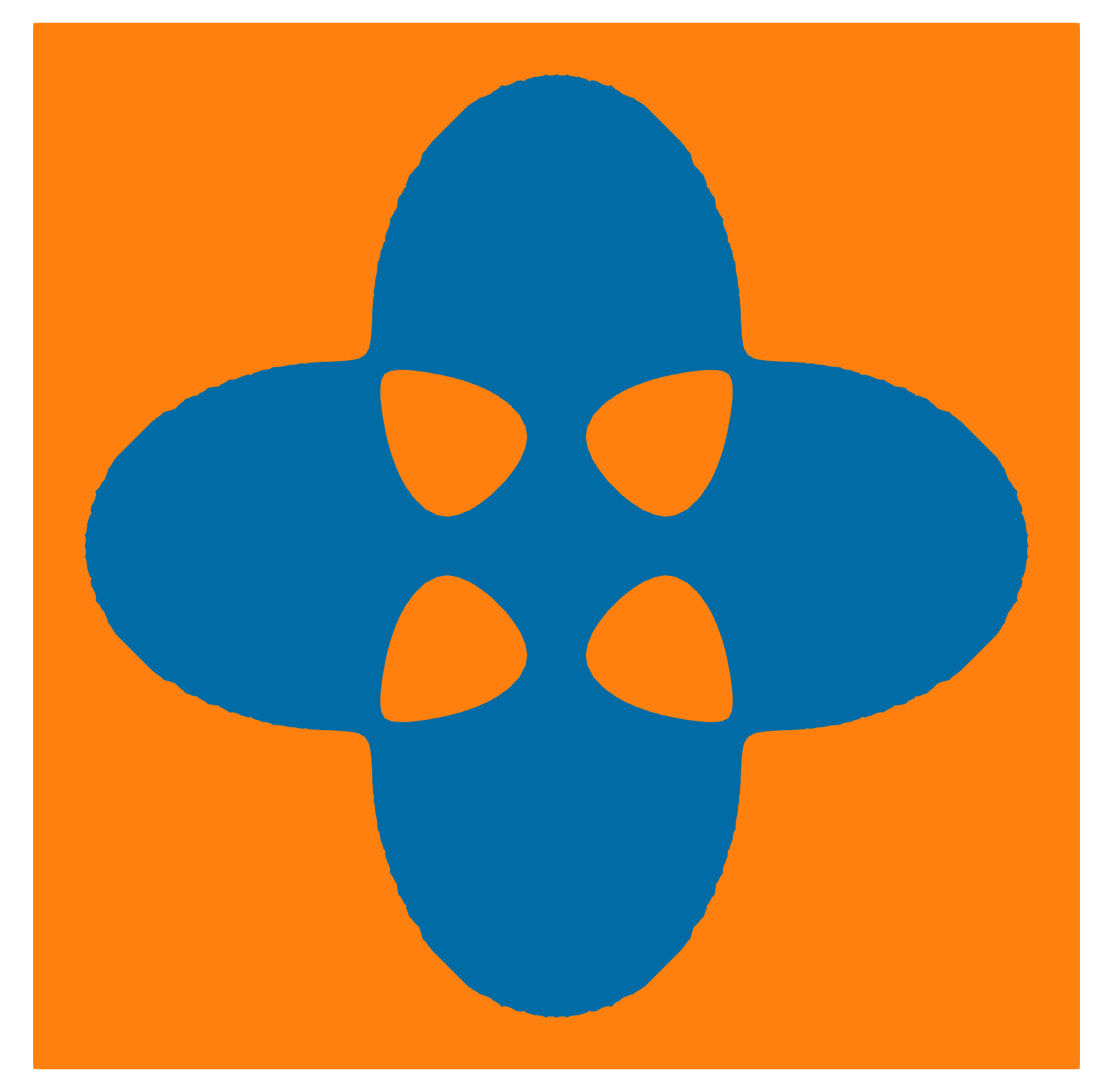} & \includegraphics[width=\sizebox]{clover} \\ 
		Gradient descent & \includegraphics[width=\sizebox]{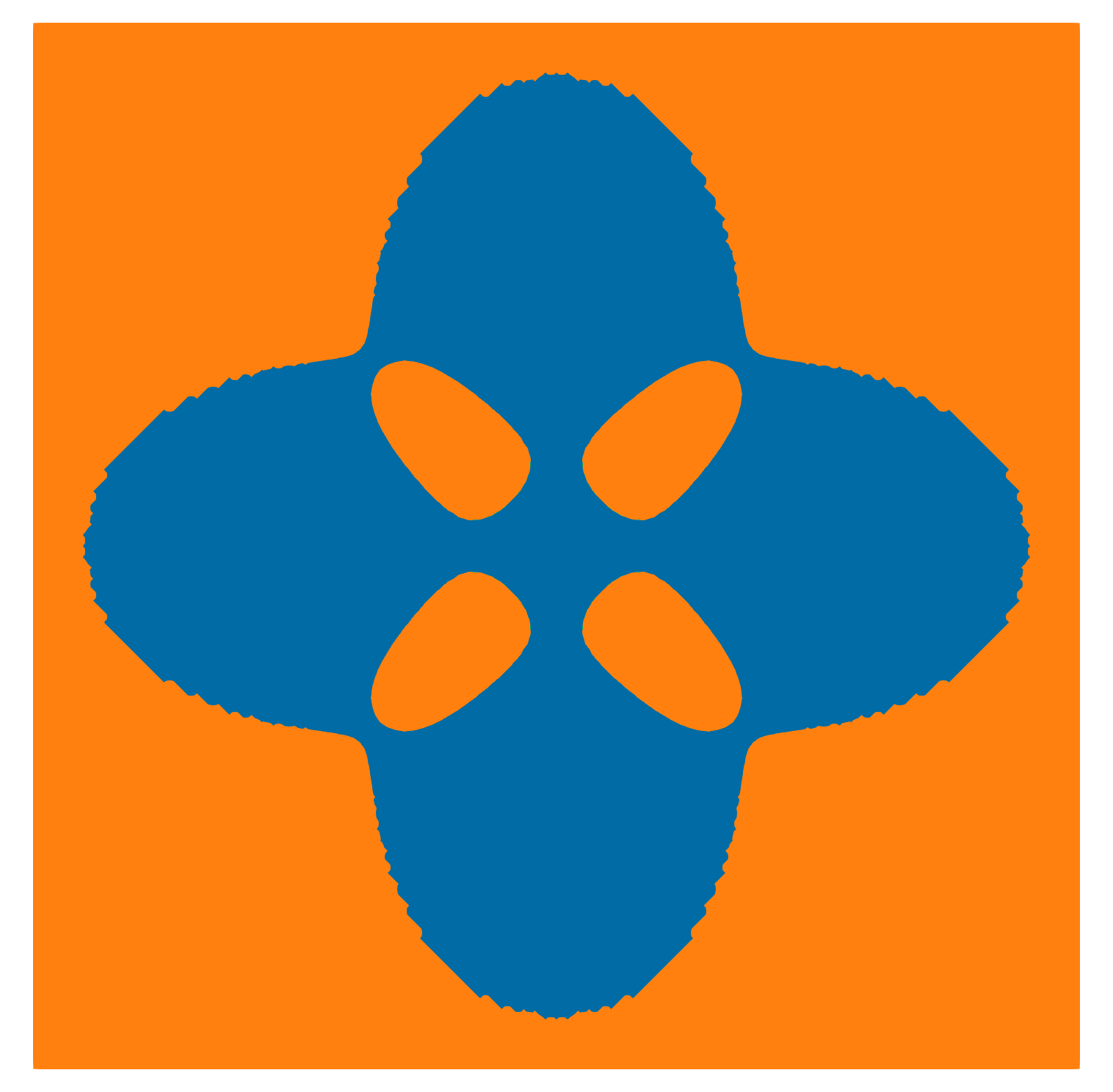} & \includegraphics[width=\sizebox]{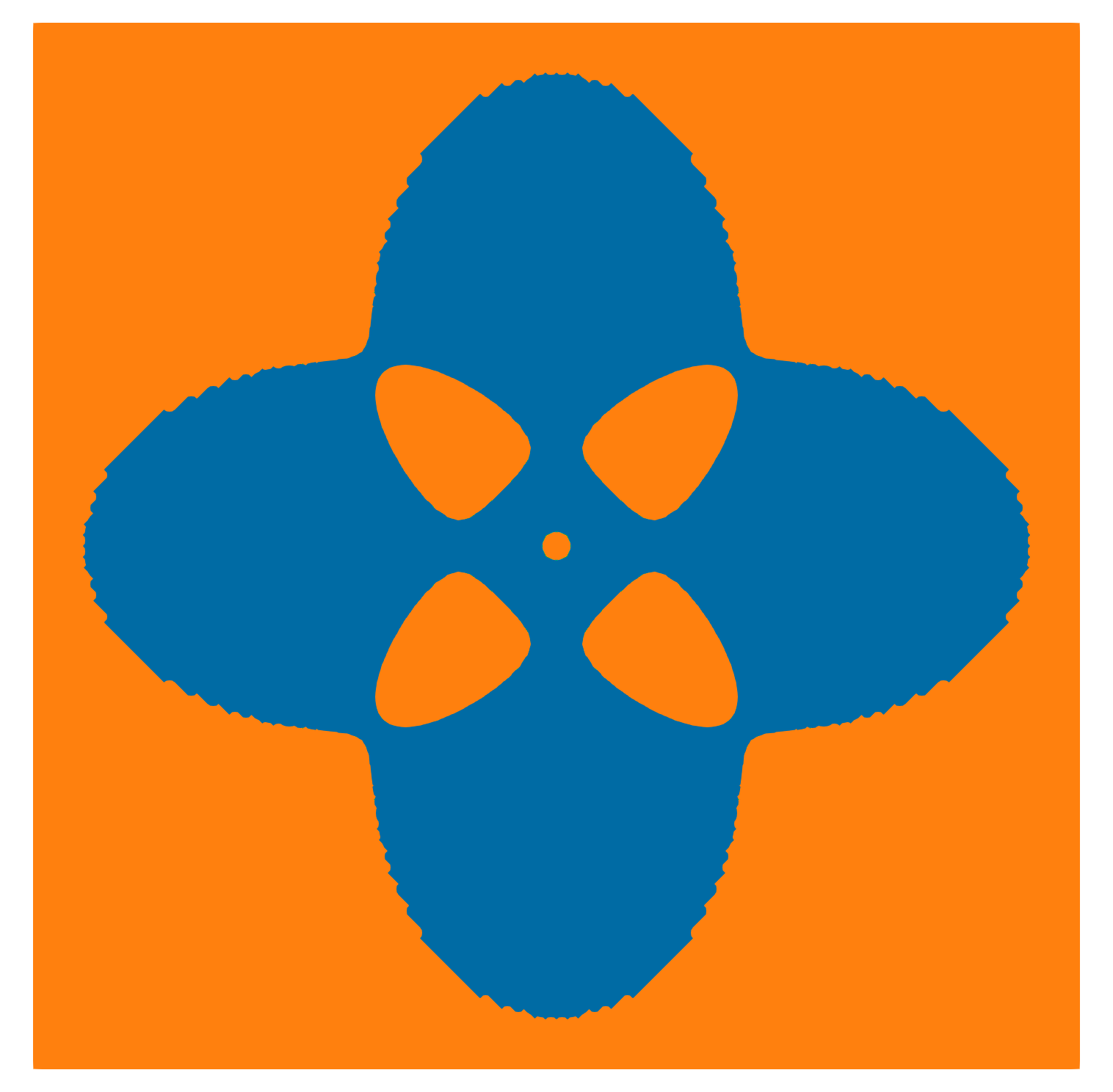} & \includegraphics[width=\sizebox]{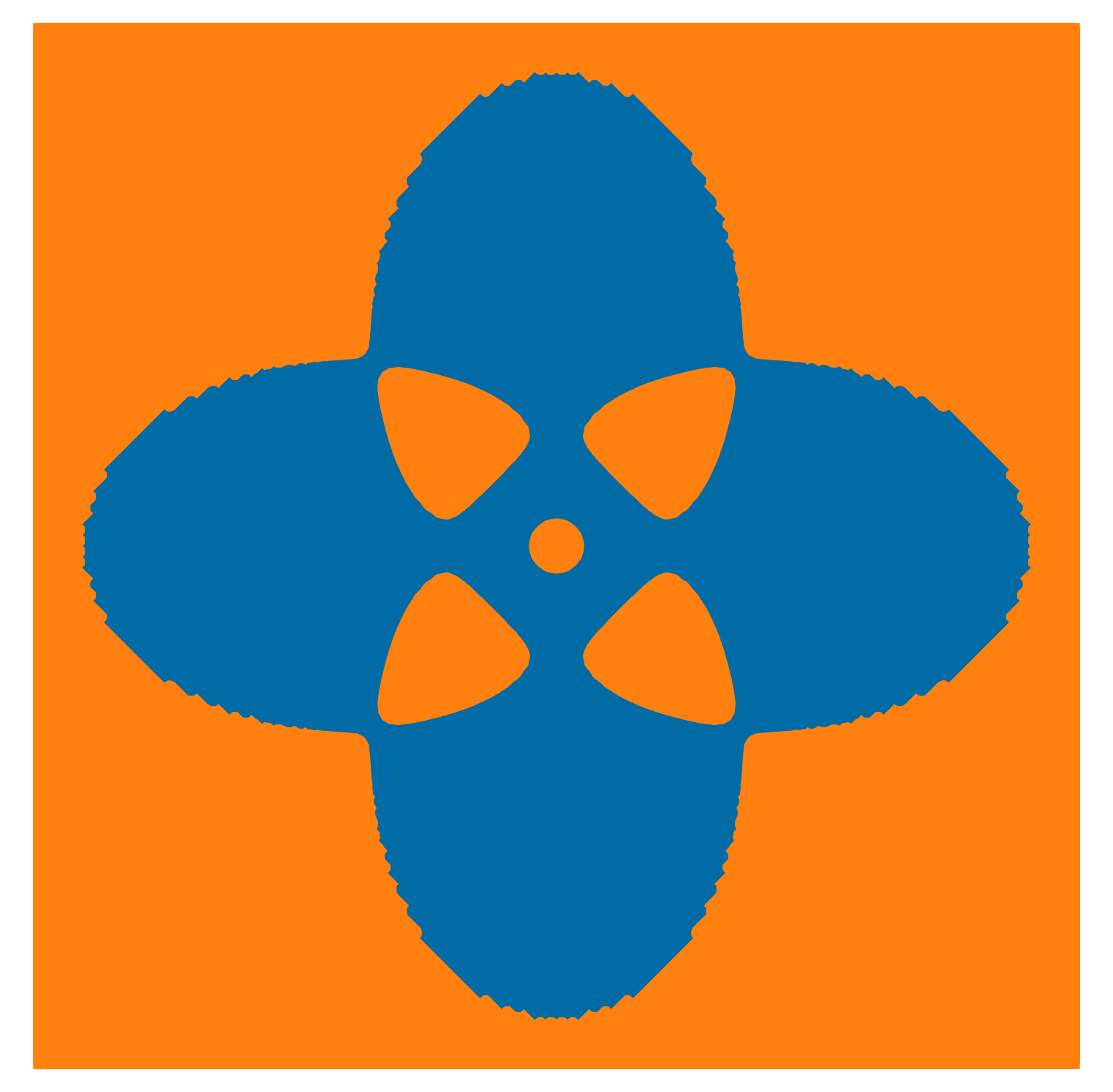} & \includegraphics[width=\sizebox]{clover} \\ 
		BFGS & \includegraphics[width=\sizebox]{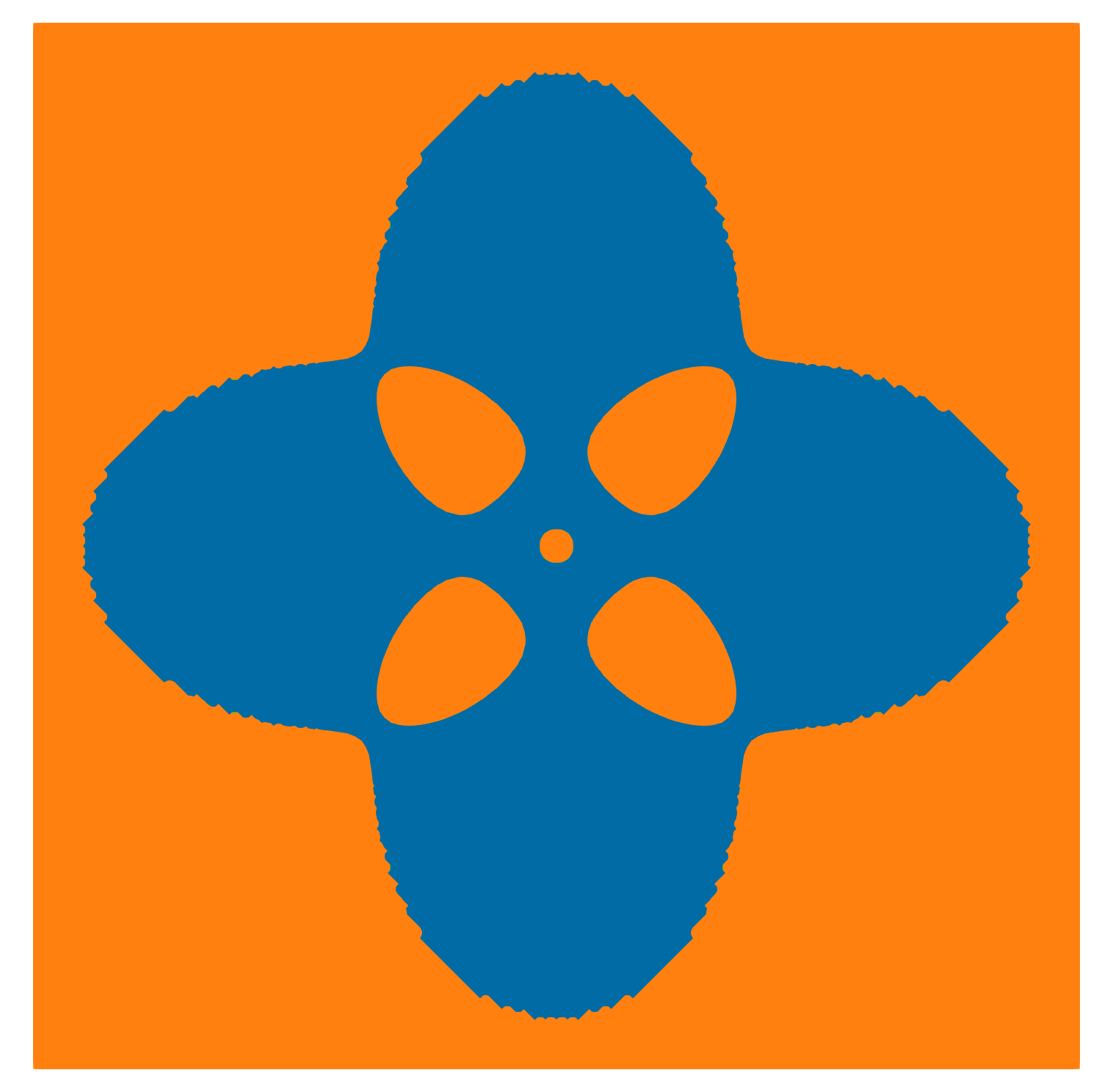} & \includegraphics[width=\sizebox]{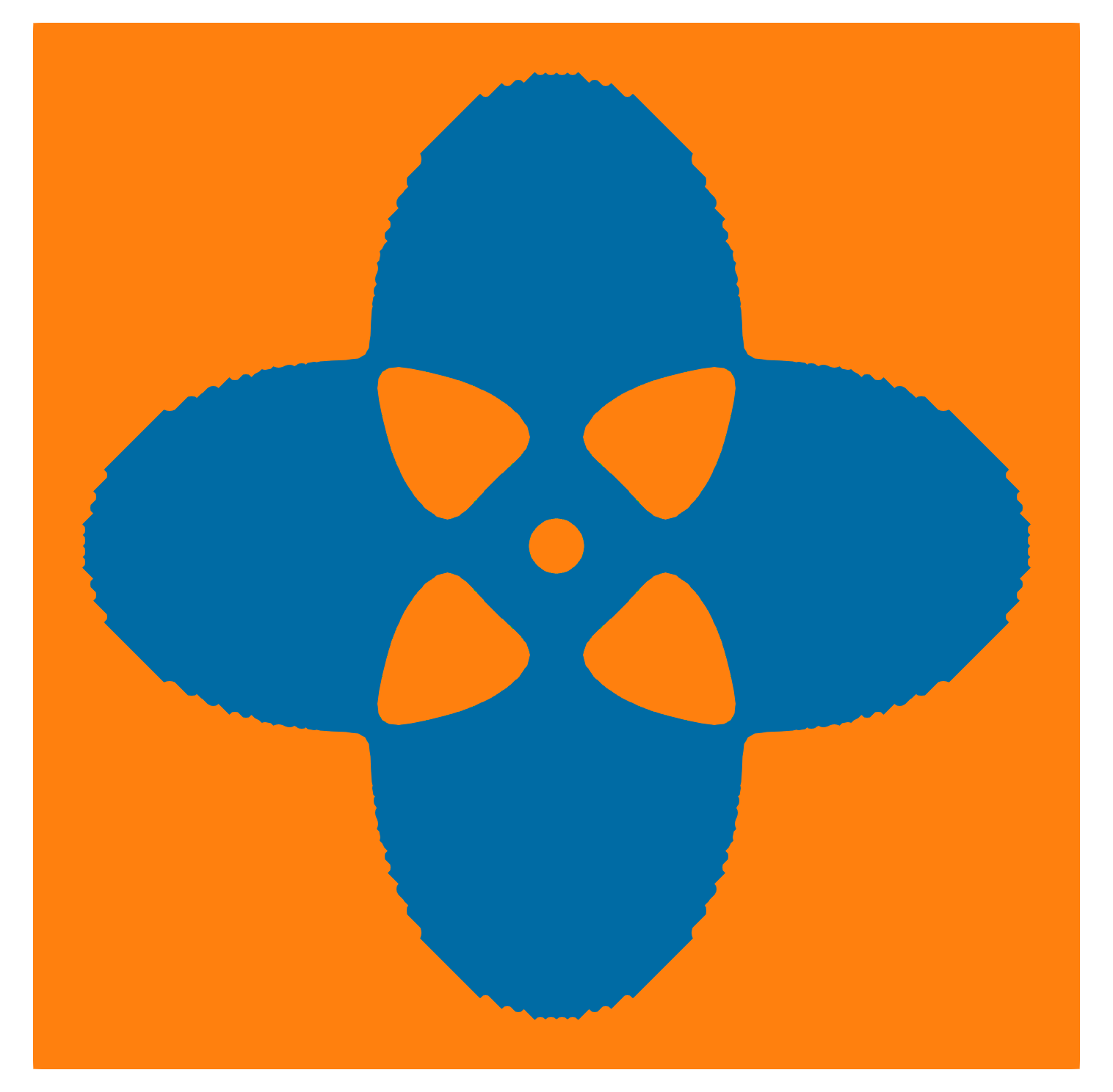} & \includegraphics[width=\sizebox]{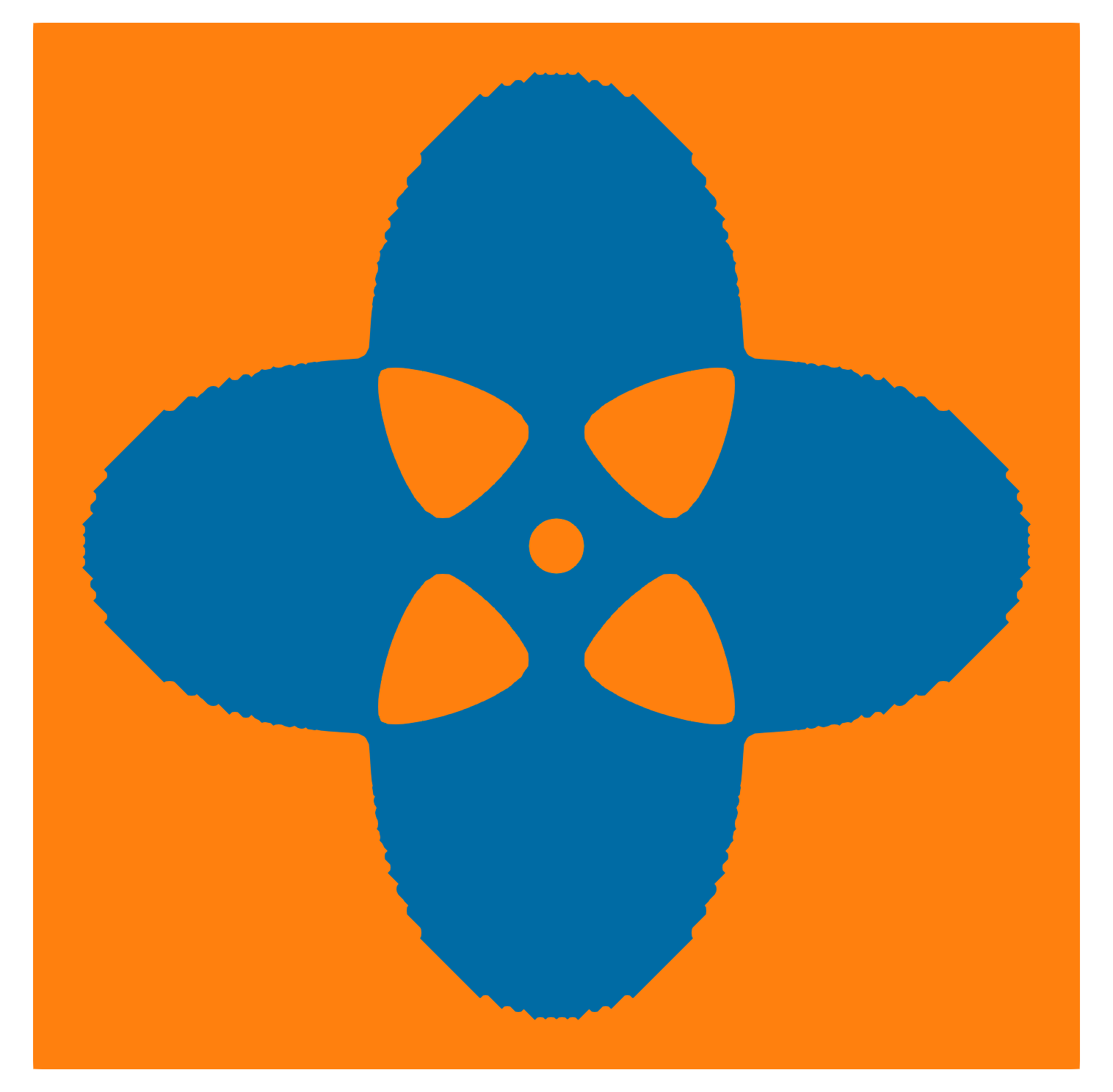} & \includegraphics[width=\sizebox]{clover} \\ 
		\bottomrule 
	\end{tabular}
\end{table}

However, none of the methods converged based on the angle criterion after \num{500}~iterations. On the contrary, the angle between level-set function and generalized topological derivate remains bounded from below by about \SI{30}{\degree}, so that none of the methods can be considered as converged by this criterion. However, if we take a look at the evolution of the relative norm of the projected gradient, we observe a steep decrease for all methods. The convex and sphere combination methods are able to decrease the norm of the projected gradient by about five and six orders of magnitude, respectively. The gradient descent and BFGS methods are even more efficient and decrease the norm of the projected gradient by about seven and eight orders of magnitude, respectively. Therefore, based on the second convergence criterion, the algorithms can be considered converged. Moreover, for this example, it seems like the angle criterion is too strict as a convergence criterion, whereas the norm of the projected gradient seems to be better suited. 

\begin{figure}[!b]
	\centering
	\begin{subfigure}{\textwidth}
		\centering
		\includegraphics[width=\textwidth, trim=35cm 0cm 35cm 0cm]{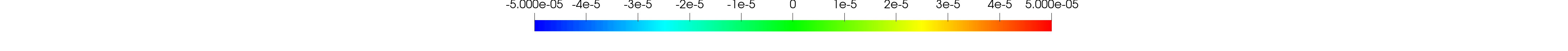}
	\end{subfigure}
	\begin{subfigure}{0.5\textwidth}
		\centering
		\includegraphics[width=\textwidth, trim=39cm 0cm 39cm 0cm, clip]{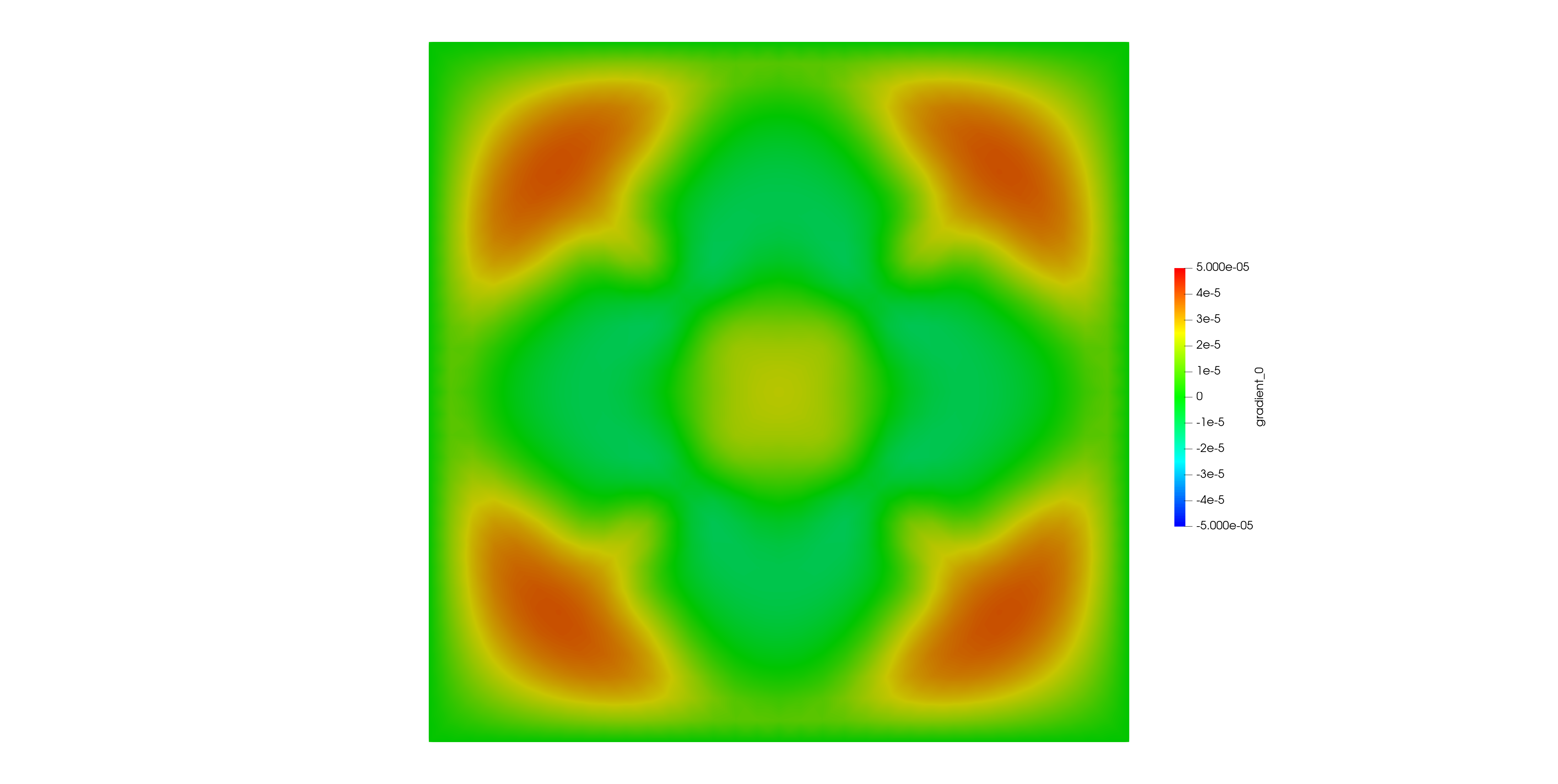}
		\caption{Iteration 203.}
		\label{fig:topo_der_first}
	\end{subfigure}%
	\hfil
	\begin{subfigure}{0.5\textwidth}
		\centering
		\includegraphics[width=\textwidth, trim=39cm 0cm 39cm 0cm, clip]{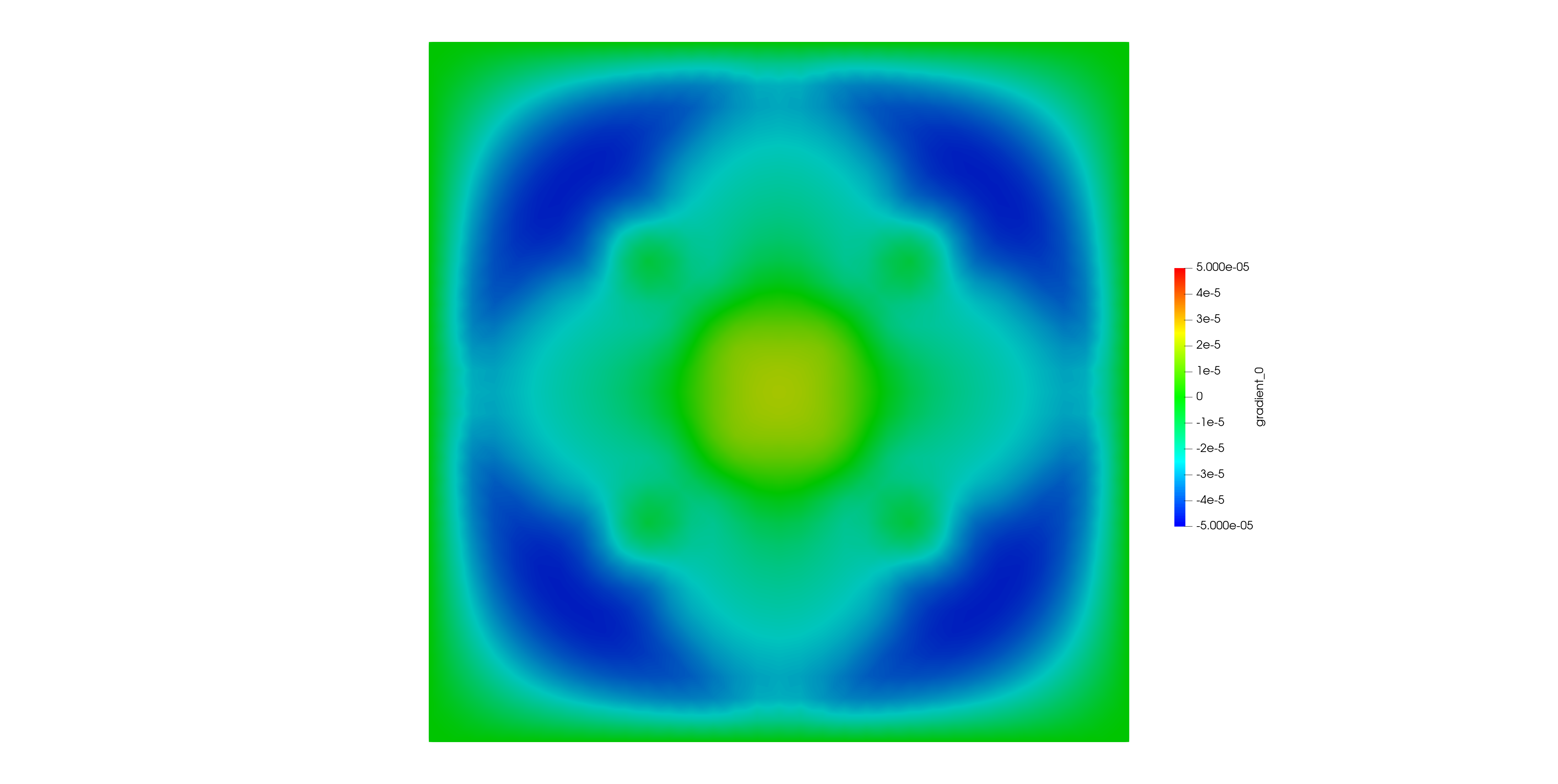}
		\caption{Iteration 204.}
		\label{fig:topo_der_second}
	\end{subfigure}
	\caption{Plots of the topological derivative at two iterations of the sphere combination algorithm for problem~\eqref{eq:state_system_poisson}.}
	\label{fig:topo_der}
\end{figure}

\begin{Remark}
	A possible explanation for this behavior is the fact, that problem \eqref{eq:state_system_poisson} is well-known to be ill-conditioned. This means that topological derivatives become flatter the closer we get to the optimal solution of the problem, which is also what we have observed in our numerical experiments. The topological derivative becomes successively smaller and flatter over the course of the optimization. 
	
	In Figure~\ref{fig:topo_der}, we show plots of the topological derivative during the middle of the optimization for the sphere combination method. In Figure~\ref{fig:topo_der_first}, the topological derivative is shown at iteration 203 and for this iterate, the angle between topological derivative and level-set function is comparatively small with \SI{29}{\degree}. On the other hand, the topological derivative in the next iteration, which is shown in Figure~\ref{fig:topo_der_second}, has a much larger angle of \SI{130}{\degree} with the level-set function. A possible reason for this behavior could be the line search, which, at first, uses successively larger step sizes up to iteration 204, after that the accepted step size drops down. This behavior repeats and causes the oscillations in the angle criterion which can be seen in Figure~\ref{fig:poisson_nonprincipal}.
	
	Moreover, we remark that the angle criterion only provides a necessary condition for optimality, as can be easily seen by the discussion in Section~\ref{ssec:topological_sensitivity_analysis} (cf.~\eqref{eq:opt_levelset} and~\eqref{eq:fonc}). There may be other minimizers that satisfy \eqref{eq:fonc} without satisfying the angle criterion. A thorough investigation of these issues is an interesting direction for future work.
\end{Remark}

The  evolution of the geometry during the optimization algorithms is depicted in Table~\ref{tab:comparison_poisson}, where we show the geometries after 50, 100, and 500 iterations for all four considered optimization algorithms on the finest discretization. In addition, we also show the reference shape for comparison. Note, that the reference shape has four larger inclusions in the middle and a smaller one in the center, making the problem particularly hard to solve. Comparing the obtained shapes after 50 iterations, we observe that the sphere combination and convex combination algorithms only start to form the four major inclusions, whereas they are already present for the gradient descent and BFGS methods. Moreover, the geometry obtained by BFGS method already has an inclusion in the center of the geometry, making it already very similar to the reference solution. After 100 iterations, the established methods have formed the major inclusions, but their position and shape is still quite wrong. The gradient descent method has improved the shapes of the major inclusions so that they are rather similar to the desired ones and is starting to locate the inclusion in the center. For the BFGS method, there is no visible difference anymore between the obtained and desired geometry, indicating that the method has already converged after 100 iterations, whereas the other methods are still quite far away from the desired shape. After 500 iterations, we observe that the sphere combination and convex combination methods have corrected the shape of the four major inclusions so that they are now quite similar to the desired ones. However, neither of these methods was able to reconstruct the inclusion in the center of the geometry. The gradient descent method, on the other hand, was able to reconstruct the center inclusion after 500 iterations, and there are no visual distinctions between the obtained and desired geometry anymore. The same is, of course, also true for the BFGS method, whose corresponding shape does not change visually between iterations 100 and 500.

\begin{figure}[!t]
	\centering
	\begin{subfigure}{0.25\textwidth}
		\centering
		\includegraphics[width=\textwidth]{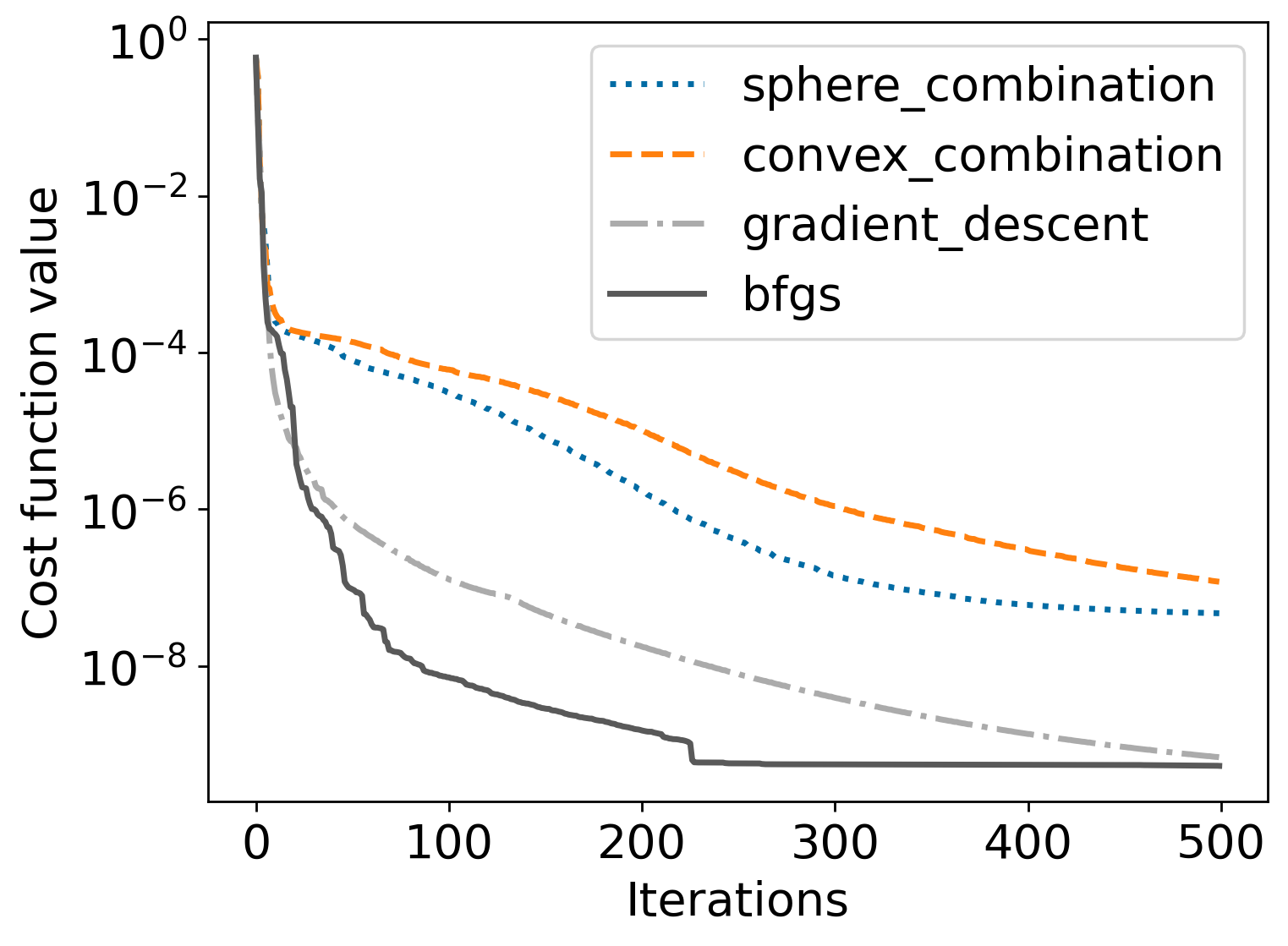}
		\caption{$32\times 32$}
	\end{subfigure}%
	\begin{subfigure}{0.25\textwidth}
		\centering
		\includegraphics[width=\textwidth]{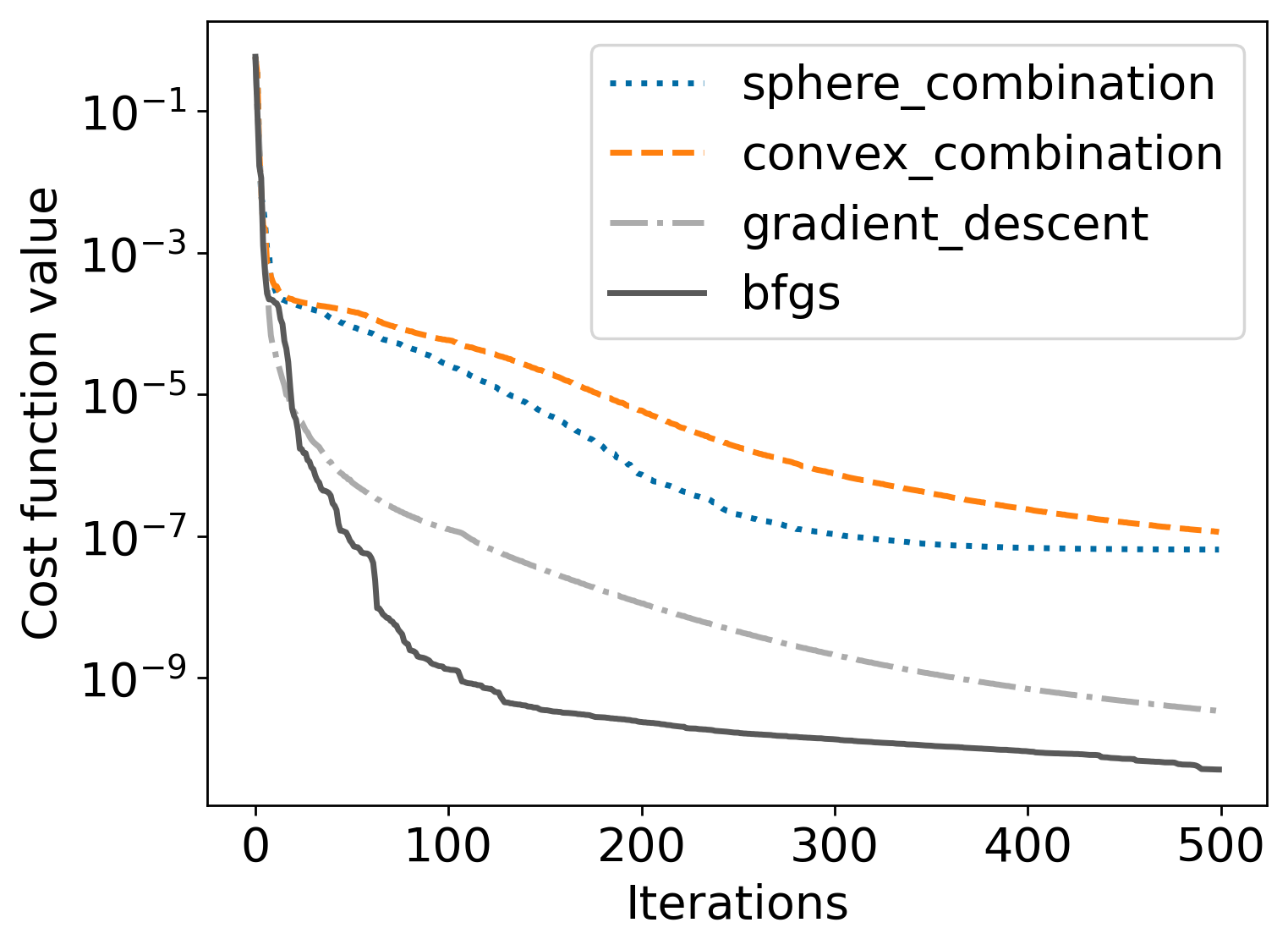}
		\caption{$48\times 48$}
	\end{subfigure}%
	\begin{subfigure}{0.25\textwidth}
		\centering
		\includegraphics[width=\textwidth]{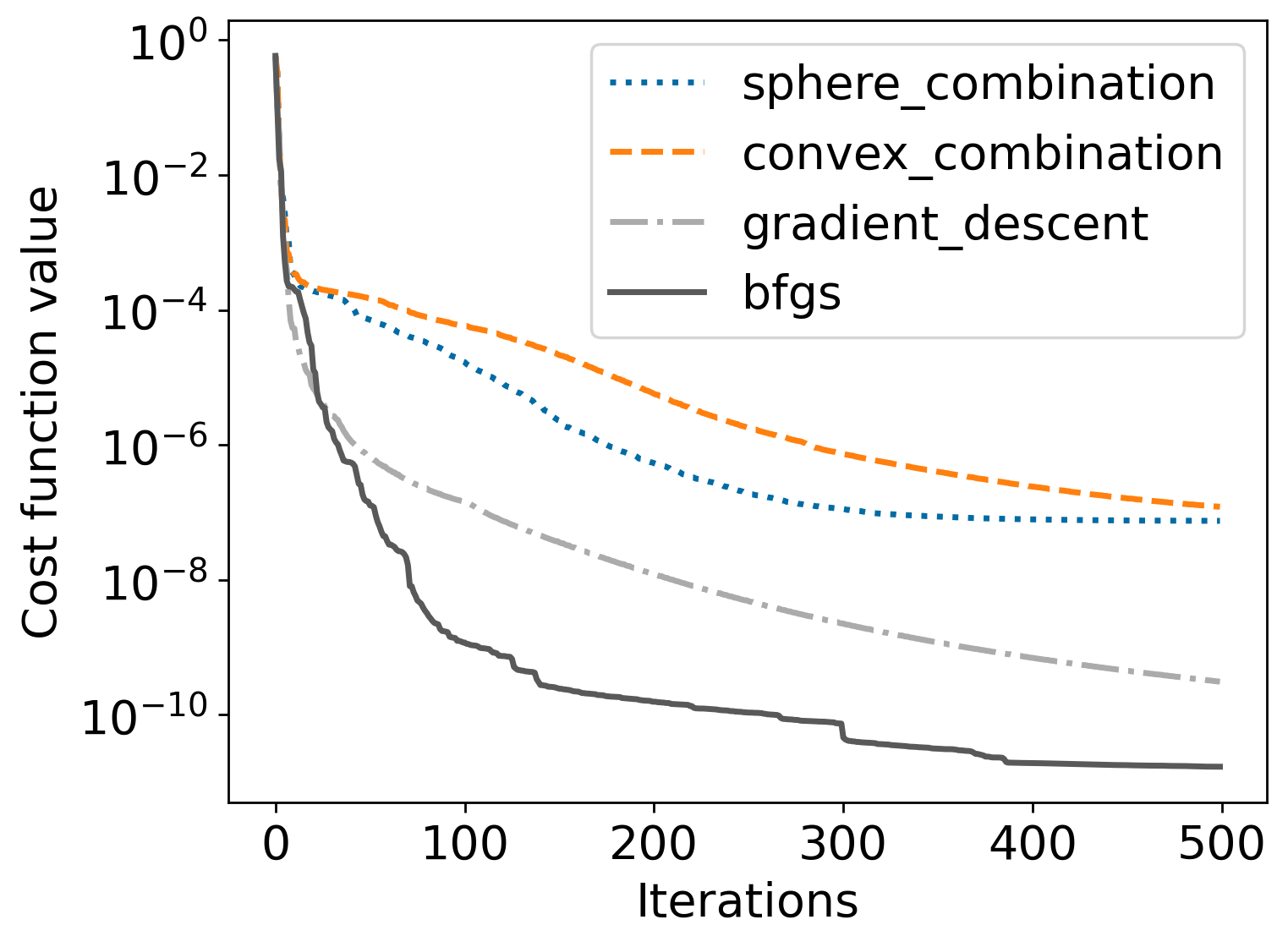}
		\caption{$64\times 64$}
	\end{subfigure}%
	\begin{subfigure}{0.25\textwidth}
		\centering
		\includegraphics[width=\textwidth]{linear_poisson_96_cost_functional}
		\caption{$96\times 96$}
	\end{subfigure}%
	\caption{Evolution of the cost functional for the linear Poisson problem \eqref{eq:state_system_poisson} for several discretization levels.}
	\label{fig:mesh_independence}
\end{figure}

Let us finally investigate the dependence of the algorithms on the discretization. We do so by comparing the evolution of the cost functional for all methods on the different discretization levels. The results are depicted in Figure~\ref{fig:mesh_independence}. Here, we observe that the performance of the algorithms, which we have discussed previously, is not dependent on the mesh size. 
In particular, the BFGS and gradient descent methods always perform best, whereas the sphere and convex combination methods perform significantly worse. Overall, the cost functional decreases further when a finer discretization is chosen. This is due to the fact, that a finer discretization also allows for a more detailed resolution of the desired shape, which in turn results in a lower cost functional value for the optimum. Therefore, we conclude that all algorithms behave mesh-independently, i.e., they do not require more iterations for finer levels of discretization to reach the same level of approximation of the optimal solution.

These results show the great potential of the proposed BFGS methods as they performed best of all considered methods, significantly outperforming the remaining algorithms. Moreover, we have shown that the methods also show mesh-independent behavior, which is due to the fact that the presented algorithms are discretizations of optimization algorithms acting in infinite-dimensional spaces.


\subsection{Semilinear Poisson Problem}
\label{ssec:semilinear_transmission}

To showcase the methods' performance for nonlinear problems, we now investigate a semilinear variant of the Poisson problem we considered in Section~\ref{ssec:numerics_model}. Our setting is like before, so $\Dsf \subset \R^d$ with $d \in \mathbb{N}_{>0}$ and $d\leq 3$ denotes the holdall domain, $\partial \Dsf$ is its boundary, $\Omega \subset \Dsf$ is an open subset of $\Dsf$, $\Gamma = \partial \Omega \setminus \partial \Dsf$ is the interior boundary of $\Omega$ in $\Dsf$, and $\Omega^c = \Dsf \setminus \closure{\Omega}$ is the complement of $\Omega$ in $\Dsf$. Our semilinear version of the problem reads
\begin{equation}
\label{eq:semilinear_poisson}
\begin{aligned}
&\min_{\Omega} J(\Omega, u) = \frac{1}{2} \integral{\Dsf} \left( u - u_\mathrm{des} \right)^2 \dmeas{x} \\
&\text{s.t. } \quad \begin{alignedat}[t]{2}
-\Delta u + \alpha_\Omega u^3 &= f_\Omega \quad &&\text{ in } \Dsf,\\
u &= 0 \quad &&\text{ on } \partial \Dsf,
\end{alignedat}
\end{aligned}
\end{equation}
where we again have $\alpha_\Omega(x) = \chi_\Omega(x) \alpha_\mathrm{in} + \chi_{\Omega^c}(x) \alpha_\mathrm{out}$ and $f_\Omega(x) = \chi_\Omega(x) f_\mathrm{in} + \chi_{\Omega^c}(x) f_\mathrm{out}$ with $\alpha_\mathrm{in}, \alpha_\mathrm{out} > 0$ as well as $f_\mathrm{in}, f_\mathrm{out} \in \mathbb{R}$. The topological derivative for this problem is given by
\begin{equation*}
\mathcal{D}J(\Omega)(x) = (\alpha_\mathrm{out} - \alpha_\mathrm{in}) u(x)^3 p(x) - (f_\mathrm{out} - f_\mathrm{in})p(x) \quad \text{ for all } x \in \Dsf \setminus \Gamma,
\end{equation*}
where $u$ solves the state equation and $p$ solves the following adjoint equation
\begin{equation*}
\begin{alignedat}{2}
-\Delta p + 3\alpha_\Omega u^2 p &= - (u - u_\mathrm{des}) \quad &&\text{ in } \Dsf,\\
p &= 0 \quad &&\text{ on } \partial \Dsf.
\end{alignedat}
\end{equation*}

\begin{figure}[!t]
	\centering
	\begin{subfigure}{0.333\textwidth}
		\centering
		\includegraphics[width=\textwidth]{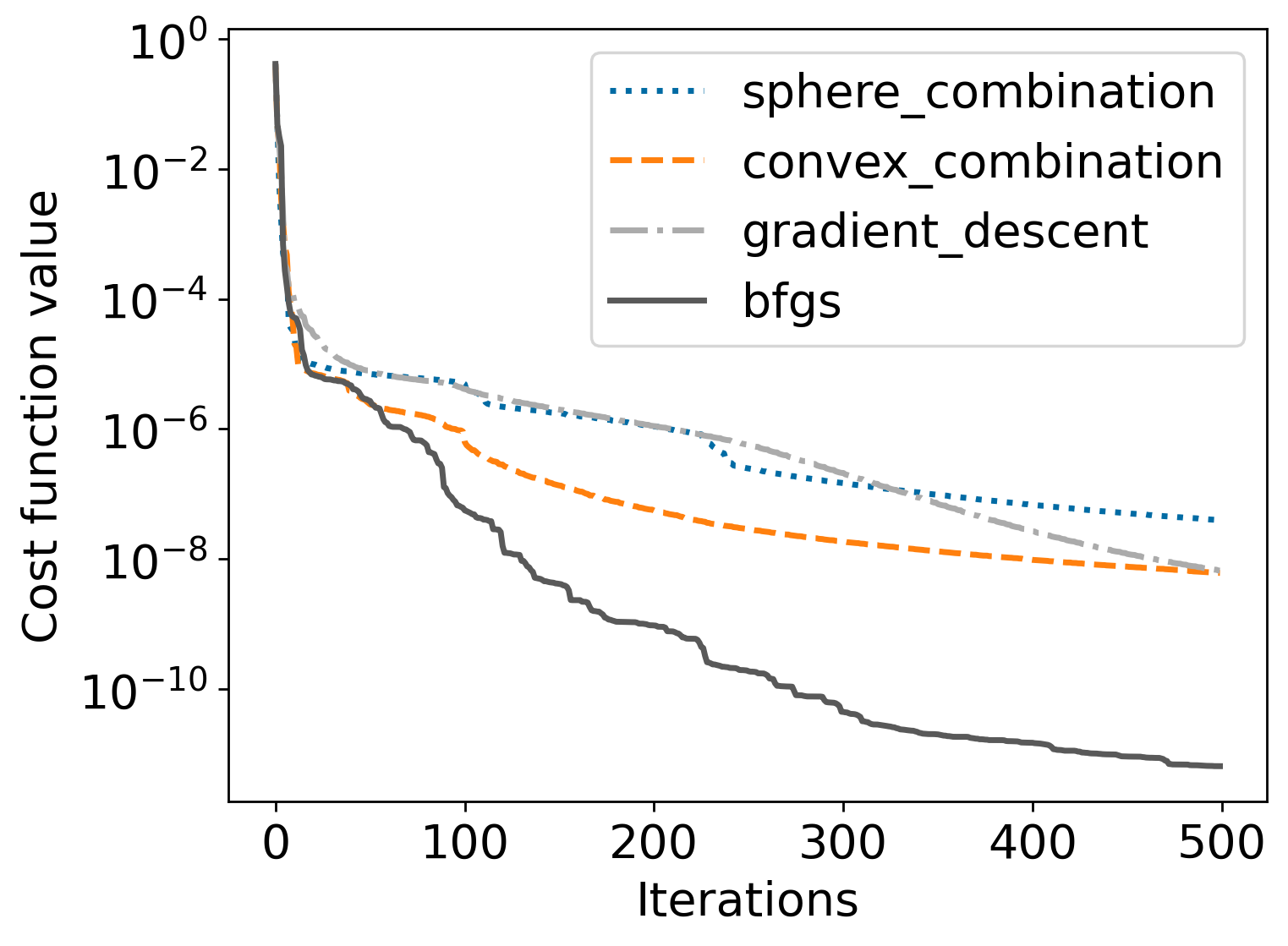}
		\caption{Cost functional.}
	\end{subfigure}%
	\begin{subfigure}{0.333\textwidth}
		\centering
		\includegraphics[width=\textwidth]{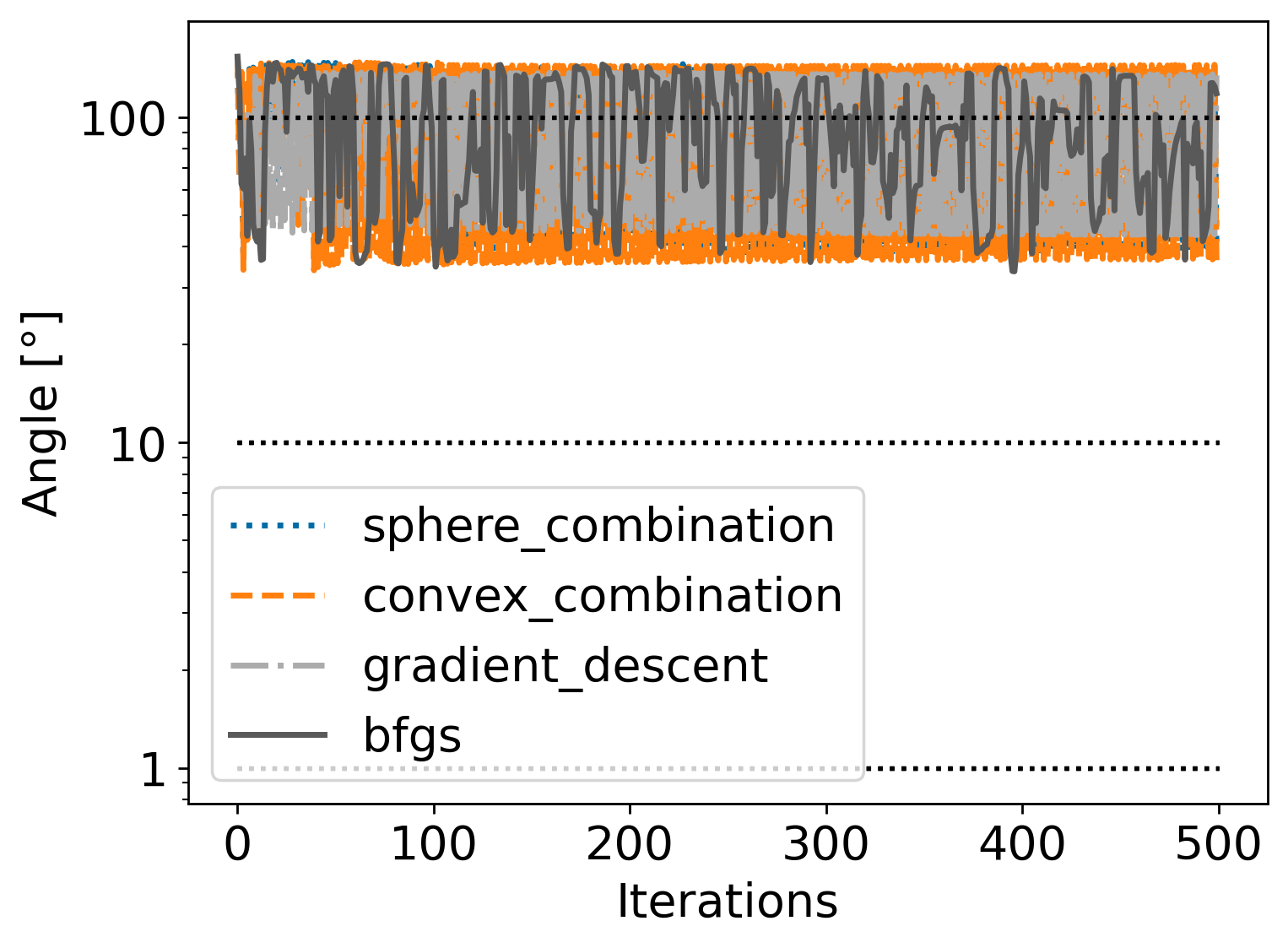}
		\caption{Angle.}
	\end{subfigure}%
	\begin{subfigure}{0.333\textwidth}
		\centering
		\includegraphics[width=\textwidth]{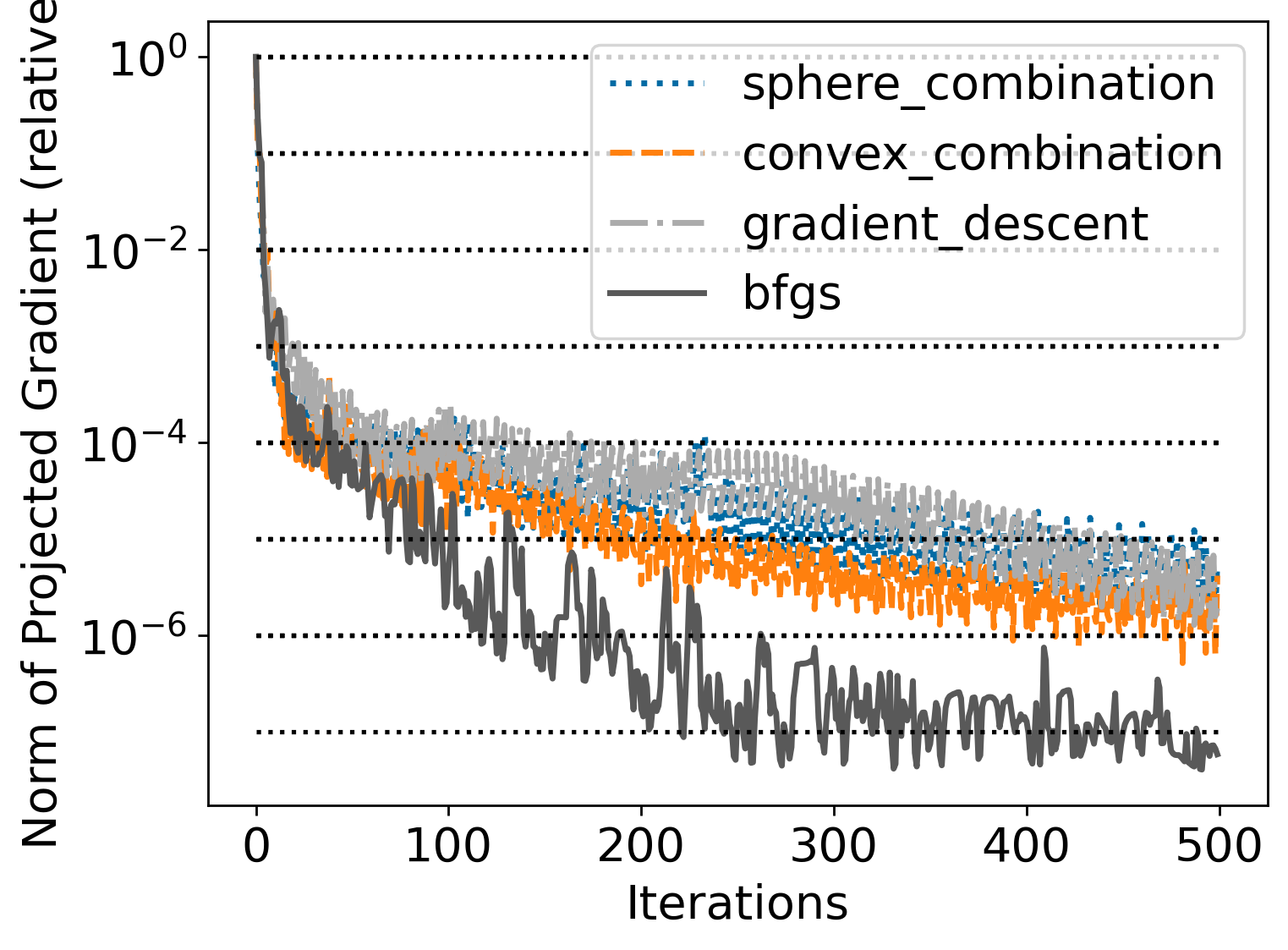}
		\caption{Norm of the projected topological derivative.}
	\end{subfigure}
	\caption{Evolution of the optimization for the semilinear Poisson problem \eqref{eq:semilinear_poisson}.}
	\label{fig:poisson_semilinear}
\end{figure}

We again solve this problem with the four solution algorithms under consideration and use the hold-all domain $\Dsf = (-2, 2)^2$. The desired state $u_\mathrm{des}$ is obtained by solving the semilinear Poisson problem on a reference domain, which is given by the same clover shape as considered in Section~\ref{ssec:numerics_model}. The geometry is discretized by dividing it into $96\times 96$ squares, which are subdivided into four triangles each, so that we use \num{18625}~nodes and \num{36864}~triangles. Again, we employ linear Lagrange elements for the discretization of the state and adjoint variables. Finally, we use the same setting for $\alpha$ and $f$ as before, so that $\alpha_\mathrm{in} = f_\mathrm{in} = 10$ and $\alpha_\mathrm{out} = f_\mathrm{out} = 1$.

The evolution of the cost functional, angle criterion, and norm of the projected topological derivative over the course of the optimization are shown in Figure~\ref{fig:poisson_semilinear}. Analogously to our previous findings, we observe that the BFGS method significantly outperforms the remaining methods as it decreases the cost functional and the norm of the projected topological derivative most. For this problem, we observe that the gradient descent method does not perform as well as previously, but is rather comparable to the sphere and convex combination methods in performance. Again, for all considered algorithms, the angle between level-set function and generalized topological derivative remained bounded away from zero.

Let us investigate the geometries obtained by the methods, which are depicted in Table~\ref{tab:comparison_semilinear} after 100, 200, and 500 iterations of the methods. Here, we observe that the BFGS method again performs substantially better than the other methods. Even after 100 iterations, all inclusions of the geometry are found and the geometry exhibits the correct topology. However, the shape of the inclusions is still a bit off. The geometry obtained with the remaining algorithms, on the other side, is still far away from the optimal geometry. The same is true after 200 iterations: There, the geometry obtained by the BFGS method is visually identical to the reference solution, whereas the other methods produce geometries that still deviate significantly, both in topology and shape, from the optimal solution. Finally, the results after 500 iterations show, that the sphere and convex combination methods are able to locate the four major inclusions of the reference solution, but their shape is still a bit off. Moreover, they were unable to detect the inclusion in the center. The gradient descent method, on the other hand, was able to find all inclusions and, hence, to reconstruct the sought topology. However, the shape of the inclusions could still be improved. For the BFGS method, there are no major changes in the geometry between 200 and 500 iterations, as the geometry showed no visual differences to the reference solution already after 200 iterations.

Again, these results highlight the potential and efficiency of the BFGS methods for solving topology optimization problems, as they significantly outperformed all other methods considered in this paper.

\begin{table}[!t]
	\newcommand{\sizebox}{0.215\textwidth}
	\newcolumntype{C}{>{\centering\arraybackslash}m{\sizebox}}
	\setlength{\tabcolsep}{0pt}
	\caption{Evolution of the geometries for the semilinear Poisson problem \eqref{eq:semilinear_poisson}.}
	\label{tab:comparison_semilinear}
	\begin{tabular}{m{0.1\textwidth} @{\hskip 0.5em} C C C @{\hskip 1em} C}
		\toprule
		Algorithm & 100 iterations & 200 iterations & 500 iterations & Reference \\ 
		\midrule
		Sphere combination & \includegraphics[width=\sizebox]{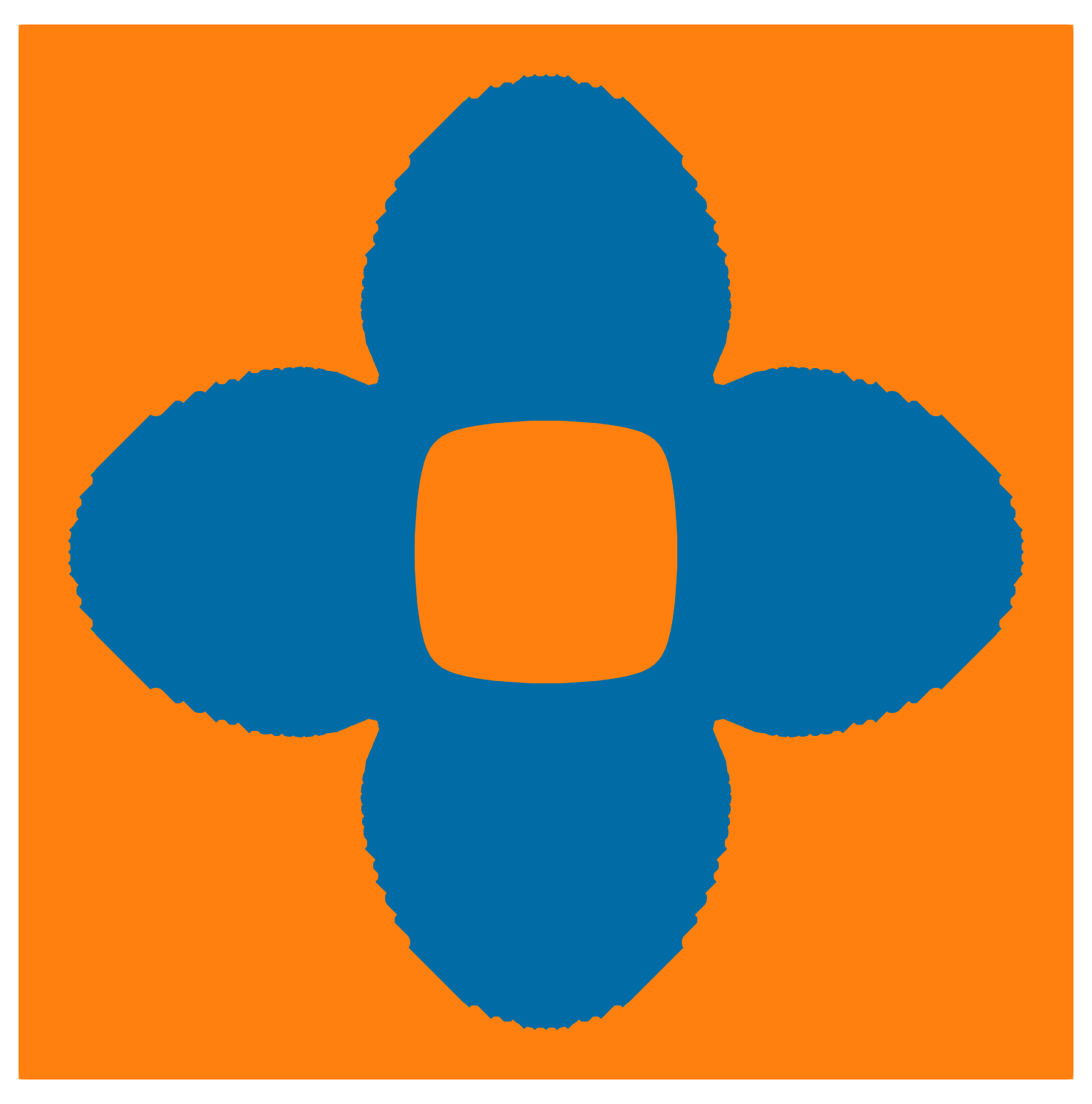} & \includegraphics[width=\sizebox]{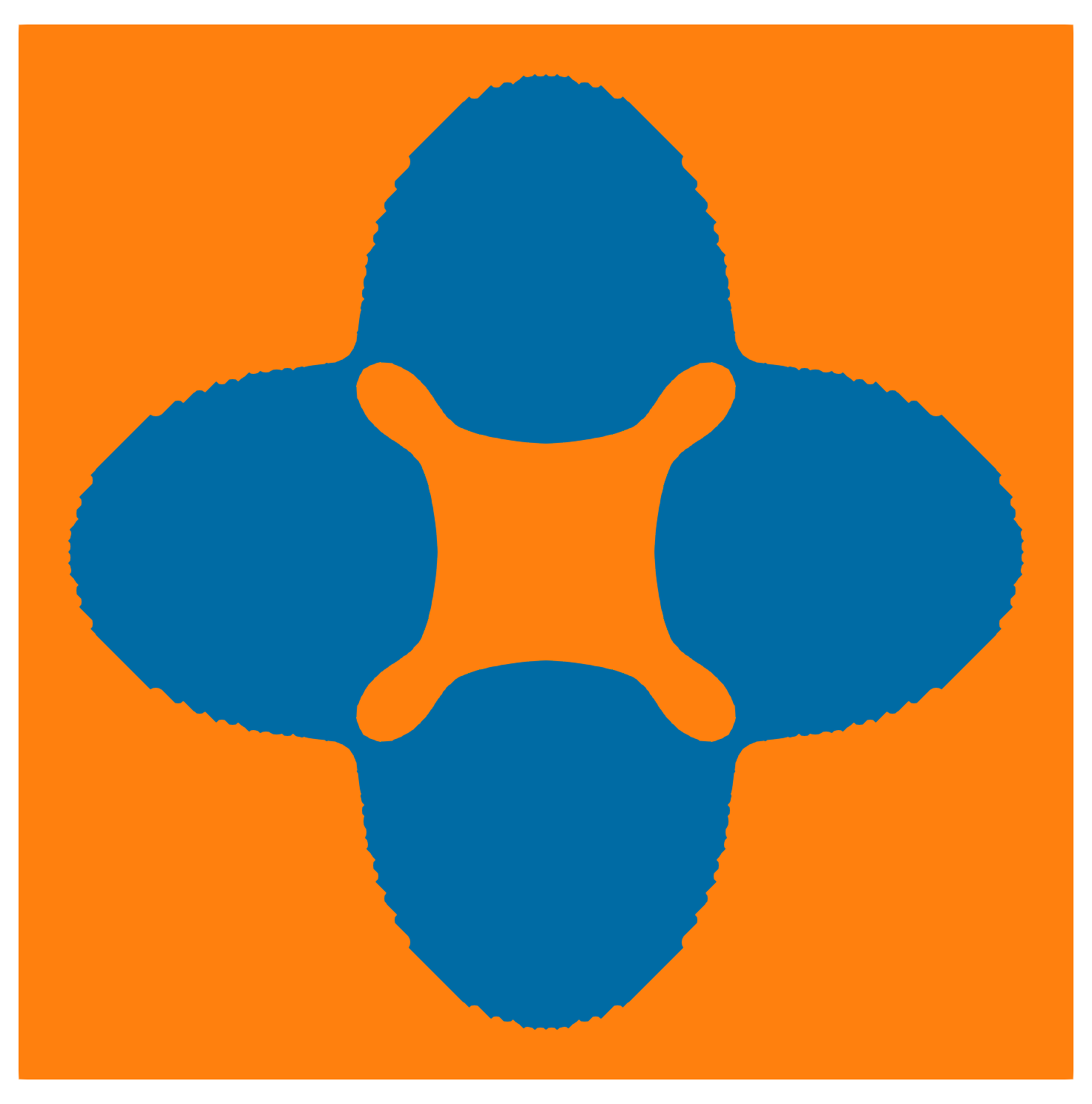} & \includegraphics[width=\sizebox]{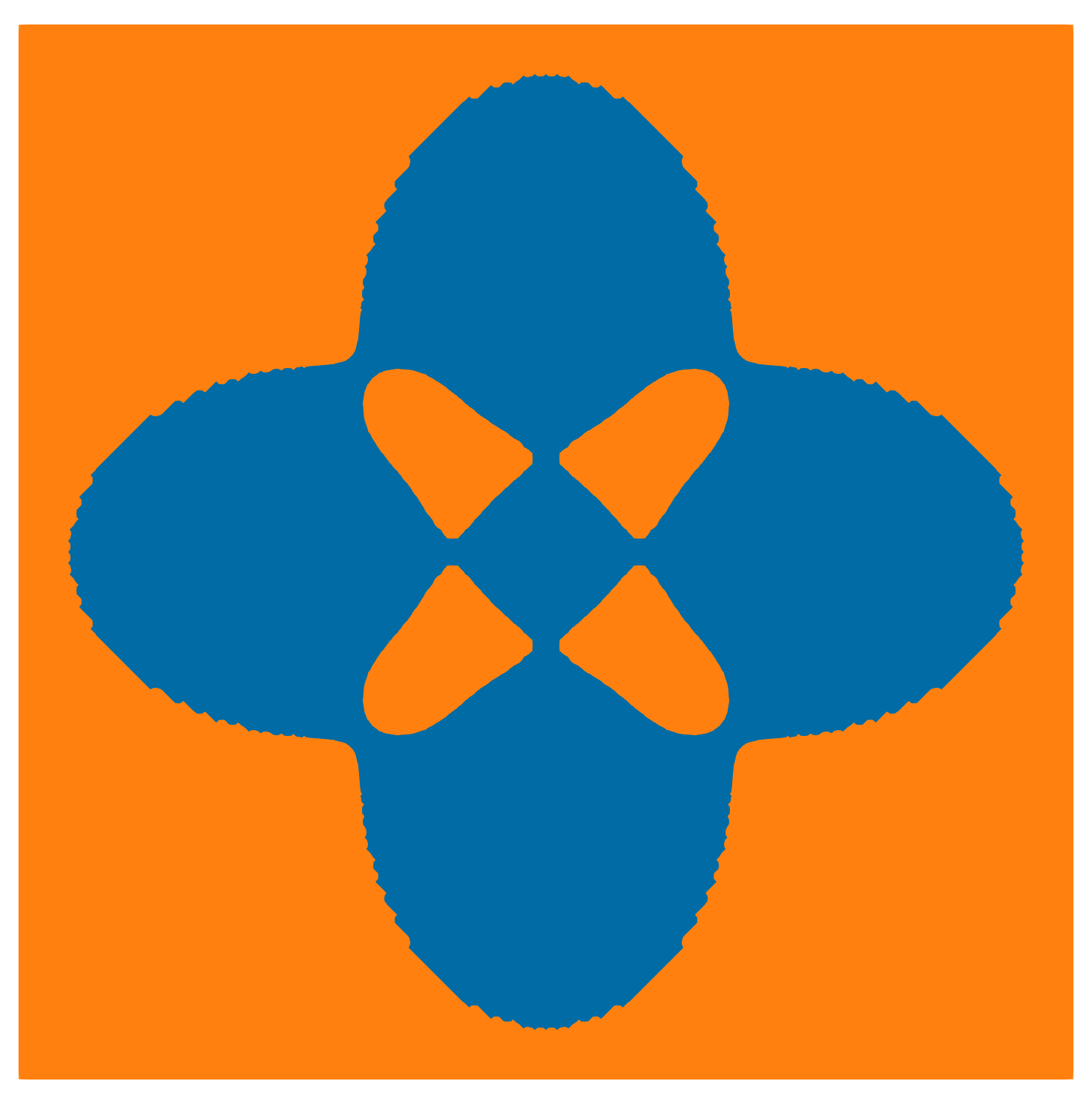} & \includegraphics[width=\sizebox]{clover} \\ 
		Convex combination & \includegraphics[width=\sizebox]{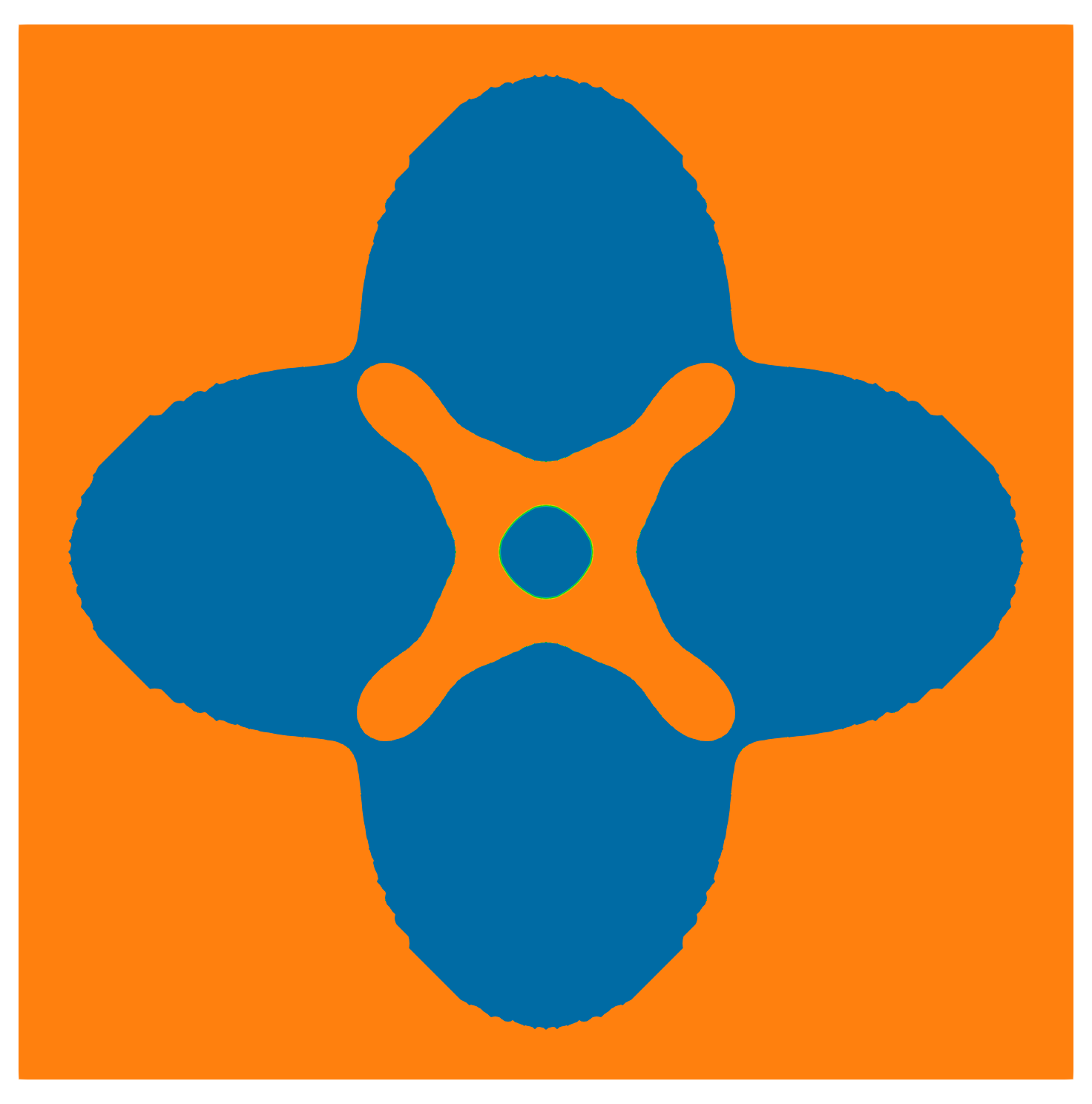} & \includegraphics[width=\sizebox]{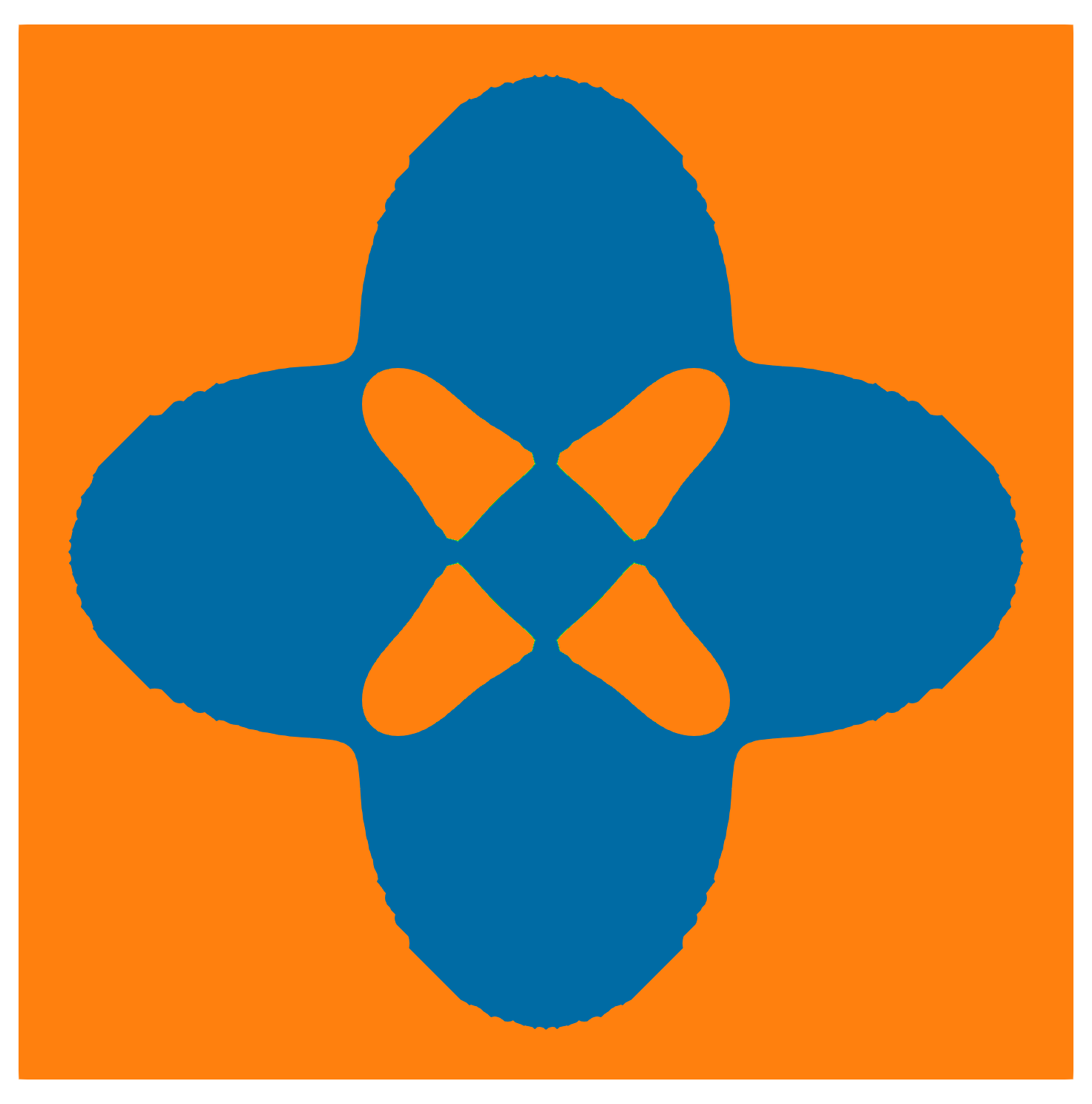} & \includegraphics[width=\sizebox]{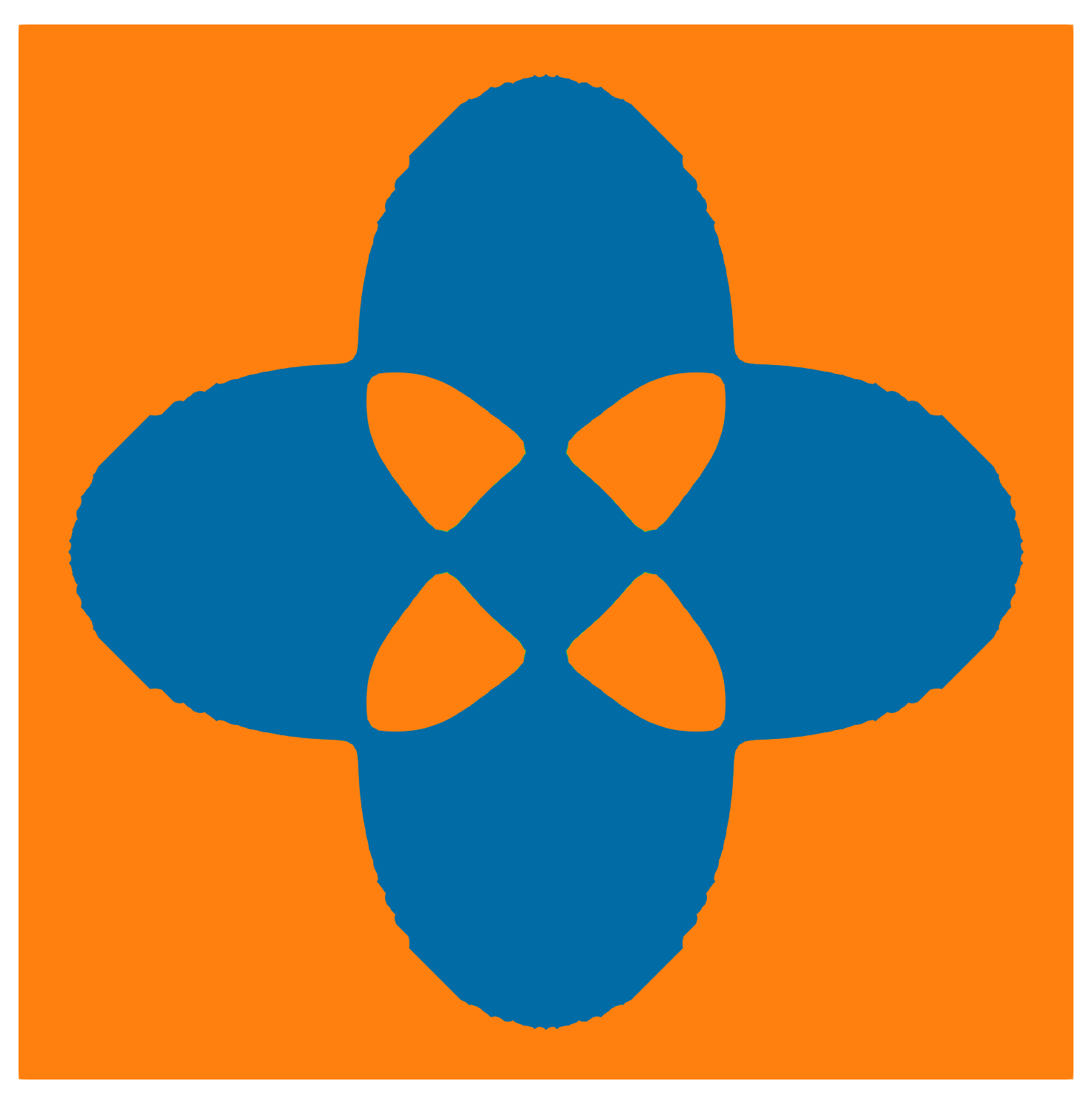} & \includegraphics[width=\sizebox]{clover} \\ 
		Gradient descent & \includegraphics[width=\sizebox]{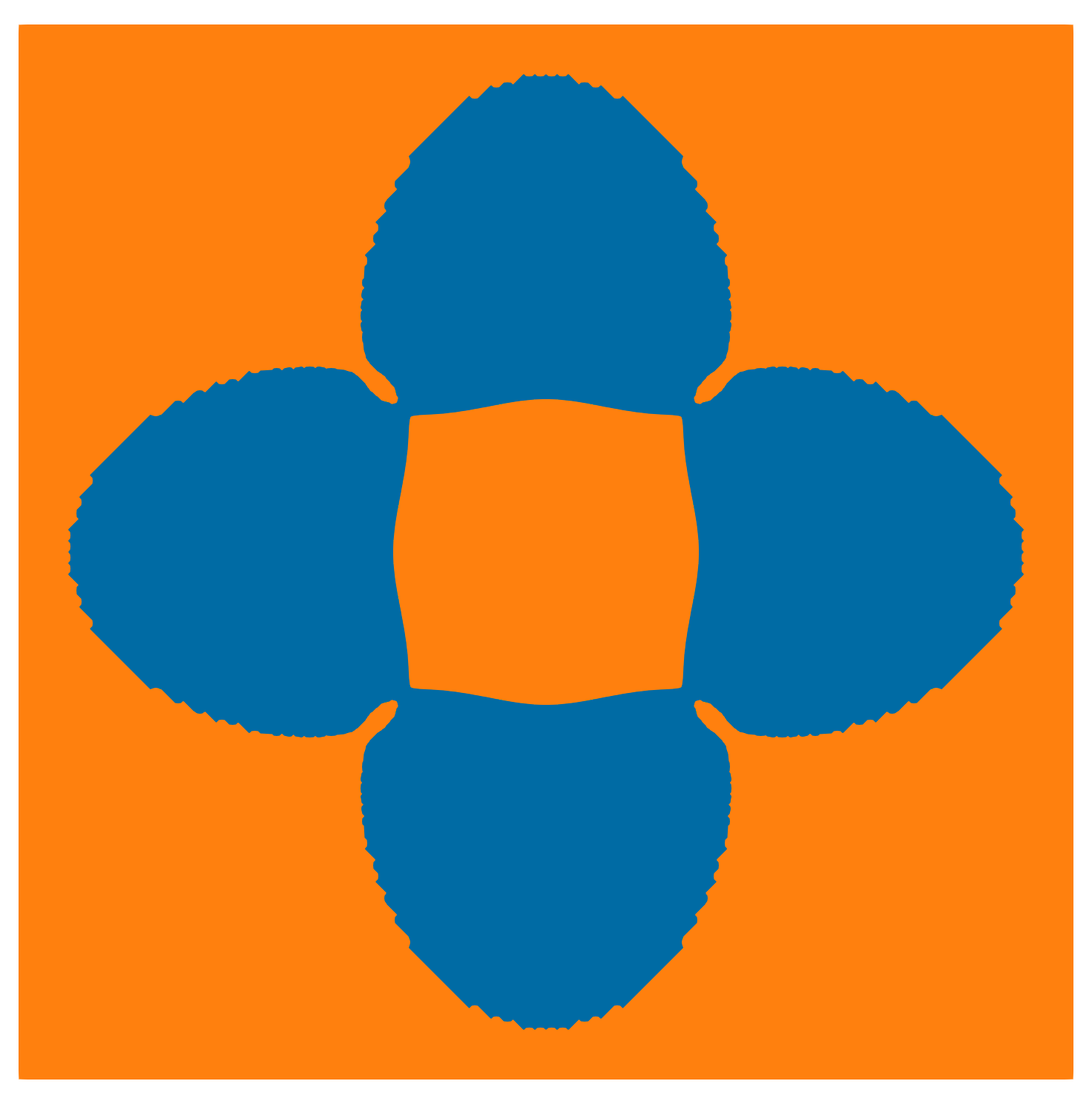} & \includegraphics[width=\sizebox]{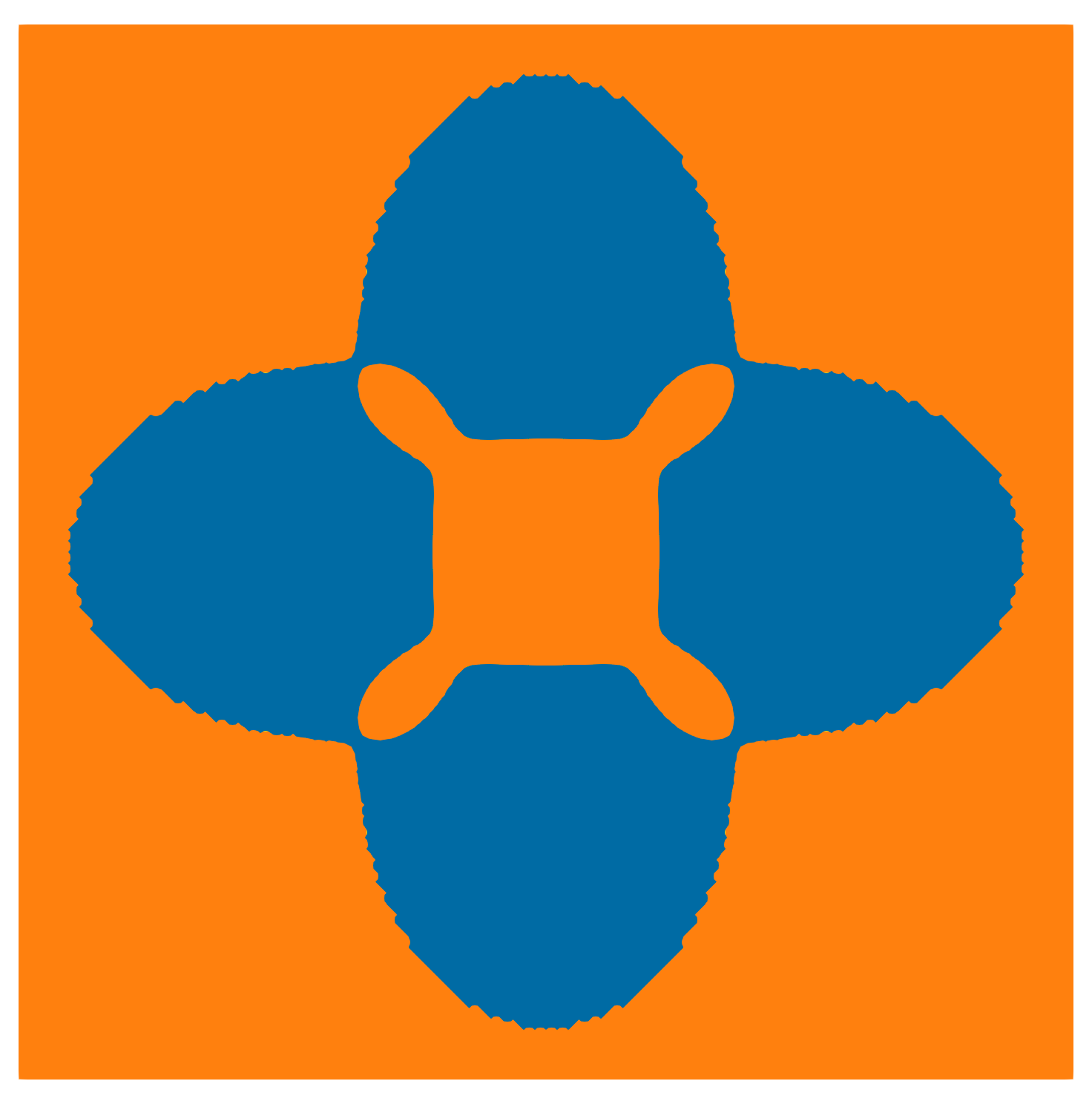} & \includegraphics[width=\sizebox]{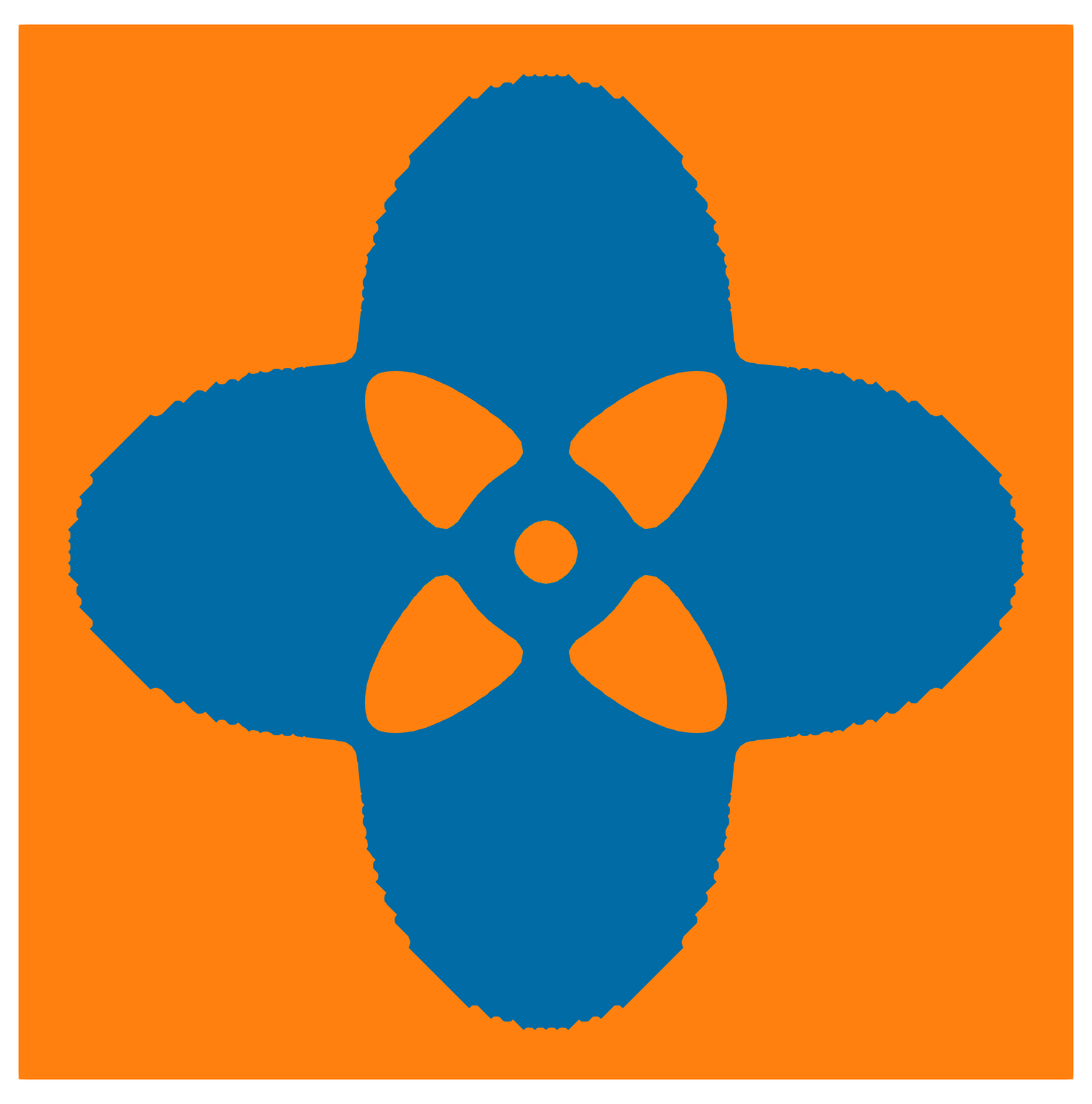} & \includegraphics[width=\sizebox]{clover} \\ 
		BFGS & \includegraphics[width=\sizebox]{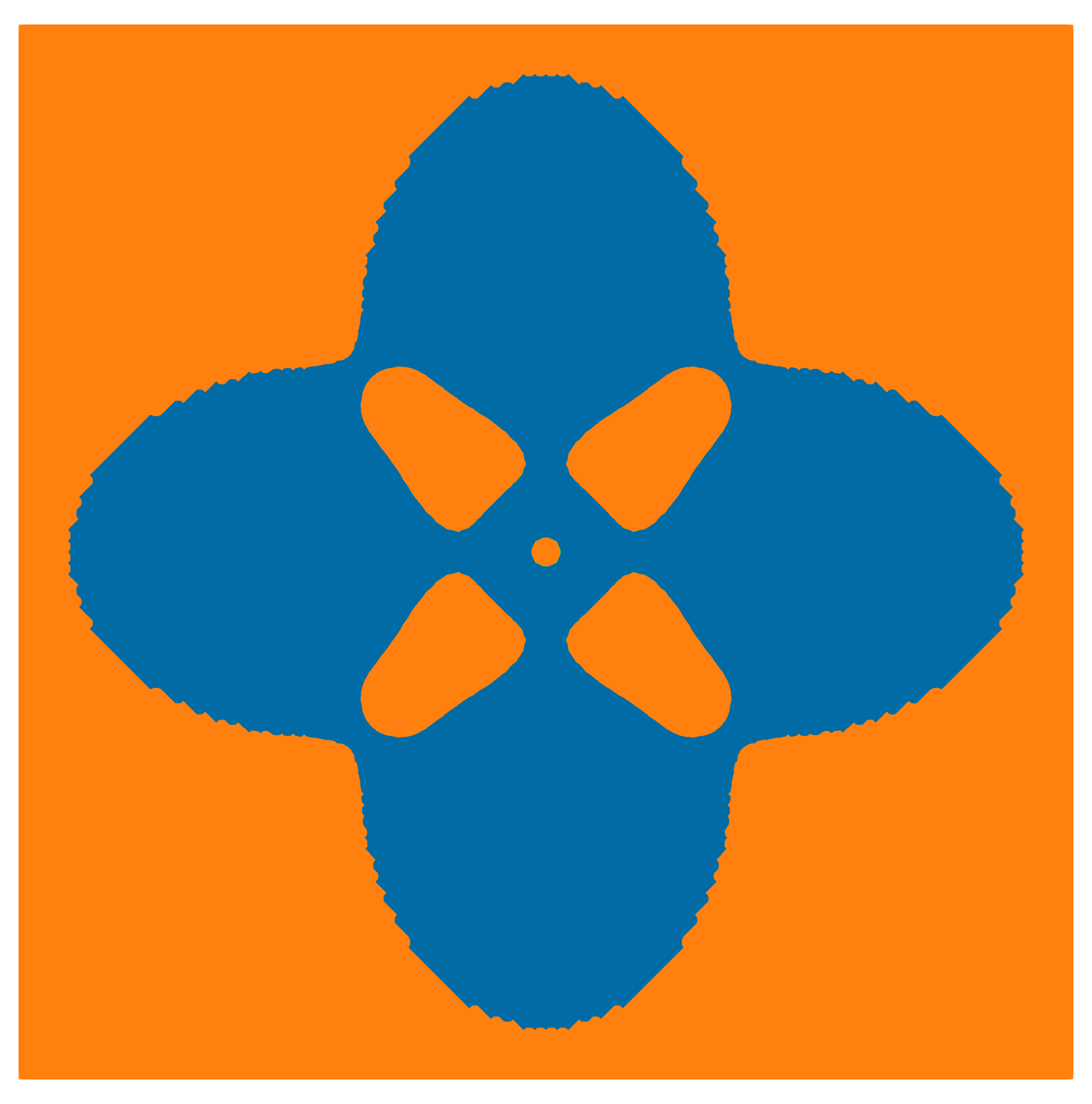} & \includegraphics[width=\sizebox]{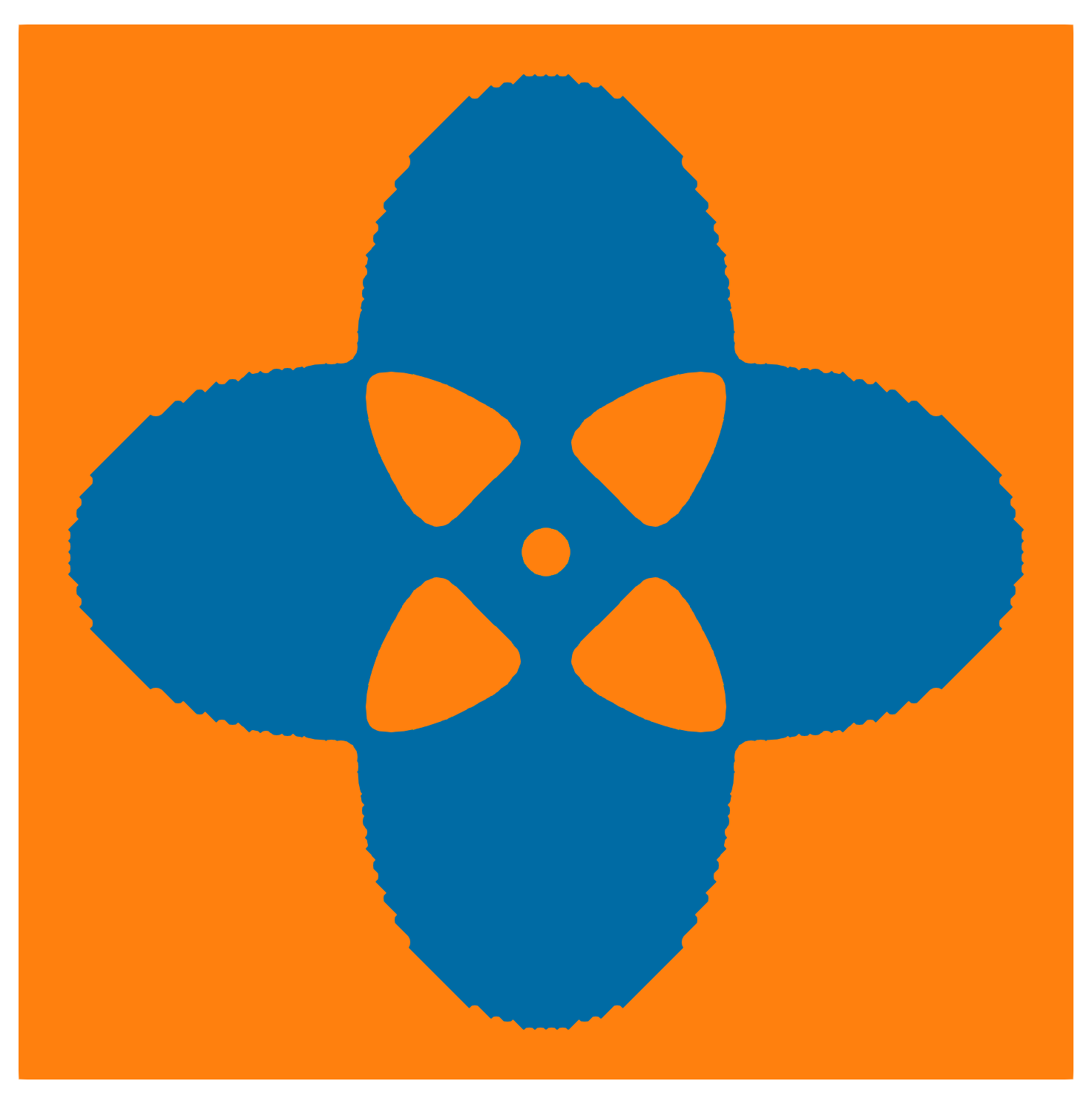} & \includegraphics[width=\sizebox]{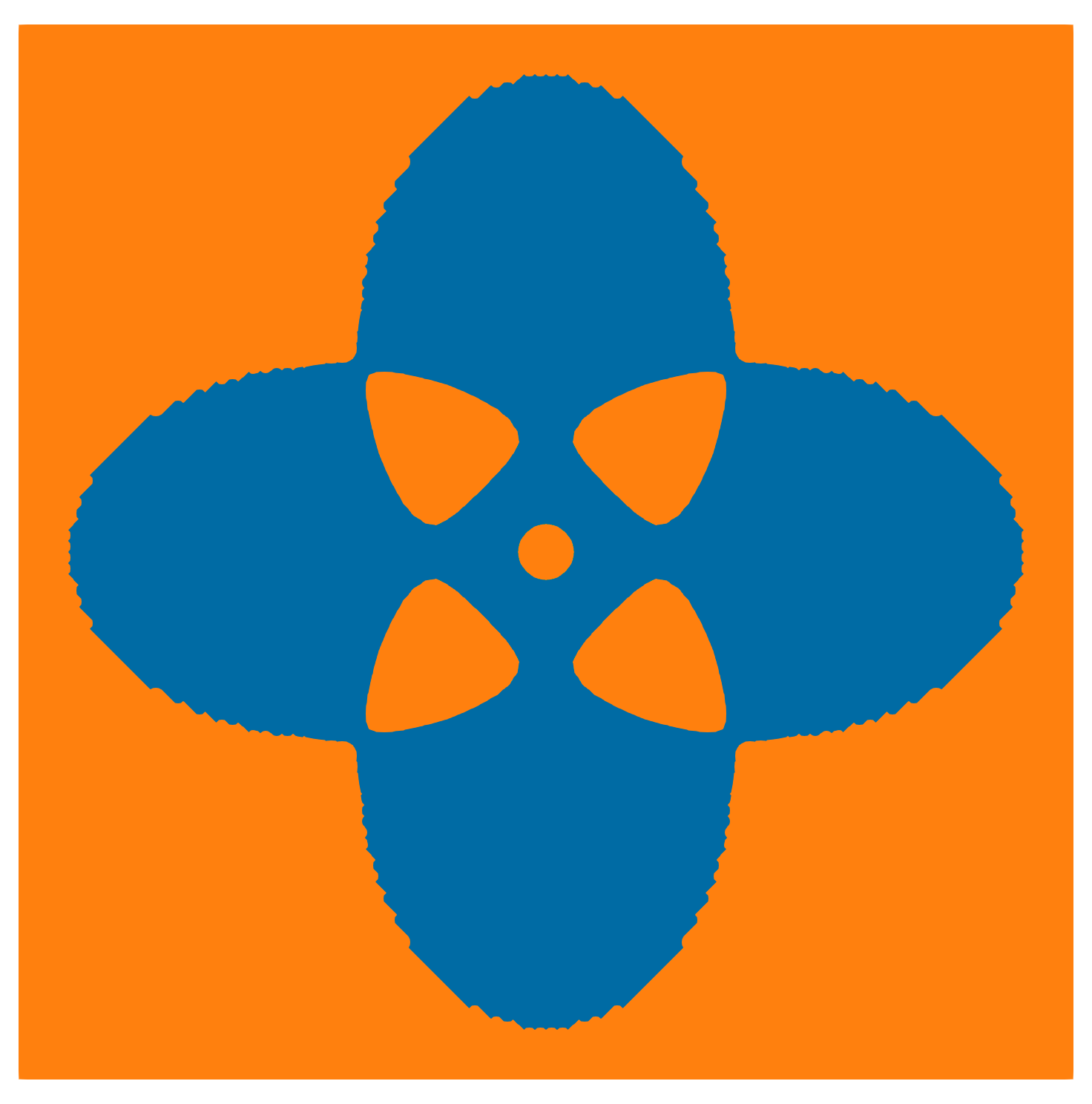} & \includegraphics[width=\sizebox]{clover} \\ 
		\bottomrule 
	\end{tabular}
\end{table}

\subsection{Compliance Minimization in Linear Elasticity}
\label{ssec:linear_elasticity}

In this section, we consider the problem of compliance minimization in linear elasticity which has been previously investigated, e.g., in \cite{Amstutz2006new, Allaire2004Structural, Allaire2005Structural}. Let $\Dsf \subset \R^d$, $d\in \mathbb{N}_{>0}$ be an open and bounded domain with boundary $\partial \Dsf$, which is the disjoint union of the Dirichlet boundary $\Gamma_D$ and Neumann boundary $\Gamma_N$. Further, let $\Omega \subset \Dsf$ and denote its complement by $\Omega^c = \Dsf \setminus \closure{\Omega}$. The compliance minimization problem is given by
\begin{equation}
\label{eq:topo_linear_elasticity}
\begin{aligned}
&\min_{\Omega} J(\Omega, u) = \int_{\Dsf} \alpha_\Omega \sigma(u) : e(u) \dmeas{x} + l \left|\Omega\right| \\
&\text{s.t.} \quad \begin{alignedat}[t]{2}
-\mathrm{div}(\alpha_\Omega \sigma(u)) &= f \quad &&\text{ in } \Dsf,\\
u &= 0 \quad &&\text{ on } \Gamma_D,\\
\alpha_\Omega \sigma(u) n &= g \quad &&\text{ on } \Gamma_N.\\
\end{alignedat}
\end{aligned}
\end{equation}
Here, $u$ is the deformation of a linear elastic material, $\sigma(u) = 2\mu e(u) + \lambda \tr{e(u)} I$ is Hooke's tensor, $e(u) = \nicefrac{1}{2}(\nabla u + (\nabla u)\transposed)$ is the symmetric gradient of $u$, and $A:B$ denotes the Frobenius inner product between matrices $A,B\in \mathbb{R}^{d\times d}$, i.e., $A:B:=\sum_{i,j=1}^d a_{ij}b_{ij}$. Here, $\mu$ and $\lambda$ are the so-called Lamé parameters for which we assume $\mu > 0$ and $2\mu + d \lambda > 0$. Moreover, $\alpha_\Omega(x) = \chi_\Omega(x) \alpha_\mathrm{in} + \chi_{\Omega^c}(x) \alpha_\mathrm{out}$ is, again, constant in $\Omega$ and $\Omega^c = \Dsf\setminus \closure{\Omega}$ with $\alpha_\mathrm{in}, \alpha_\mathrm{out} > 0$. The first term in the cost functional measures the compliance of the structure and the second term is a regularization parameter which penalizes large domains $\Omega$, so that the optimization is not trivial.

The topological derivative for problem \eqref{eq:topo_linear_elasticity} can be found, e.g., in \cite{Amstutz2006new} and is given by
\begin{equation*}
DJ(\Omega)(x) = \begin{cases}
-\alpha_\mathrm{in} \frac{r_\mathrm{in} - 1}{\kappa r_\mathrm{in} + 1} \frac{\kappa + 1}{2} \left( 2 \sigma(u) : e(u) + \frac{(r_\mathrm{in} - 1)(\kappa - 2)}{\kappa + 2 r_\mathrm{in} - 1} \tr{\sigma(u)} \tr{e(u)} \right) - l & \quad \text{ for } x\in \Omega, \\
-\alpha_\mathrm{out} \frac{r_\mathrm{out} - 1}{\kappa r_\mathrm{out} + 1} \frac{\kappa + 1}{2} \left( 2\sigma(u) : e(u) + \frac{(r_\mathrm{out} - 1)(\kappa - 2)}{\kappa + 2 r_\mathrm{out} -1} \tr{\sigma(u)} \tr{e(u)} \right) + l & \quad \text{ for } x \in \Omega^c = \Dsf \setminus \closure{\Omega},
\end{cases}
\end{equation*}
where $r_\mathrm{in} = \frac{\alpha_\mathrm{out}}{\alpha_\mathrm{in}}$, $r_\mathrm{out} = \frac{\alpha_\mathrm{in}}{\alpha_\mathrm{out}}$, and $\kappa = \frac{\lambda + 3\mu}{\lambda + \mu}$.

In the following, we consider two test cases for this problem, namely the so-called cantilever and bridge problems, which are taken from \cite{Amstutz2006new}. For these test cases, we follow \cite{Amstutz2006new} and use $\alpha_\mathrm{in} = 1$ and $\alpha_\mathrm{out} = \num{1e-3}$ as well as $f = 0$. Additionally, as we consider the case of plane stress, the Lam\'e parameters are given by $\mu = \frac{E}{2(1 + \nu)}$ and $\lambda = \frac{2 \mu \lambda^*}{\lambda^* + 2 \mu}$ with $\lambda^* = \frac{E \nu}{(1+ \nu) (1- 2\nu)}$ and we use a Young's modulus of $E = \num{1}$ and a Poisson's ratio of $\nu = 0.3$ for the following numerical experiments. Finally, we discretize the state variable using linear Lagrange finite elements.

\subsubsection{Cantilever}

\begin{figure}[!b]
	\centering
	\begin{tikzpicture}[scale=3]
	\draw [line width=0.4pt] (0,0) -- (2,0) -- (2,1) -- (0,1) -- cycle;
	
	\draw [dimen] (0, -0.1) -- (2, -0.1) node {2};
	\draw [dimen] (2.45, 0) -- (2.45, 1) node {1};
	
	\draw [-latex, line width=2.5pt] (2,0.5) -- (2, 0.125);
	\draw[pattern=north west lines, pattern color=black, draw=none] (0,0) rectangle (-0.25,1);
	
	\node at (0.125, 0.5) {{\LARGE $\Gamma_D$}};
	\node at (2.2, 0.5) {{\LARGE $\Gamma_N$}};
	\node at (1.0, 0.5) {{\LARGE $D$}};
	\end{tikzpicture}
	\caption{Schematic of the cantilever problem.}
	\label{fig:cantilever}
\end{figure}
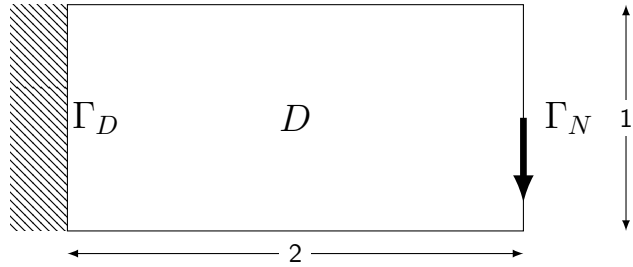

\begin{figure}[!b]
	\centering
	\begin{subfigure}{0.333\textwidth}
		\centering
		\includegraphics[width=\textwidth]{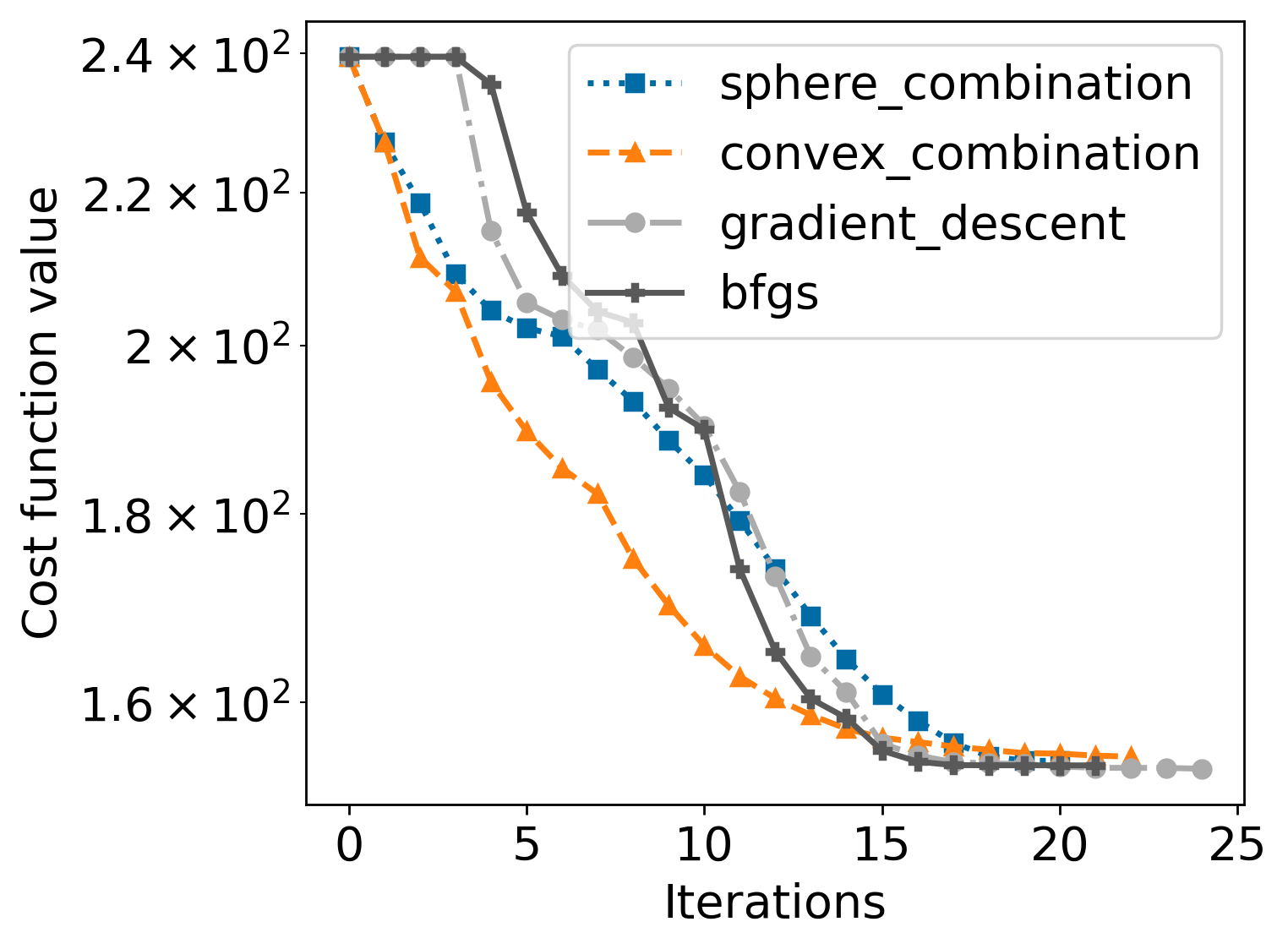}
		\caption{Cost functional.}
	\end{subfigure}%
	\begin{subfigure}{0.333\textwidth}
		\centering
		\includegraphics[width=\textwidth]{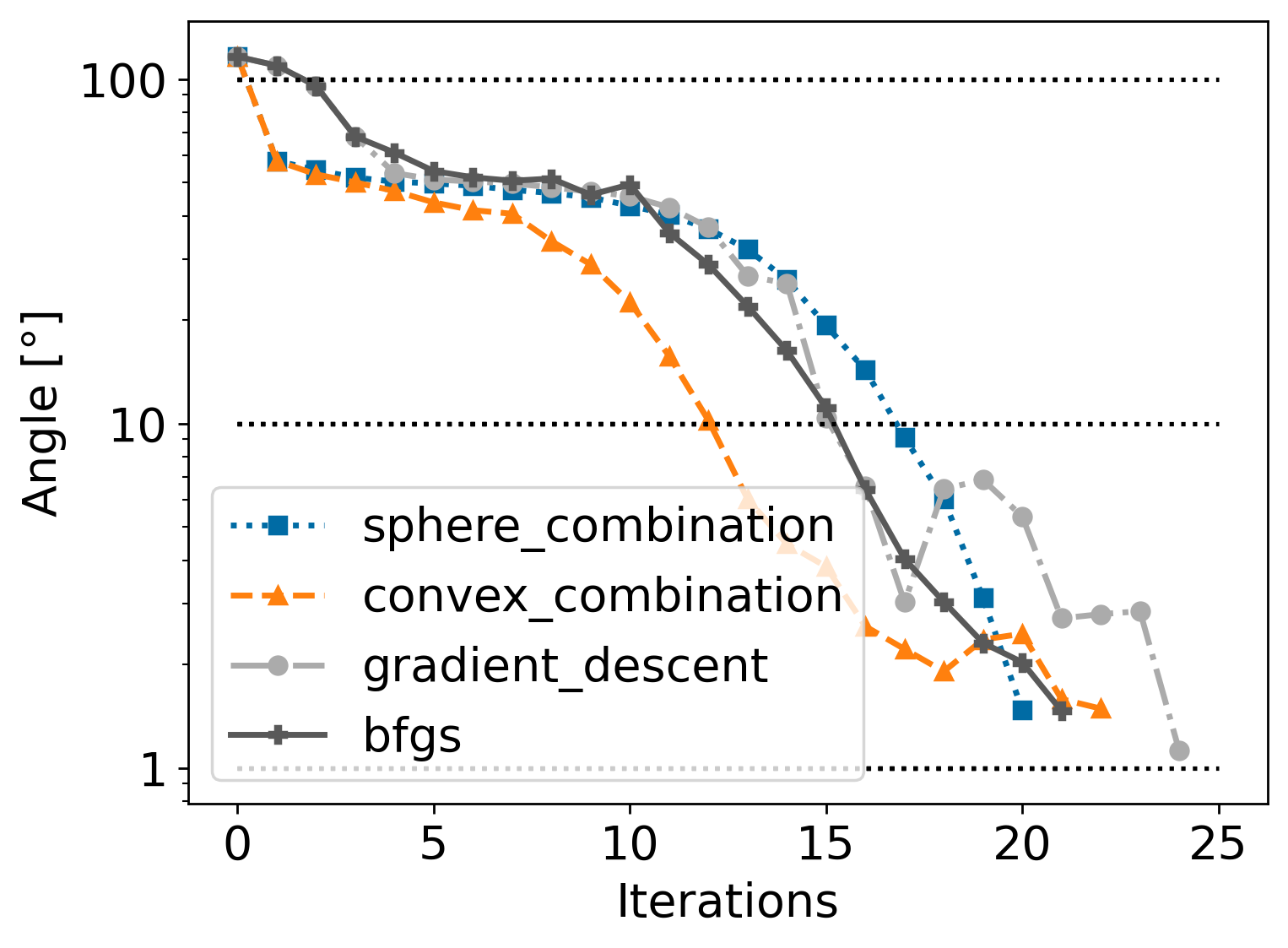}
		\caption{Angle.}
	\end{subfigure}%
	\begin{subfigure}{0.333\textwidth}
		\centering
		\includegraphics[width=\textwidth]{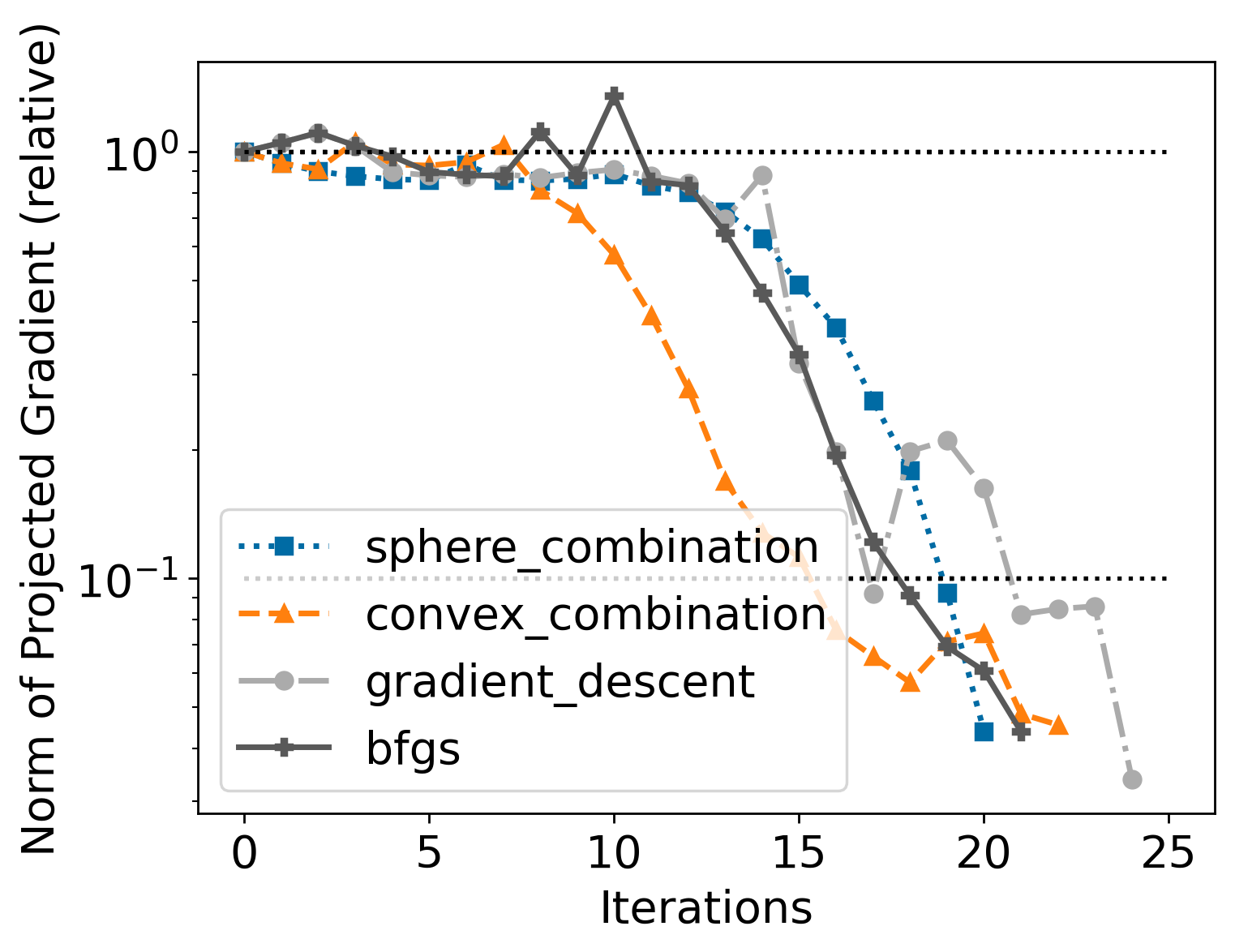}
		\caption{Norm of the projected topological derivative.}
	\end{subfigure}
	\caption{Evolution of the optimization for the cantilever problem.}
	\label{fig:evolution_cantilever}
\end{figure}

\begin{figure}[!t]
	\centering
	\begin{subfigure}{0.25\textwidth}
		\centering
		\includegraphics[width=0.95\textwidth]{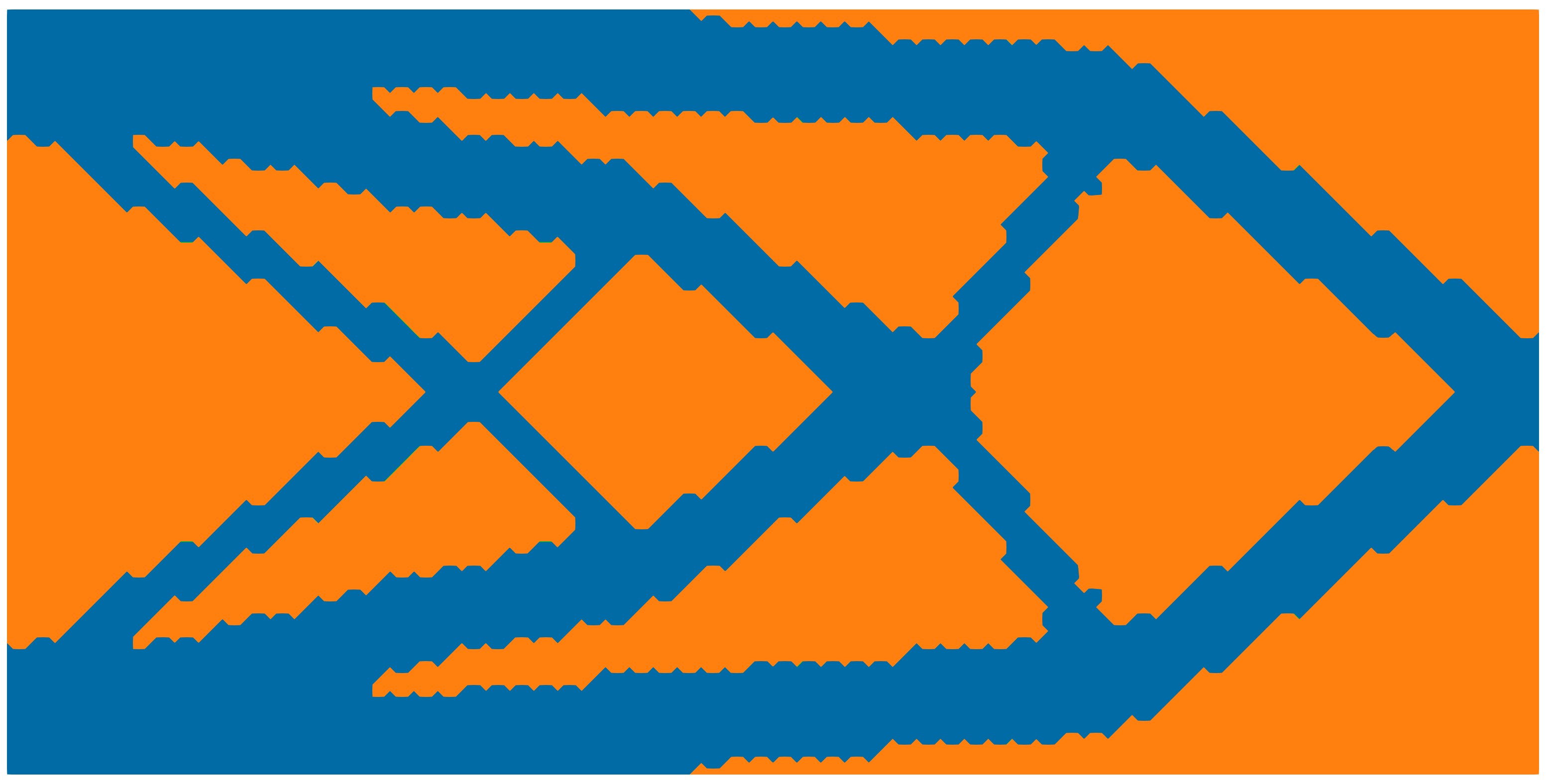}
		\caption{Sphere combination.}
	\end{subfigure}%
	\begin{subfigure}{0.25\textwidth}
		\centering
		\includegraphics[width=0.95\textwidth]{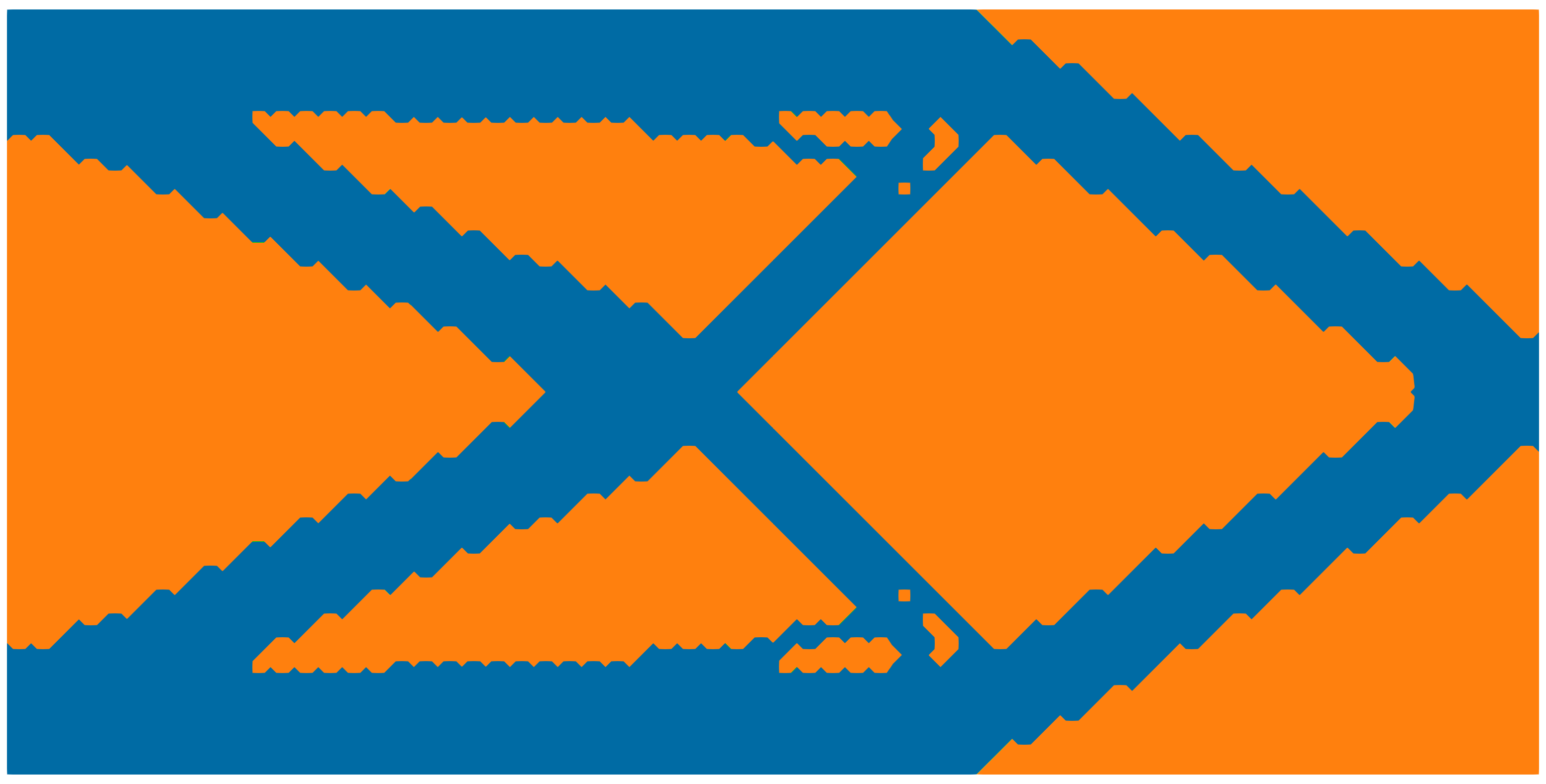}
		\caption{Convex combination.}
	\end{subfigure}%
	\begin{subfigure}{0.25\textwidth}
		\centering
		\includegraphics[width=0.95\textwidth]{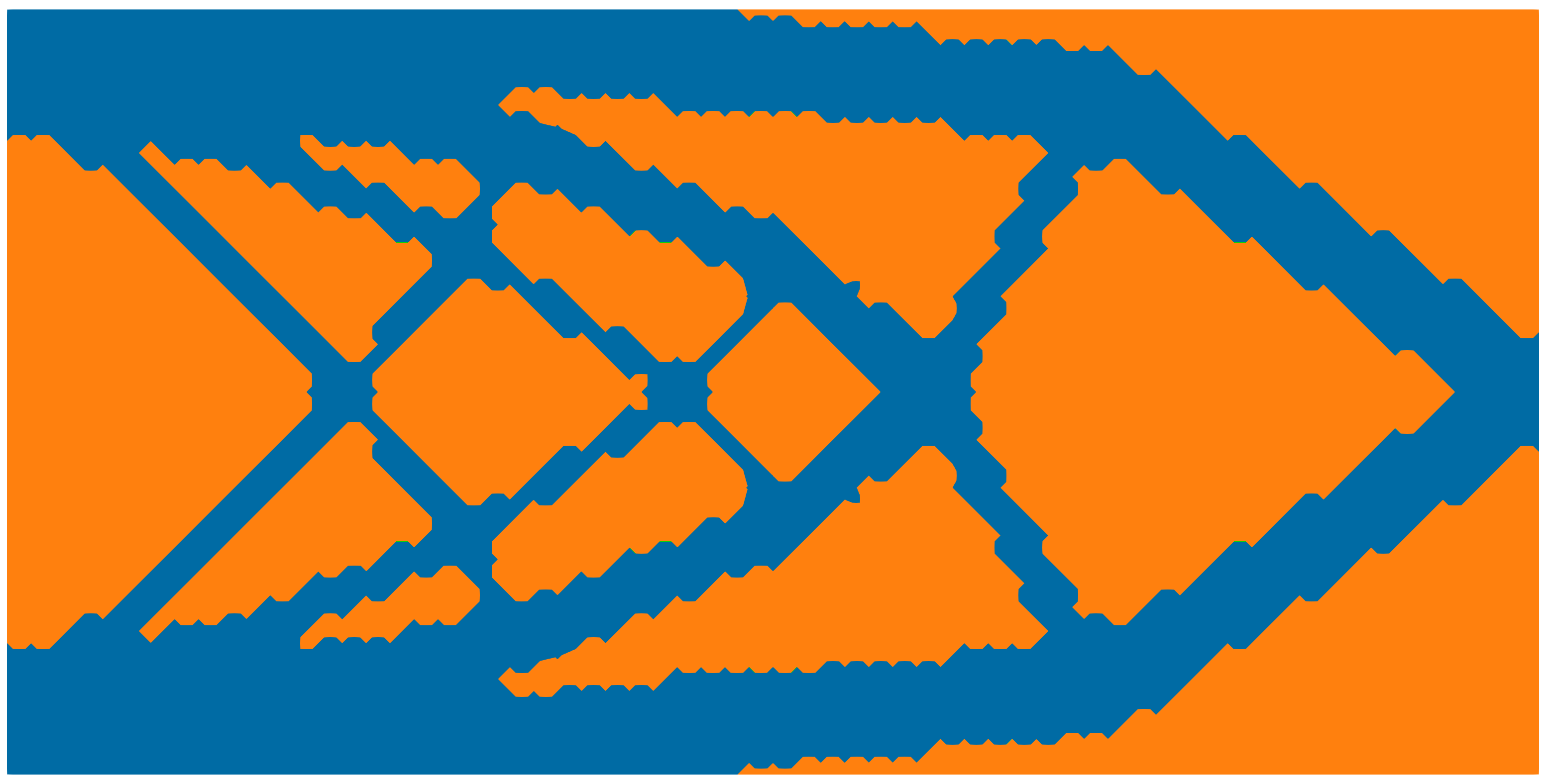}
		\caption{Gradient descent.}
	\end{subfigure}%
	\begin{subfigure}{0.25\textwidth}
		\centering
		\includegraphics[width=0.95\textwidth]{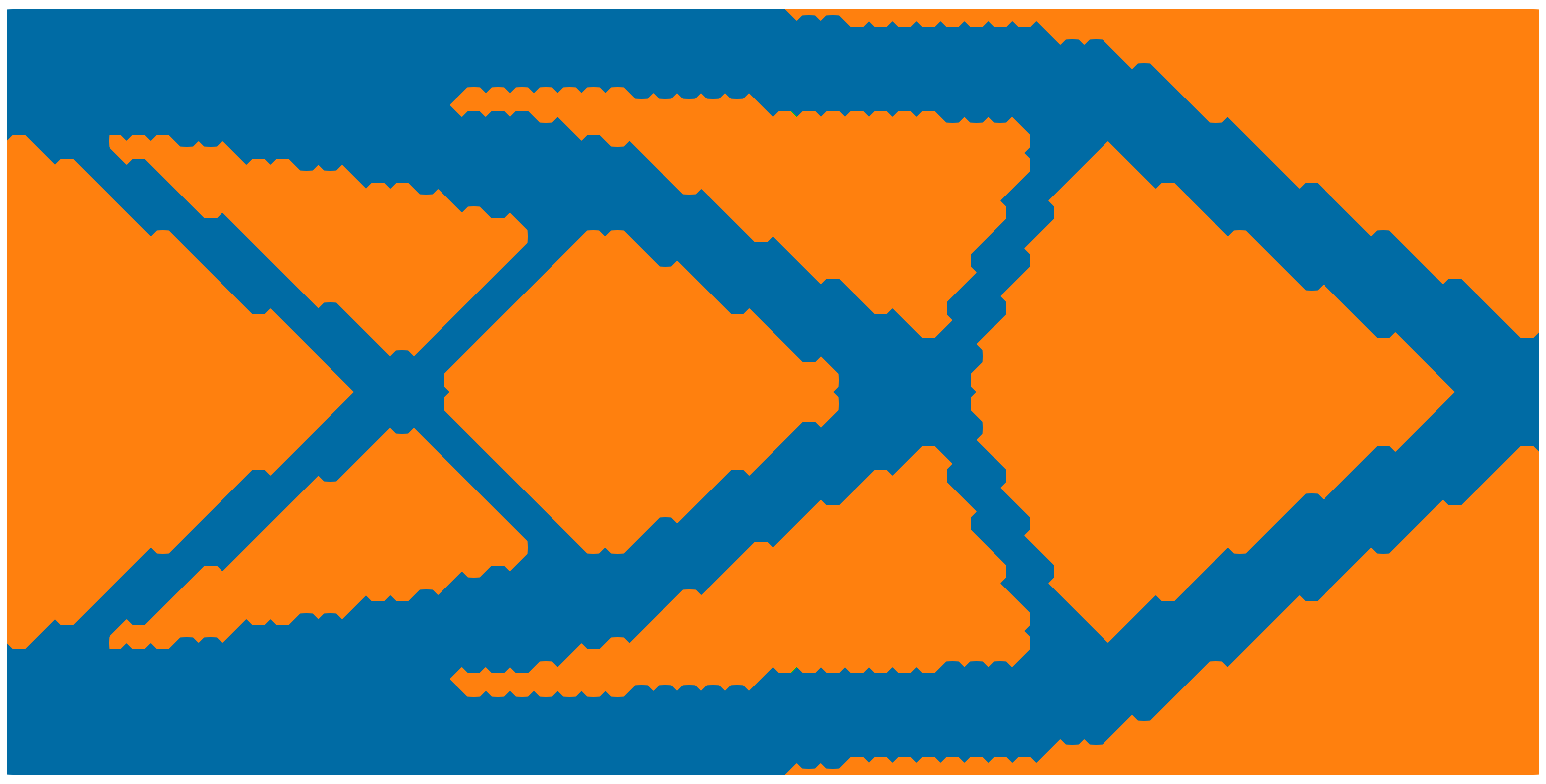}
		\caption{BFGS.}
	\end{subfigure}%
	\caption{Optimized geometries for the cantilever problem.}
	\label{fig:geometries_cantilever}
\end{figure}

For our first example, we consider the so-called cantilever problem (see, e.g., \cite{Amstutz2006new, Allaire2004Structural, Allaire2005Structural}). Here, the holdall domain is given by $\Dsf = (0,2) \times (0,1)$. The Dirichlet boundary $\Gamma_D = \Set{x = 0}$ is the left side of the rectangle and for the Neumann load $g$ we consider a unitary point load at $(2, 0.5)$. As in \cite{Amstutz2006new}, we discretize this domain with a uniform triangular mesh consisting of \num{4193}~nodes and \num{8192}~triangles and choose $l = 100$. A schematic of the problem setting can be seen in Figure~\ref{fig:cantilever}.

The evolution of the cost functional, angle criterion, and norm of the projected topological derivative are shown in Figure~\ref{fig:evolution_cantilever}. Here, we observe that all methods converge very fast, requiring a maximum of 25 iterations to satisfy the angle criterion with a tolerance of \SI{1.5}{\degree}. The performance of all methods is very comparable for this problem and no method performs significantly better or worse than the others. 

However, when we investigate the optimized geometries obtained with the methods, there are some differences, as they converged to different local minimizers of the problems. Whereas the sphere combination and BFGS methods converged to the same solution that was reported in \cite{Amstutz2006new}, the gradient descent and convex combination method converged to different geometries with finer beam structures.

\subsubsection{Bridge}

\begin{figure}[!b]
	\centering
	\begin{subfigure}[t]{0.5\textwidth}
		\centering
		\begin{tikzpicture}[scale=3]
		\draw [line width=0.4pt] (0,0) -- (2,0) -- (2,1.2) -- (0,1.2) -- cycle;
		
		\draw [dimen] (0, -0.4) -- (2, -0.4) node {2};
		\draw [dimen] (2.2, 0) -- (2.2, 1.2) node {1.2};
		
		\draw [-latex, line width=2.5pt] (1,0) -- (1, -0.25);
		
		\draw (0, -0.1) circle (0.1);
		\draw (2, -0.1) circle (0.1);
		\draw[pattern=north west lines, pattern color=black, draw=none] (-0.2,-0.2) rectangle (0.2,-0.3);
		\draw[pattern=north west lines, pattern color=black, draw=none] (1.8,-0.2) rectangle (2.2,-0.3);
		
		\node at (1.0, 0.6) {{\LARGE $D$}};
		\end{tikzpicture}
		\caption{Single load case.}
	\end{subfigure}%
	\begin{subfigure}[t]{0.5\textwidth}
		\centering
		\begin{tikzpicture}[scale=3]
		\draw [line width=0.4pt] (0,0) -- (2,0) -- (2,1.2) -- (0,1.2) -- cycle;
		
		\draw [dimen] (0, -0.4) -- (2, -0.4) node {2};
		\draw [dimen] (2.2, 0) -- (2.2, 1.2) node {1.2};
		
		\draw [-latex, line width=2.5pt] (1,0) -- (1, -0.25);
		\draw [-latex, line width=2.5pt] (0.5,0) -- (0.5, -0.25);
		\draw [-latex, line width=2.5pt] (1.5,0) -- (1.5, -0.25);
		
		\draw (0, -0.1) circle (0.1);
		\draw (2, -0.1) circle (0.1);
		
		\draw[pattern=north west lines, pattern color=black, draw=none] (-0.2,-0.2) rectangle (0.2,-0.3);
		\draw[pattern=north west lines, pattern color=black, draw=none] (1.8,-0.2) rectangle (2.2,-0.3);
		
		\node at (1.0, 0.6) {{\LARGE $D$}};
		\end{tikzpicture}
		\caption{Multiple load case.}
	\end{subfigure}%
	\caption{Schematic setup of the bridge problem.}
	\label{fig:bridge}
\end{figure}
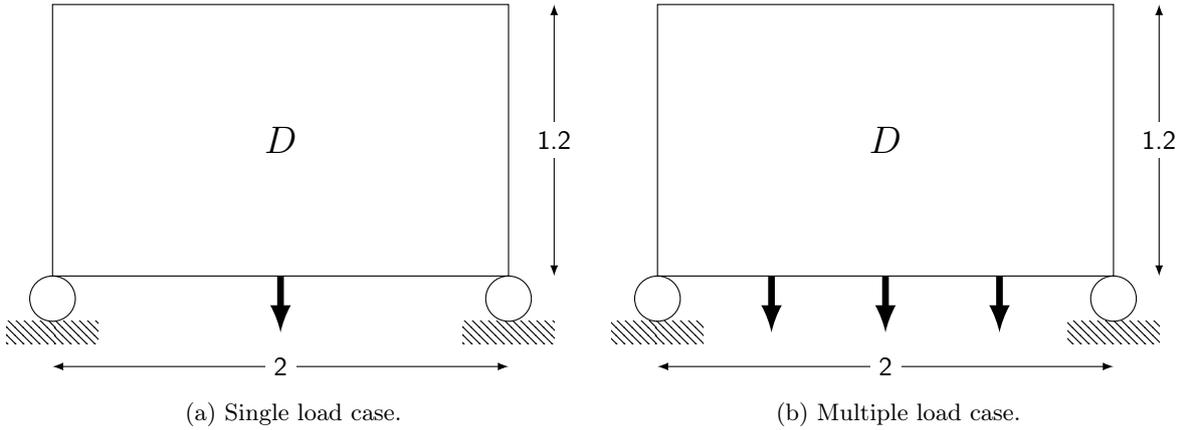

\begin{figure}[!b]
	\centering
	\begin{subfigure}{0.333\textwidth}
		\centering
		\includegraphics[width=\textwidth]{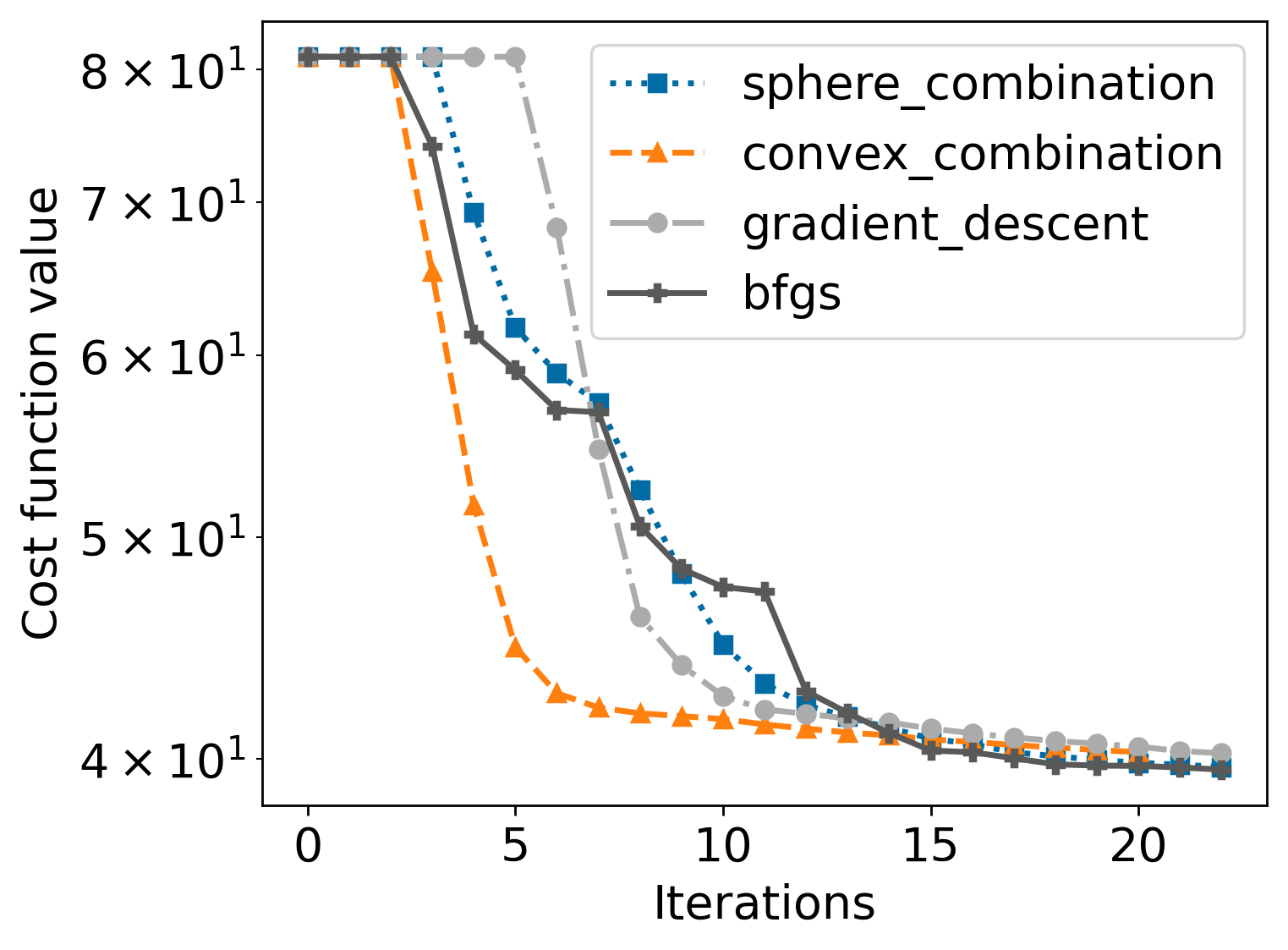}
		\caption{Cost functional.}
	\end{subfigure}%
	\begin{subfigure}{0.333\textwidth}
		\centering
		\includegraphics[width=\textwidth]{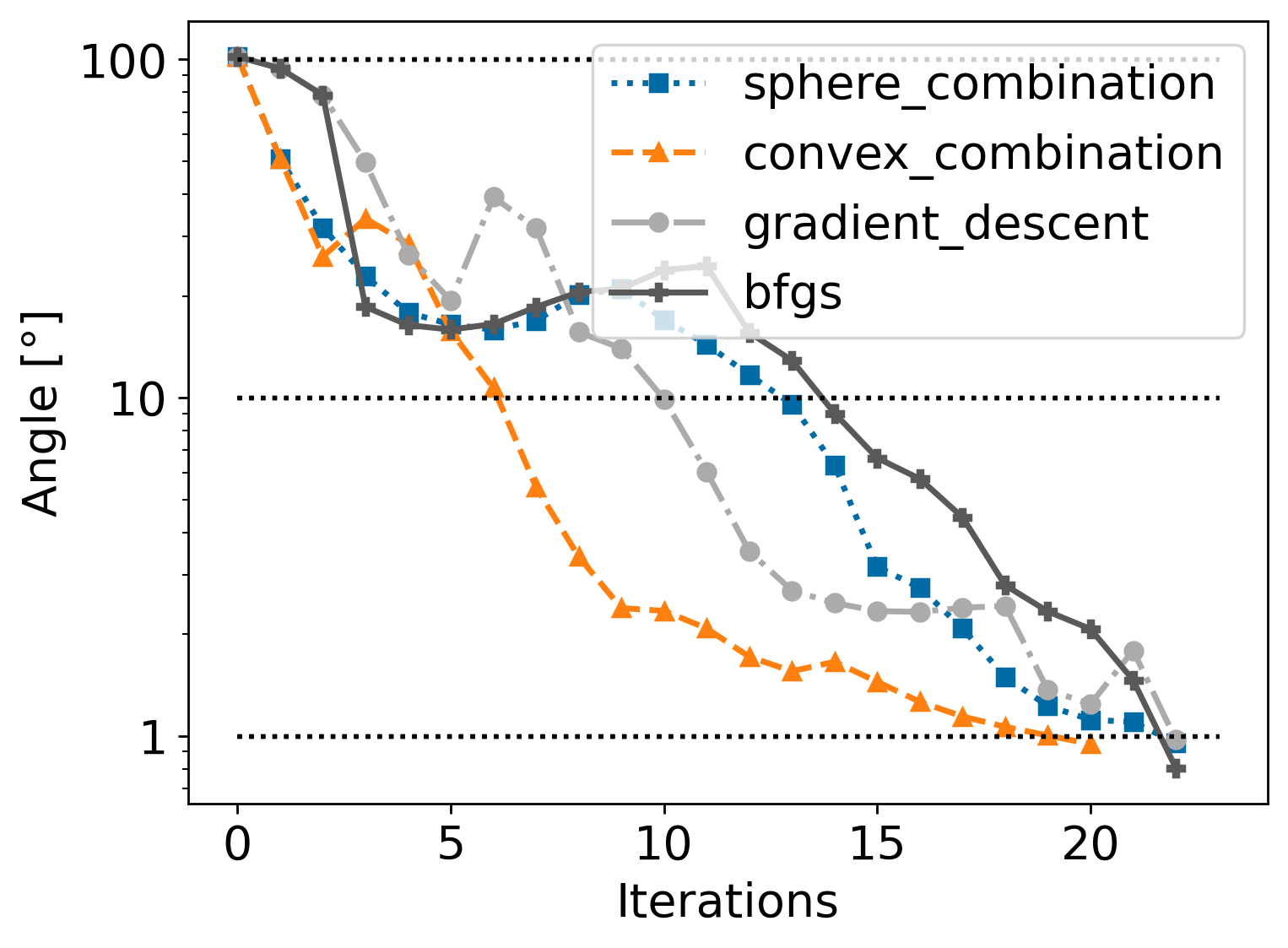}
		\caption{Angle.}
	\end{subfigure}%
	\begin{subfigure}{0.333\textwidth}
		\centering
		\includegraphics[width=\textwidth]{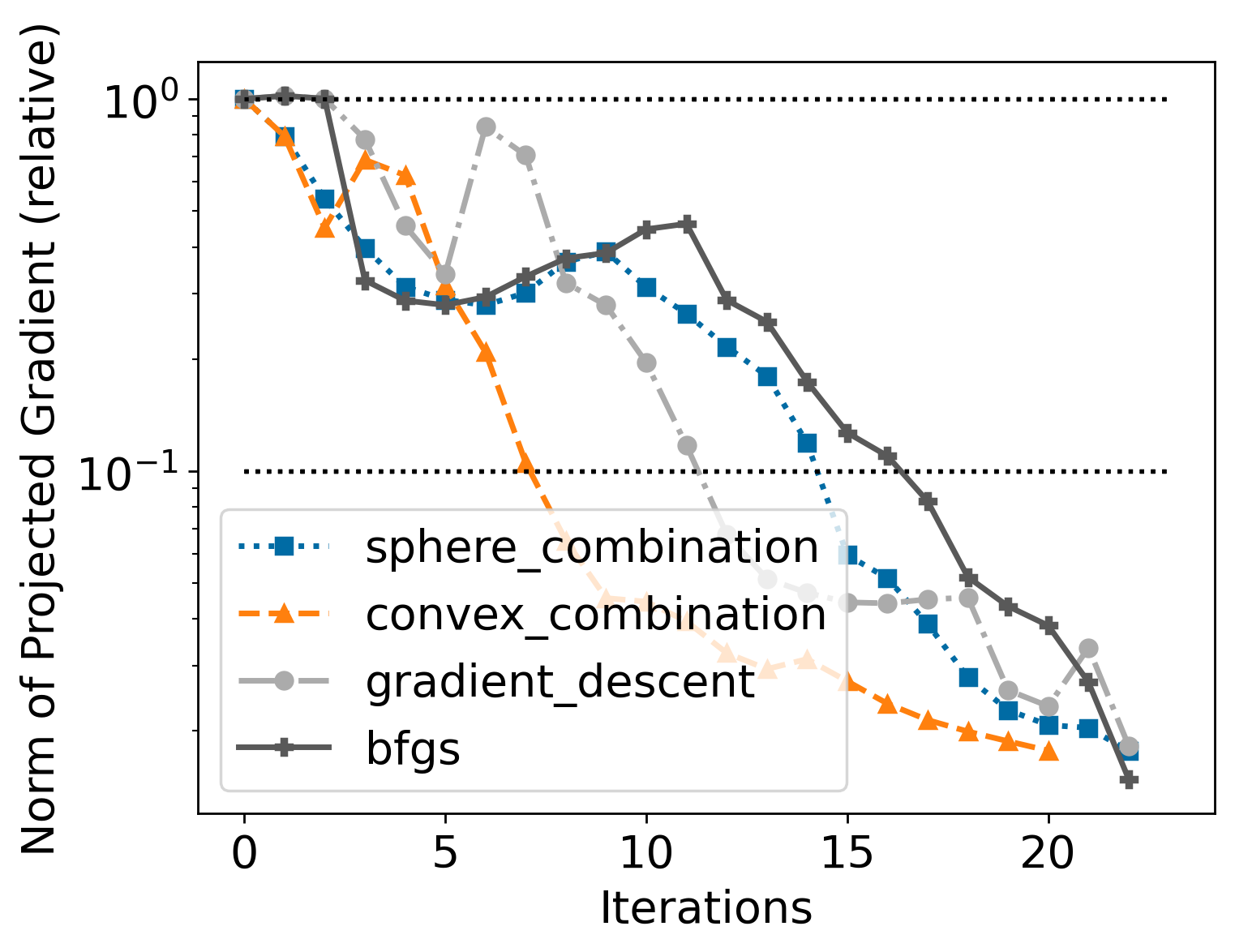}
		\caption{Norm of the projected topological derivative.}
	\end{subfigure}
	\caption{Evolution of the optimization for the bridge with a single load.}
	\label{fig:evolution_bridge_single}
\end{figure}

In this section, we consider another problem of compliance minimization in linear elasticity which corresponds to a bridge with a single and multiple loads. As before, these problems are taken from the literature \cite{Amstutz2006new, Allaire2004Structural, Allaire2005Structural}. Here, the hold-all domain is given by $\Dsf = (0,2) \times (0, 1.2)$. As boundary conditions, we have zero vertical displacement at the bottom left and right corners. For the single load case, we apply a vertical downwards unitary force at $(1, 0)$ and for the multiple load case we apply three vertical unitary loads at $(0.5, 0)$, $(1, 0)$, and $(1.5, 0)$. We discretize this setup with a mesh consisting of \num{7809}~vertices and \num{15360}~triangles. For the volume regularization, we use the parameter $l = 30$ in the single load case and use $l = 120$ for the multiple load case, in analogy to \cite{Amstutz2006new}. A schematic of the problem can be seen in Figure~\ref{fig:bridge}.

\begin{figure}[!t]
	\centering
	\begin{subfigure}{0.25\textwidth}
		\centering
		\includegraphics[width=\textwidth]{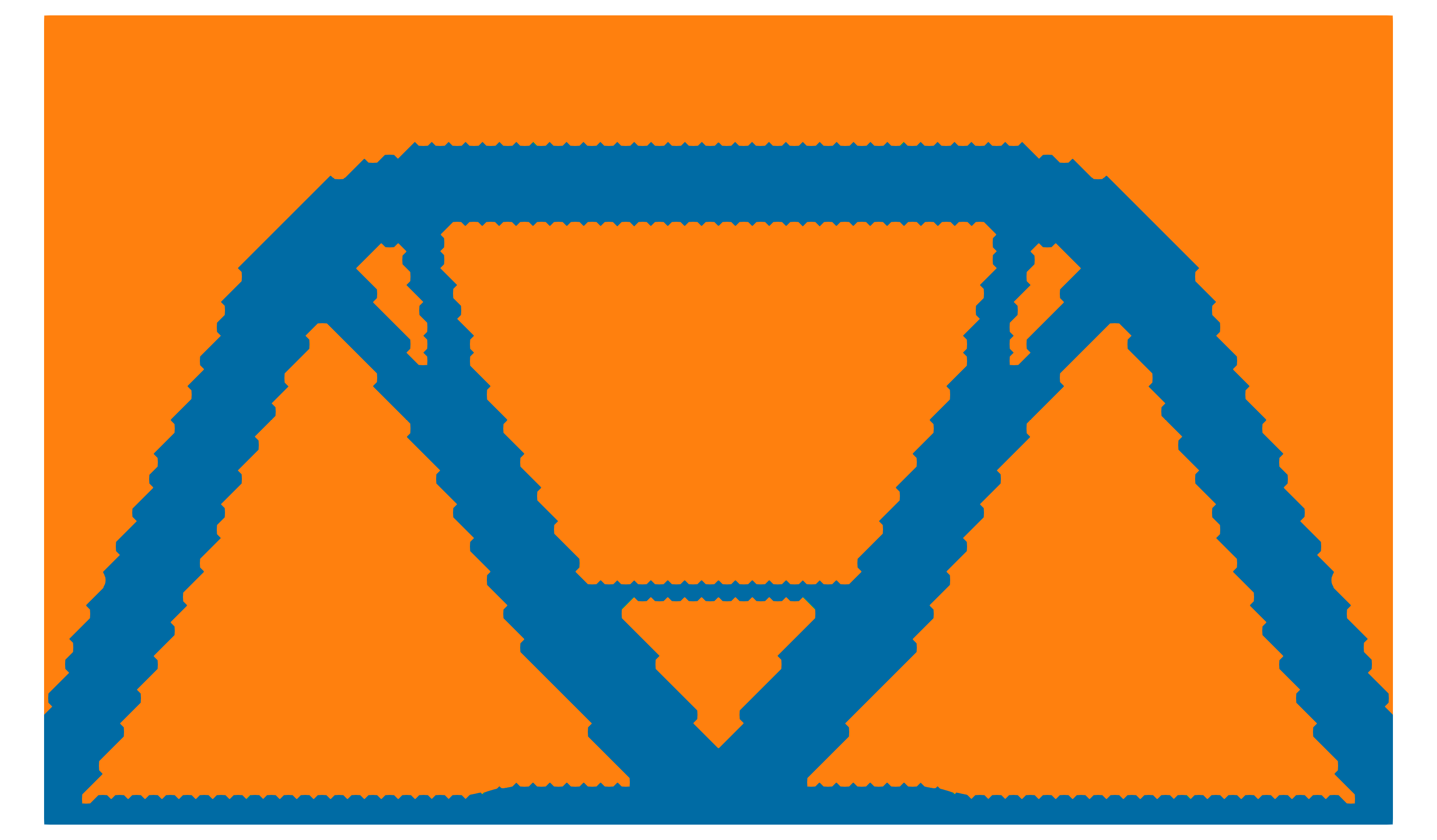}
		\caption{Sphere combination.}
	\end{subfigure}%
	\begin{subfigure}{0.25\textwidth}
		\centering
		\includegraphics[width=\textwidth]{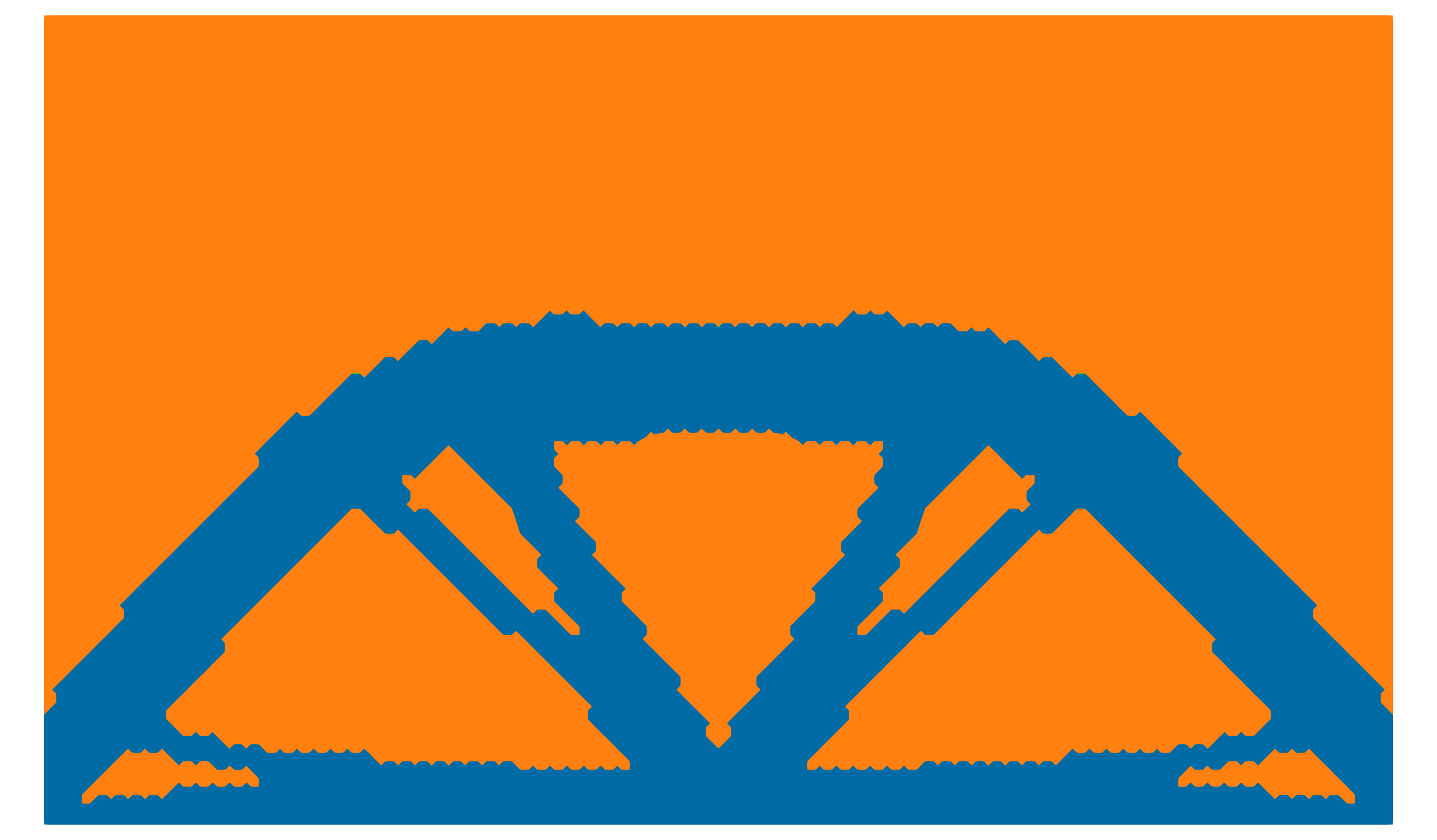}
		\caption{Convex combination.}
	\end{subfigure}%
	\begin{subfigure}{0.25\textwidth}
		\centering
		\includegraphics[width=\textwidth]{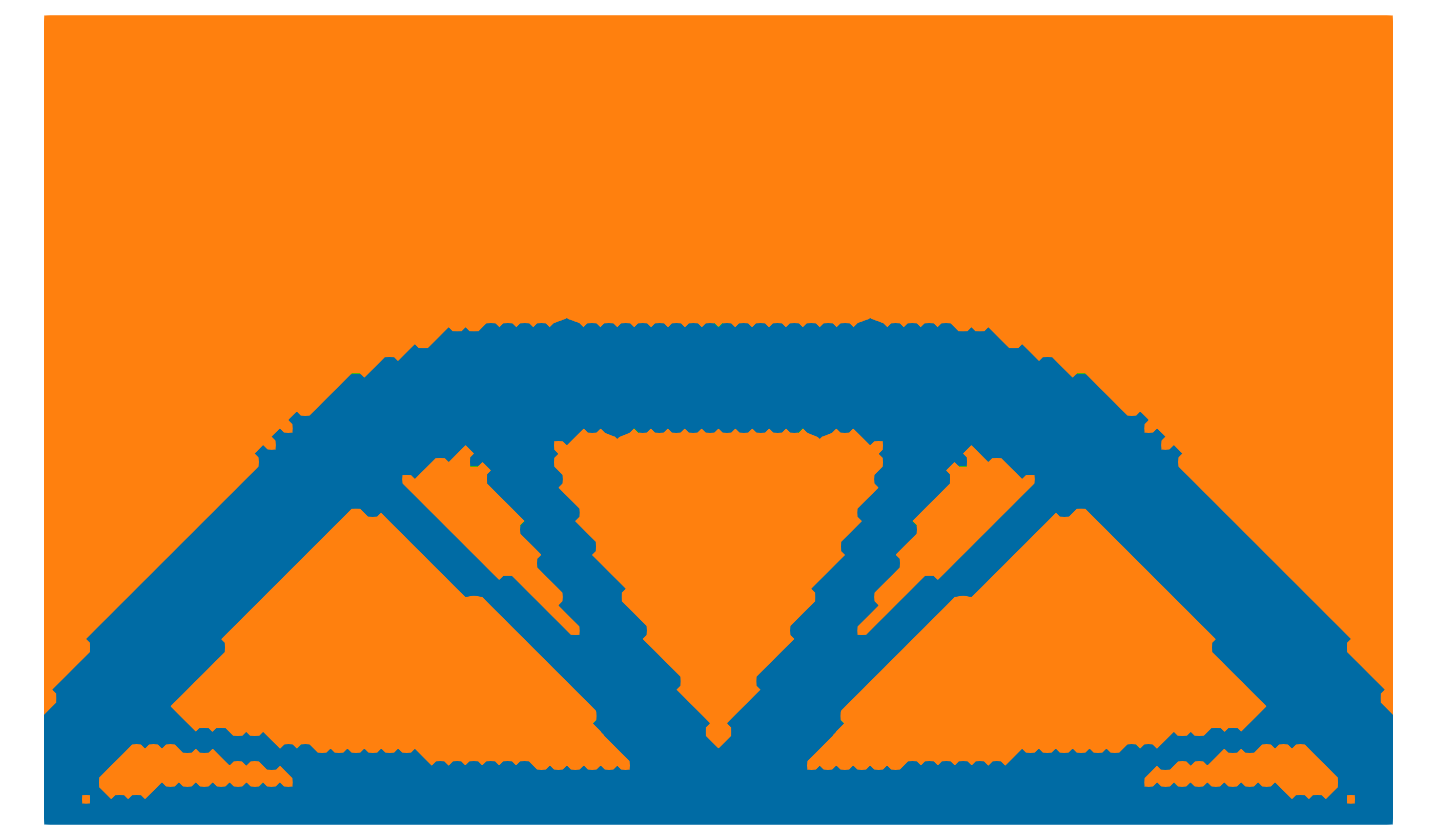}
		\caption{Gradient descent.}
	\end{subfigure}%
	\begin{subfigure}{0.25\textwidth}
		\centering
		\includegraphics[width=\textwidth]{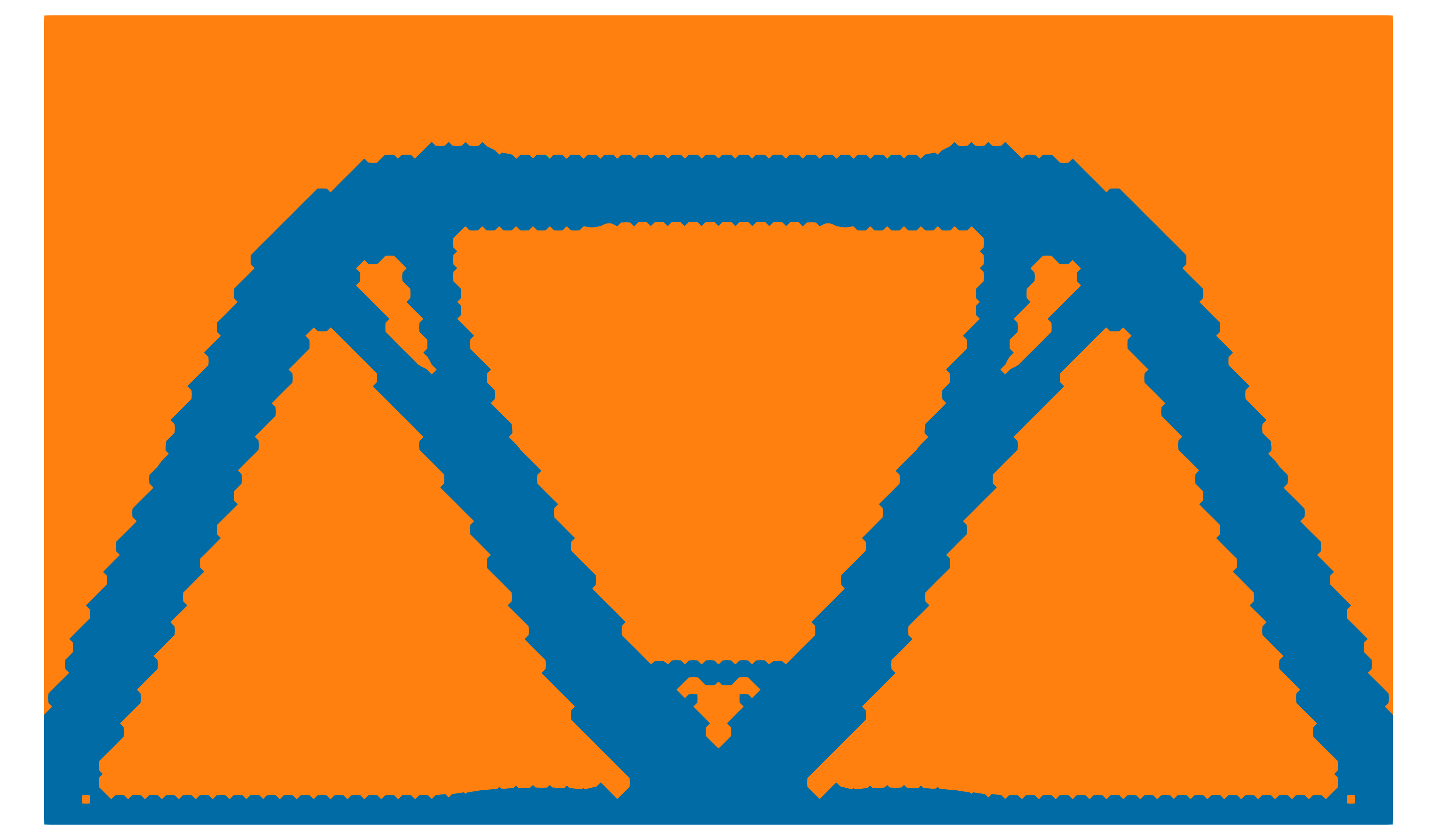}
		\caption{BFGS.}
	\end{subfigure}%
	\caption{Optimized geometries for the bridge with a single load.}
	\label{fig:geom_bridge_single}
\end{figure}

\begin{figure}[!b]
	\centering
	\begin{subfigure}{0.333\textwidth}
		\centering
		\includegraphics[width=\textwidth]{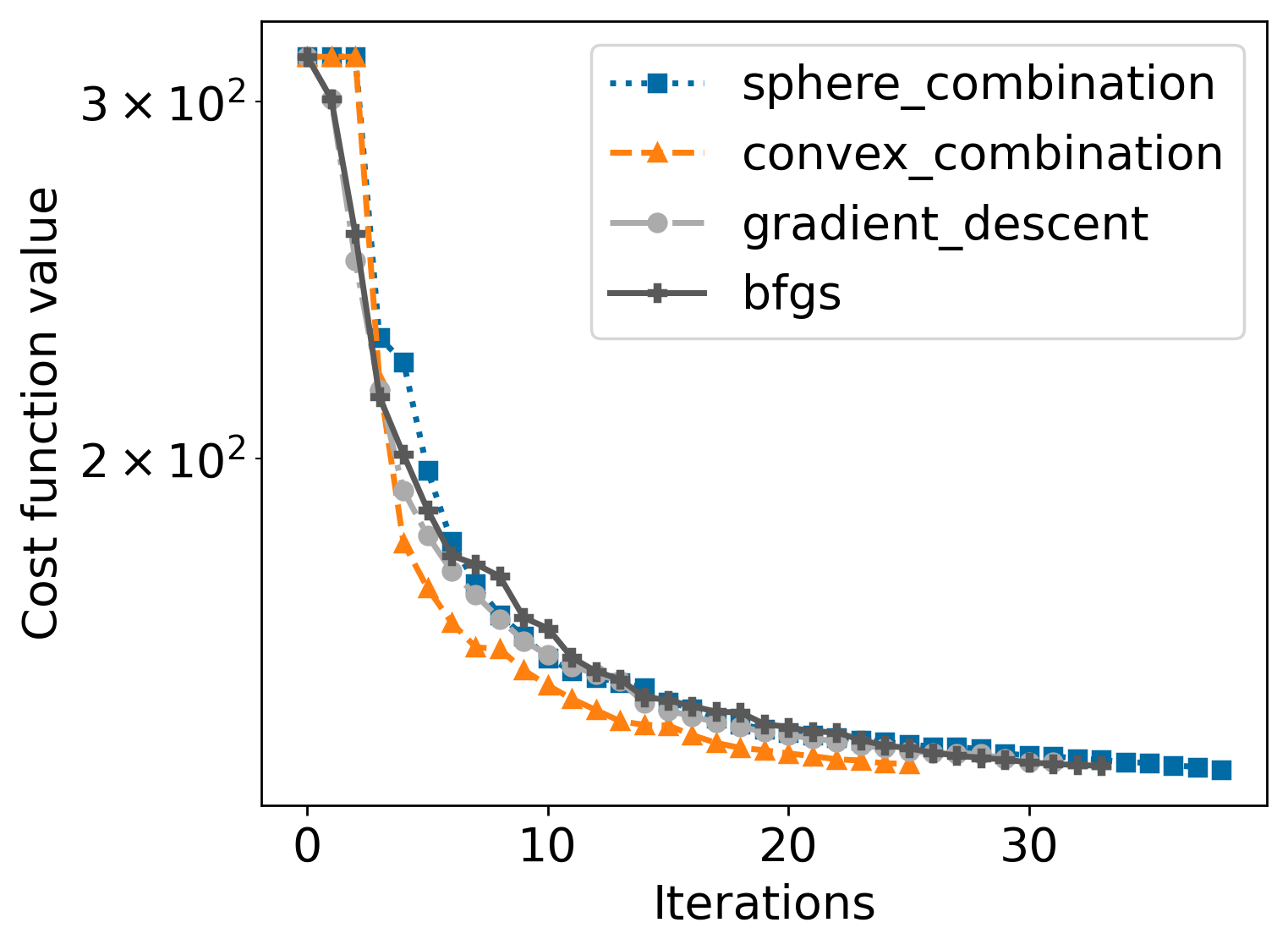}
		\caption{Cost functional.}
	\end{subfigure}%
	\begin{subfigure}{0.333\textwidth}
		\centering
		\includegraphics[width=\textwidth]{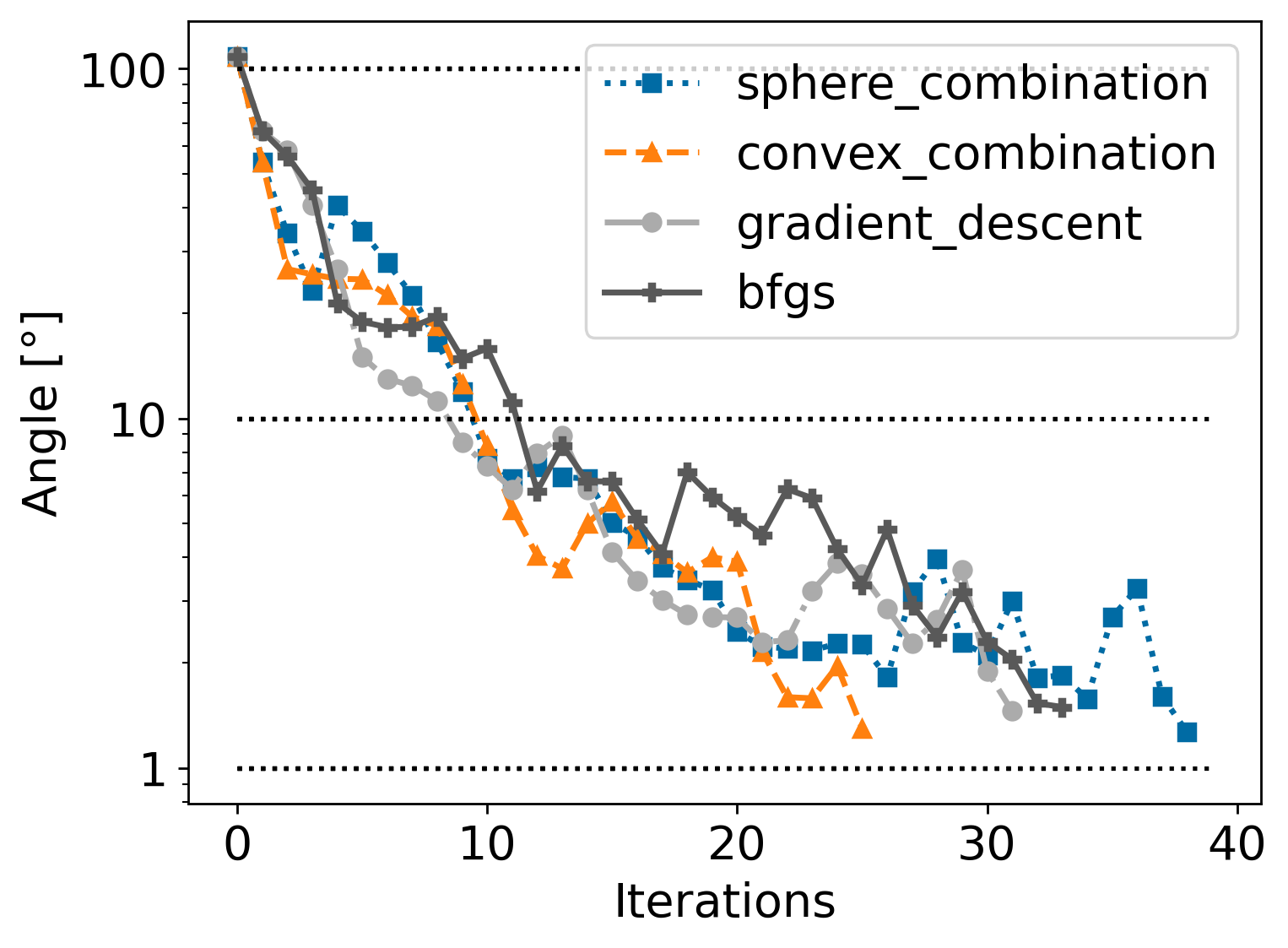}
		\caption{Angle.}
	\end{subfigure}%
	\begin{subfigure}{0.333\textwidth}
		\centering
		\includegraphics[width=\textwidth]{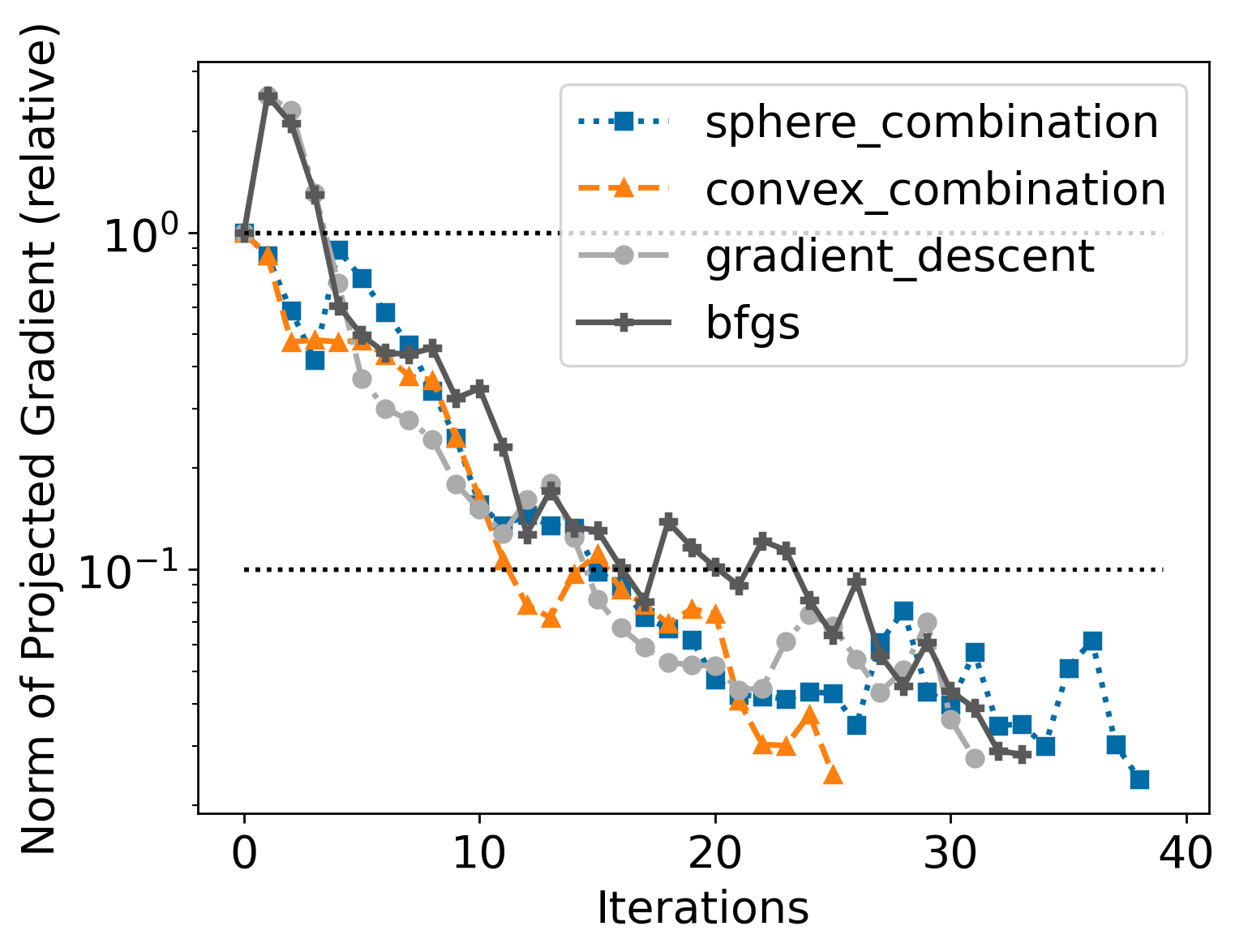}
		\caption{Norm of the projected topological derivative.}
	\end{subfigure}
	\caption{Evolution of the optimization for the bridge with multiple loads.}
	\label{fig:evolution_bridge_multiple}
\end{figure}

\begin{figure}[!b]
	\centering
	\begin{subfigure}{0.25\textwidth}
		\centering
		\includegraphics[width=\textwidth]{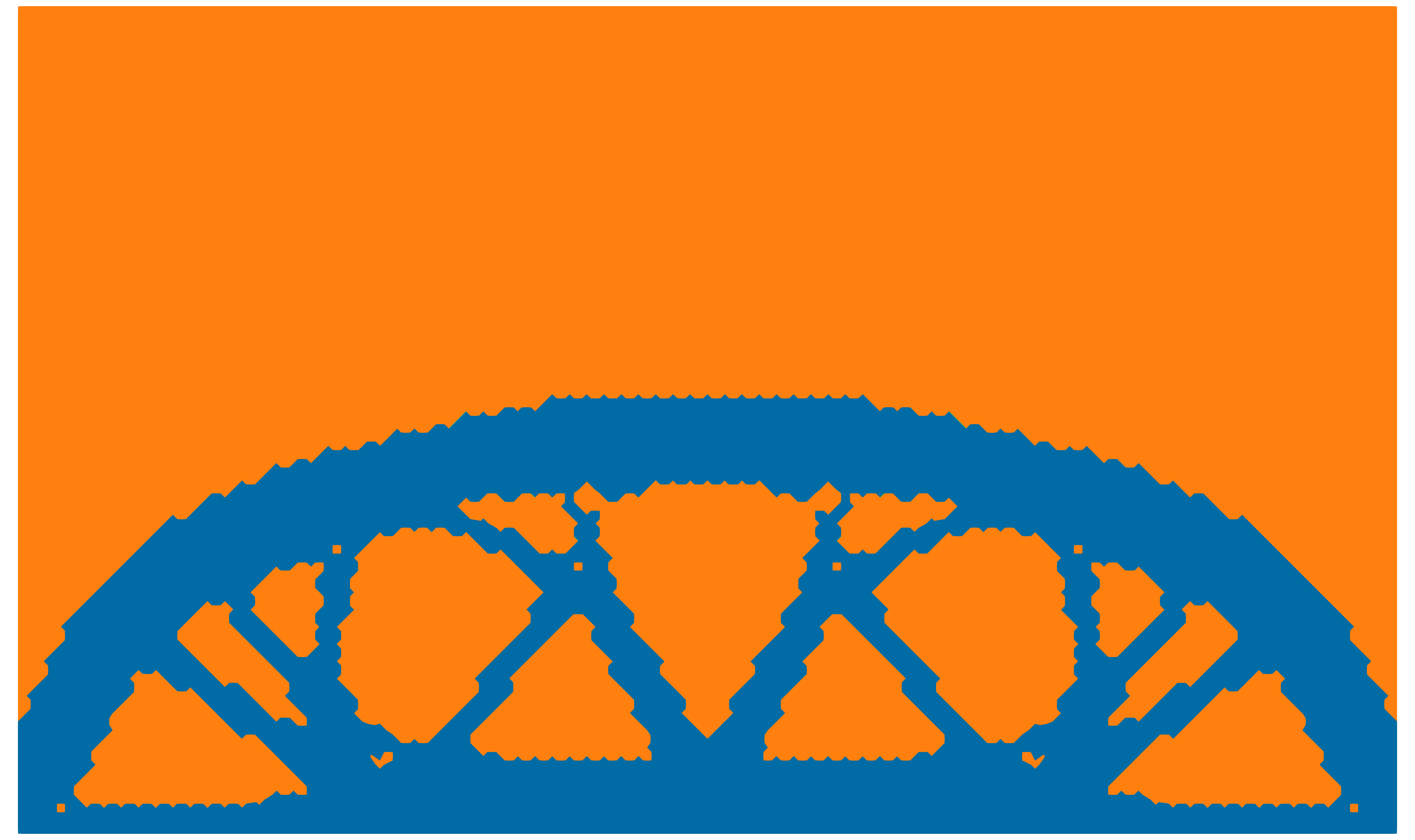}
		\caption{Sphere combination.}
	\end{subfigure}%
	\begin{subfigure}{0.25\textwidth}
		\centering
		\includegraphics[width=\textwidth]{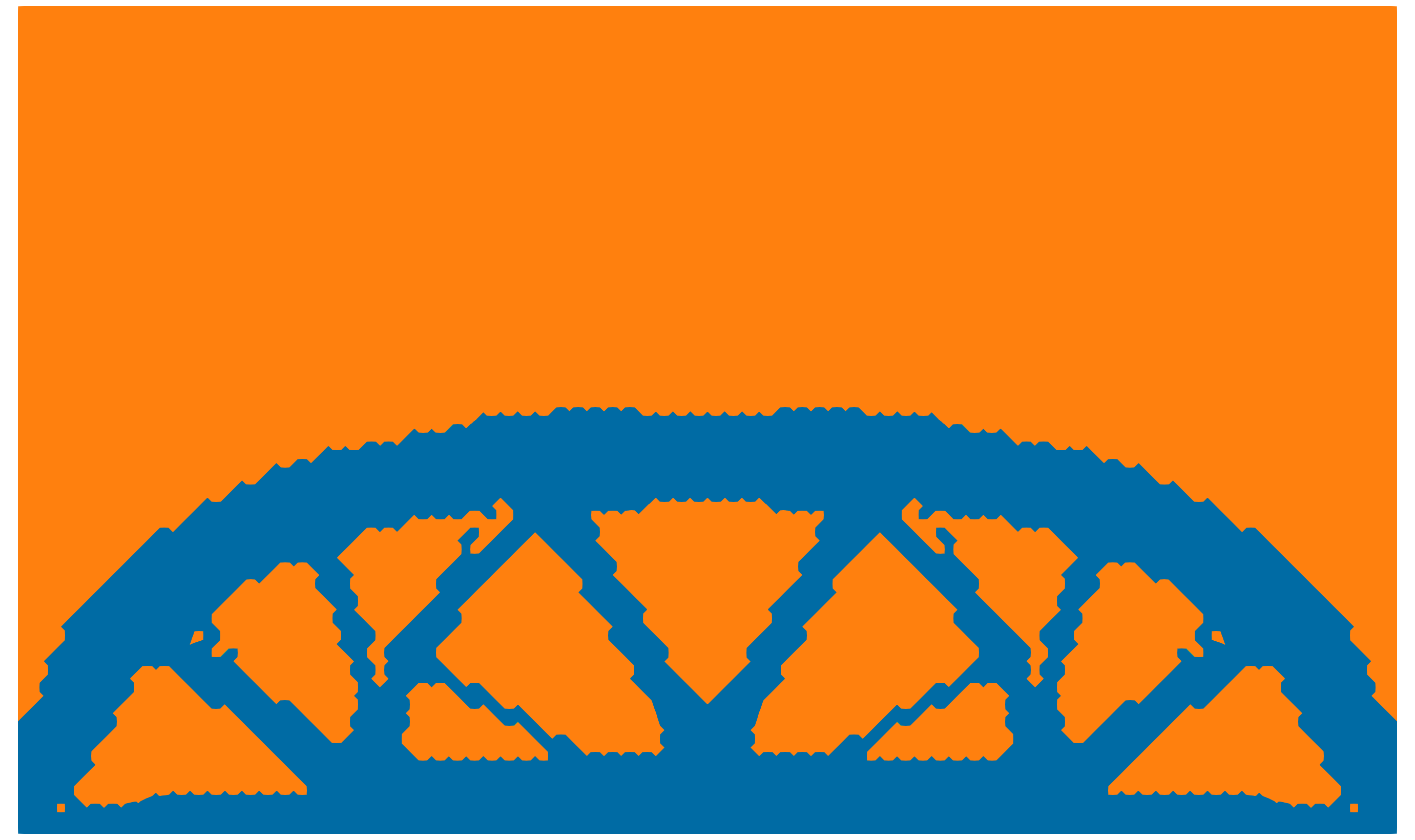}
		\caption{Convex combination.}
	\end{subfigure}%
	\begin{subfigure}{0.25\textwidth}
		\centering
		\includegraphics[width=\textwidth]{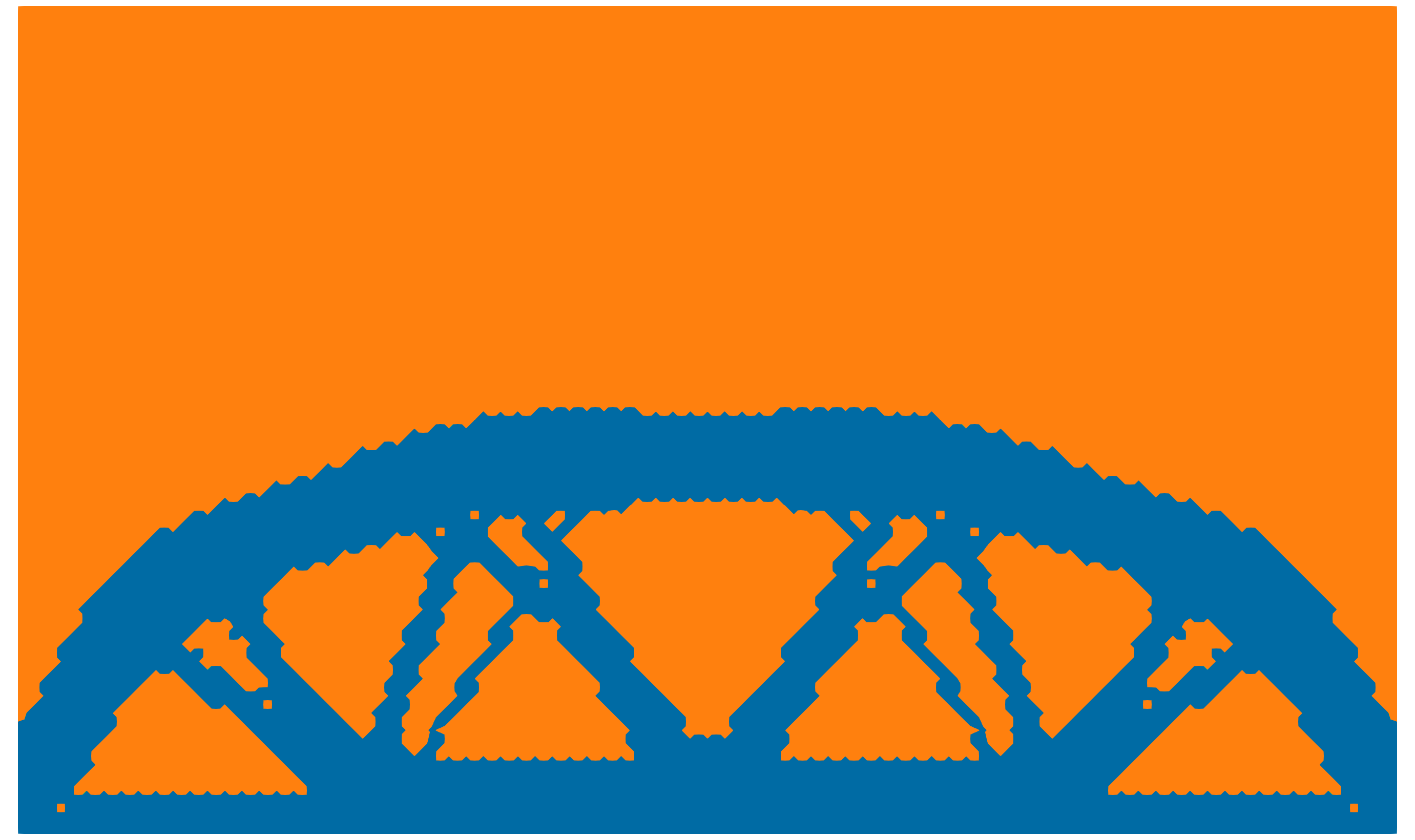}
		\caption{Gradient descent.}
	\end{subfigure}%
	\begin{subfigure}{0.25\textwidth}
		\centering
		\includegraphics[width=\textwidth]{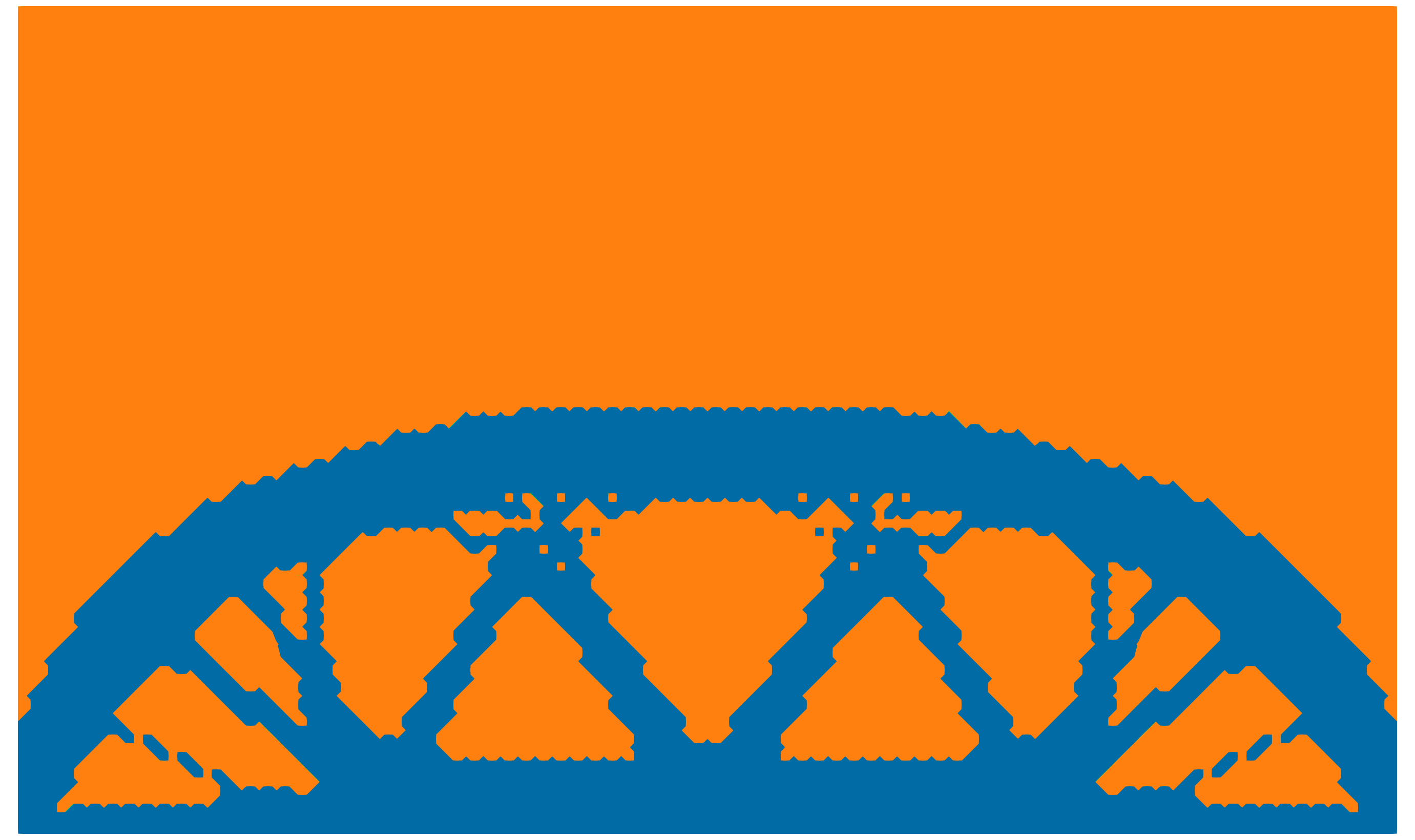}
		\caption{BFGS.}
	\end{subfigure}%
	\caption{Optimized geometries for the bridge with multiple load.}
	\label{fig:geom_bridge_multiple}
\end{figure}

Let us first discuss the single load case. The history of the cost functional, angle criterion, and norm of the projected topological derivative are shown in Figure~\ref{fig:evolution_bridge_single}. Similarly to the cantilever problem, there are no significant differences in the performance of the optimization algorithms. All methods required between 20 and 22 iterations to reach the prescribed tolerance of \SI{1}{\degree} of the angle criterion. The only notable difference between the methods is that the convex combination algorithm decreases the quality measures slightly earlier than the other methods, but they \qe{catch up} during later iterations. The optimized geometries obtained with the methods are shown in Figure~\ref{fig:geom_bridge_single}. Here, we also observe, that the methods converge to different local minimizers. Whereas the sphere combination and BFGS method find a similar geometry, which is the same as reported in \cite{Amstutz2006new}, the convex combination and gradient descent method find a different local minimizer with less height and a different supporting beam structure.

The results are very similar for the multiple loads case. In Figure~\ref{fig:evolution_bridge_multiple} the evolution of the considered quality measures is depicted over the history of the optimization. Here, we observe some slight differences in the performance of the algorithms. The convex combination method performs best as it requires only 25 iterations to satisfy the angle criterion, where we have chosen a tolerance of \SI{1.5}{\degree} for this example. The gradient descent and BFGS methods showed the second best performance, requiring about 30 iterations each to satisfy the stopping criterion. The sphere combination method performed worst and required 40 iterations to reach the stopping tolerance. Finally, let us briefly investigate the optimized geometries, which are shown in Figure~\ref{fig:geom_bridge_multiple}. Here, we see that all methods converge to a similar solution, which has a slightly more complicated beam structure than the one reported in \cite{Amstutz2006new}.

Altogether, for the case of linear elasticity, we observe no major differences between all considered optimization algorithms. The proposed BFGS methods show a very similar performance to the already established methods on all considered test problems. A possible reason for this behavior is that the problems are comparatively easy to solve, at least in comparison to the problems considered in Sections~\ref{ssec:numerics_model} and~\ref{ssec:semilinear_transmission}, as the already established methods only require 20 to 50 iterations to solve these problems to a desired tolerance, so that there may not be enough iterations for the BFGS methods to show their potential, particularly if they have to be restarted often due to large changes in the topology. This topic is of interest for future research.

\subsection{Optimization of Fluids in Navier-Stokes Flow}
\label{ssec:stokes}

Let us now consider another application, the optimization of fluids in Navier-Stokes flow. Our setting is similar to before, i.e., let $\Dsf \subset \R^d$ be an open and bounded hold-all domain with boundary $\partial \Dsf$. Let $\Omega \subset \Dsf$, denote by $\Omega^c = \Dsf \setminus \closure{\Omega}$ the complement of $\Omega$ in $\Dsf$, and let $\Gamma = \partial \Omega \setminus \partial \Dsf$. We consider the following optimization problem
\begin{equation}
\label{eq:topo_stokes}
\begin{aligned}
&\min_{\Omega} J(\Omega, u) = \int_\Dsf \mu \nabla u : \nabla u + \alpha_\Omega u \cdot u \dmeas{x} \\
&\text{s.t.} \quad \begin{alignedat}[t]{2}
- \mu \Delta u + \rho \left( u \cdot \nabla \right) u + \nabla p + \alpha_\Omega u &= 0 \quad &&\text{ in } \Dsf, \\
\nabla \cdot u &= 0 \quad &&\text{ in } \Dsf,\\
u &= u_\mathrm{D} \quad &&\text{ on } \partial \Dsf,\\
\int_\Dsf p \dmeas{x} &= 0, \\
\abs{\Omega} &= \text{vol}_\mathrm{des},
\end{alignedat}
\end{aligned}
\end{equation}
where $\alpha_\Omega(x) = \chi_\Omega(x) \alpha_\mathrm{in} + \chi_{\Omega^c}(x) \alpha_\mathrm{out}$ with $\alpha_\mathrm{in}, \alpha_\mathrm{out} > 0$ and $\abs{\Omega}$ denotes the Lebesgue measure in $\mathbb{R}^d$. Here, $u$ denotes the fluid's velocity and $p$ its pressure, $\mu$ is its viscosity, $\rho$ its density, and $\alpha$ is the inverse permeability. The cost functional of the above problem models the energy dissipation of the fluid, which should be minimized. Moreover, we have a volume equality constraint, which ensures that only the desired volume of the domain is occupied by the fluid. For the topology optimization problem \eqref{eq:topo_stokes}, the sought domain $\Omega$ is the domain of the fluid, whereas its complement $\Omega^c$ plays the role of a  solid region. In particular, the inverse permeability $\alpha$ is small inside $\Omega$ and large outside of it. For a more detailed discussion, we refer the reader to \cite{Borrvall2003Topology}, where this model was first introduced and used for the topology optimization of fluids. 

For our numerical experiments, we regularize the equality constraint with a quadratic penalty term, leading to the following optimization problem
\begin{equation}
\label{eq:topo_stokes_reg}
\begin{aligned}
&\min_{\Omega} J(\Omega, u) = \int_\Dsf \mu \nabla u : \nabla u + \alpha_\Omega u\cdot u \dmeas{x} + \frac{l}{2} \left( \left|\Omega\right| - \text{vol}_\mathrm{des} \right)^2 \\
&\text{s.t.} \quad \begin{alignedat}[t]{2}
- \mu \Delta u + \rho \left( u \cdot \nabla \right) u + \nabla p + \alpha_\Omega u &= 0 \quad &&\text{ in } \Dsf, \\
\nabla \cdot u &= 0 \quad &&\text{ in } \Dsf,\\
u &= u_\mathrm{D} \quad &&\text{ on } \partial \Dsf,\\
\int_\Dsf p \dmeas{x} &= 0.
\end{alignedat}
\end{aligned}
\end{equation}
The topological derivative for problem \eqref{eq:topo_stokes_reg} can be found, e.g., in \cite{N.Sa2016Topological} and is given by
\begin{equation*}
\mathcal{D}J(\Omega)(x) = \left( \alpha_\mathrm{out} - \alpha_\mathrm{in} \right) u(x) \cdot \left( u(x) + v(x) \right) + l \left( \left|\Omega\right| - \text{vol}_\mathrm{des} \right) \quad \text{ for } x \in \Dsf \setminus \Gamma,
\end{equation*}
where $(v, q)$ solves the adjoint Navier-Stokes system
\begin{equation*}
\begin{alignedat}{2}
-\mu \Delta v + \rho (Du)\transposed v - \rho \left( u \cdot \nabla \right) v + \nabla q + \alpha_\Omega v &= 2\left( \mu \Delta u - \alpha_\Omega u\right) \quad &&\text{ in } \Dsf, \\
\text{div}(v) &= 0 \quad &&\text{ in } \Dsf, \\
v &= 0 \quad &&\text{ on } \partial \Dsf,\\
\int_\Dsf q \dmeas{x} &= 0.
\end{alignedat}
\end{equation*}

To model an actual solid region, $\alpha_\mathrm{out}$ should tend to $+\infty$. For our numerical investigation, however, we follow \cite{N.Sa2016Topological} and chose a finite value for $\alpha_\mathrm{in}$ and $\alpha_\mathrm{out}$, namely we choose
\begin{equation*}
\alpha_\mathrm{in} = \frac{2.5 \mu}{100^2} \qquad \alpha_\mathrm{out} = \frac{2.5 \mu}{0.01^2}.
\end{equation*}
Moreover, we choose a value of $\rho = 1$ for the fluid's density and $\mu = \num{1e-2}$ for its viscosity. In the following, we will consider two common benchmark problems, which are taken from \cite{Borrvall2003Topology,N.Sa2016Topological} and consider the problem of constructing a pipe bend and the drag minimization of an obstacle. Note, that for both problems, the hold-all domain is given by $\Dsf = (0,1) \times (0,1)$, and that we discretize this with a uniform mesh consisting of \num{20201}~nodes and \num{40000}~triangles. Moreover, we use the LBB-stable Taylor-Hood finite element pair of quadratic Lagrange elements for the velocity and adjoint velocity and linear Lagrange elements for the pressure and adjoint pressure.

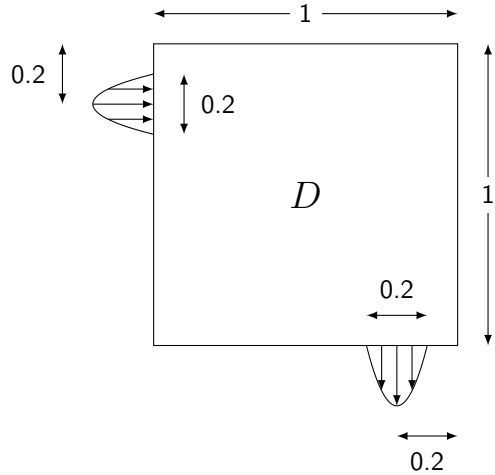
\begin{figure}[!b]
	\centering
	\begin{tikzpicture}[scale=4]
	\draw [line width=0.4pt] (0,0) -- (1,0) -- (1,1) -- (0,1) -- cycle;
	
	\draw [dimen] (0, 1.1) -- (1, 1.1) node {1};
	\draw [dimen] (1.1, 0) -- (1.1, 1) node {1};
	\draw [dimen] (0.1, 0.7) -- (0.1, 0.9) node [right=0.1] {0.2};
	\draw [dimen] (0.7, 0.1) -- (0.9, 0.1) node [above=0.1] {0.2};
	\draw [dimen] (-0.3, 0.8) -- (-0.3, 1) node [left=0.1] {0.2};
	\draw [dimen] (0.8, -0.3) -- (1, -0.3) node [below=0.1] {0.2};
	
	\draw[domain=0.7:0.9, smooth, variable=\y] plot (-{0.2/0.01*(\y-0.7)*(0.9-\y)}, {\y});
	
	\draw[-latex] (-0.15, 0.75) -- (0, 0.75);
	\draw[-latex] (-0.2, 0.8) -- (0, 0.8);
	\draw[-latex] (-0.15, 0.85) -- (0, 0.85);
	
	\draw[domain=0.7:0.9, smooth, variable=\x] plot ({\x}, -{0.2/0.01*(\x-0.7)*(0.9-\x)});
	
	\draw[-latex] (0.75, 0) -- (0.75, -0.15);
	\draw[-latex] (0.8, 0) -- (0.8, -0.2);
	\draw[-latex] (0.85, 0) -- (0.85, -0.15);
	
	\node at (0.5, 0.5) {{\LARGE $D$}};
	\end{tikzpicture}
	\caption{Schematic setup of the pipe bend problem.}
	\label{fig:pipe_bend}
\end{figure}

\subsubsection{Pipe Bend}

Let us first investigate the problem of designing a pipe bend, which is taken from \cite{Borrvall2003Topology}. For Dirichlet boundary conditions, we prescribe the inlet velocity with the parabolic profile 
\begin{equation*}
u_\mathrm{D}(x) = \begin{bsmallmatrix}
1 - 100 \left(x_2 - 0.8\right)^2, \\
0
\end{bsmallmatrix} \qquad \text{ for } x_1 = 0 \text{ and } 0.7 \leq x_2 \leq 0.9
\end{equation*}
on the top part of the left boundary of $\Dsf$, and for the outlet velocity we use the profile
\begin{equation*}
u_\mathrm{D}(x) = \begin{bsmallmatrix}
0, \\
-\left( 1 - 100 \left(x_1 - 0.8\right)^2\right)
\end{bsmallmatrix} \qquad \text{ for } x_2 = 0 \text{ and } 0.7 \leq x_1 \leq 0.9
\end{equation*}
on the right side of the bottom boundary of $\Dsf$. On all other parts of the boundary, we use a no-slip boundary condition, so that $u_\mathrm{D} = [0, 0]\transposed$. The corresponding setup is shown schematically in Figure~\ref{fig:pipe_bend}. For the volume constraint, we proceed according to \cite{Borrvall2003Topology} and use a value of $\text{vol}_\mathrm{des} = 0.08 \pi$. Additionally, we choose a regularization parameter of $l = \num{1e4}$ to enforce the volume constraint.

\begin{figure}[!t]
	\centering
	\begin{subfigure}{0.333\textwidth}
		\centering
		\includegraphics[width=\textwidth]{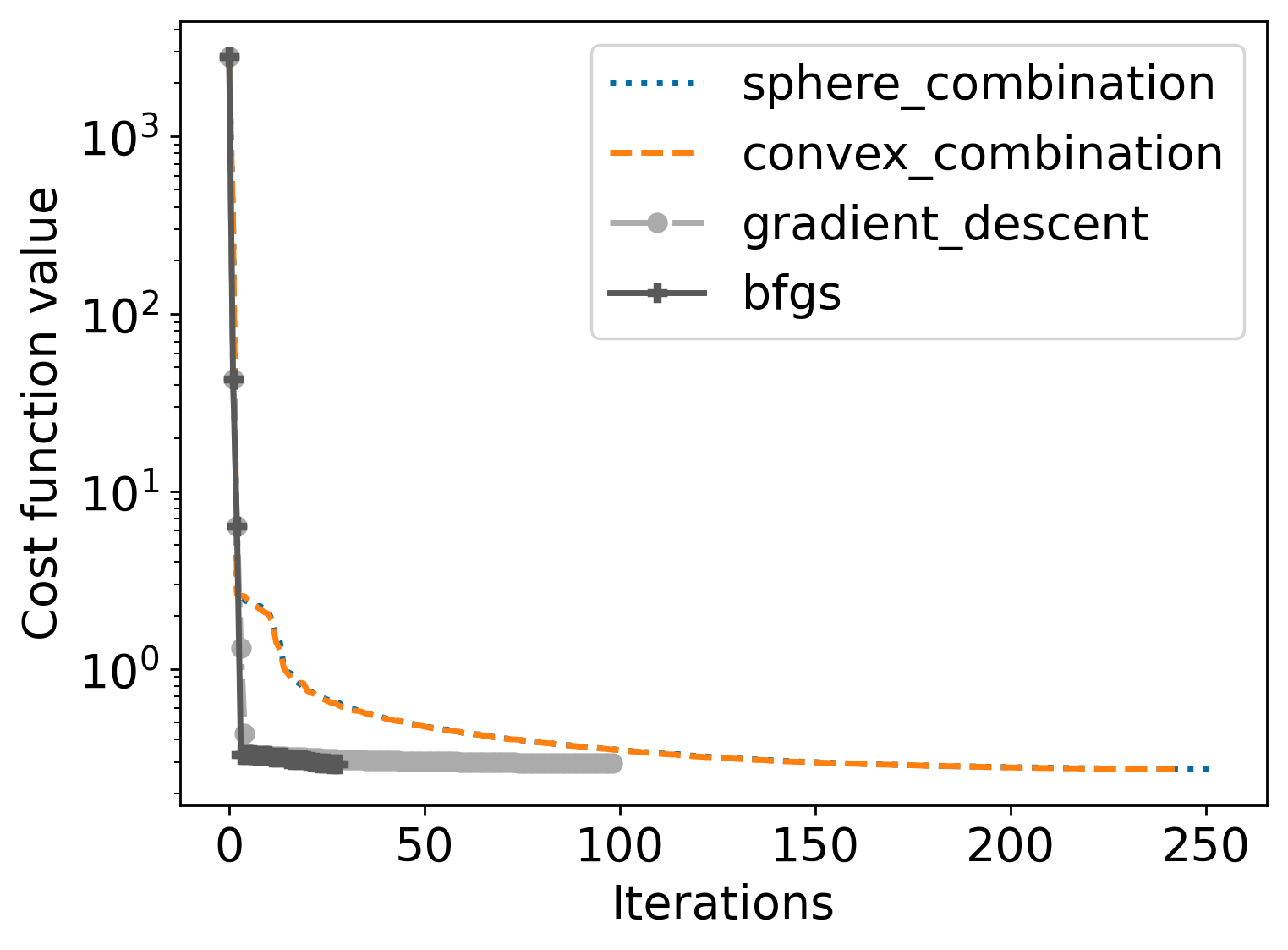}
		\caption{Cost functional.}
	\end{subfigure}%
	\begin{subfigure}{0.333\textwidth}
		\centering
		\includegraphics[width=\textwidth]{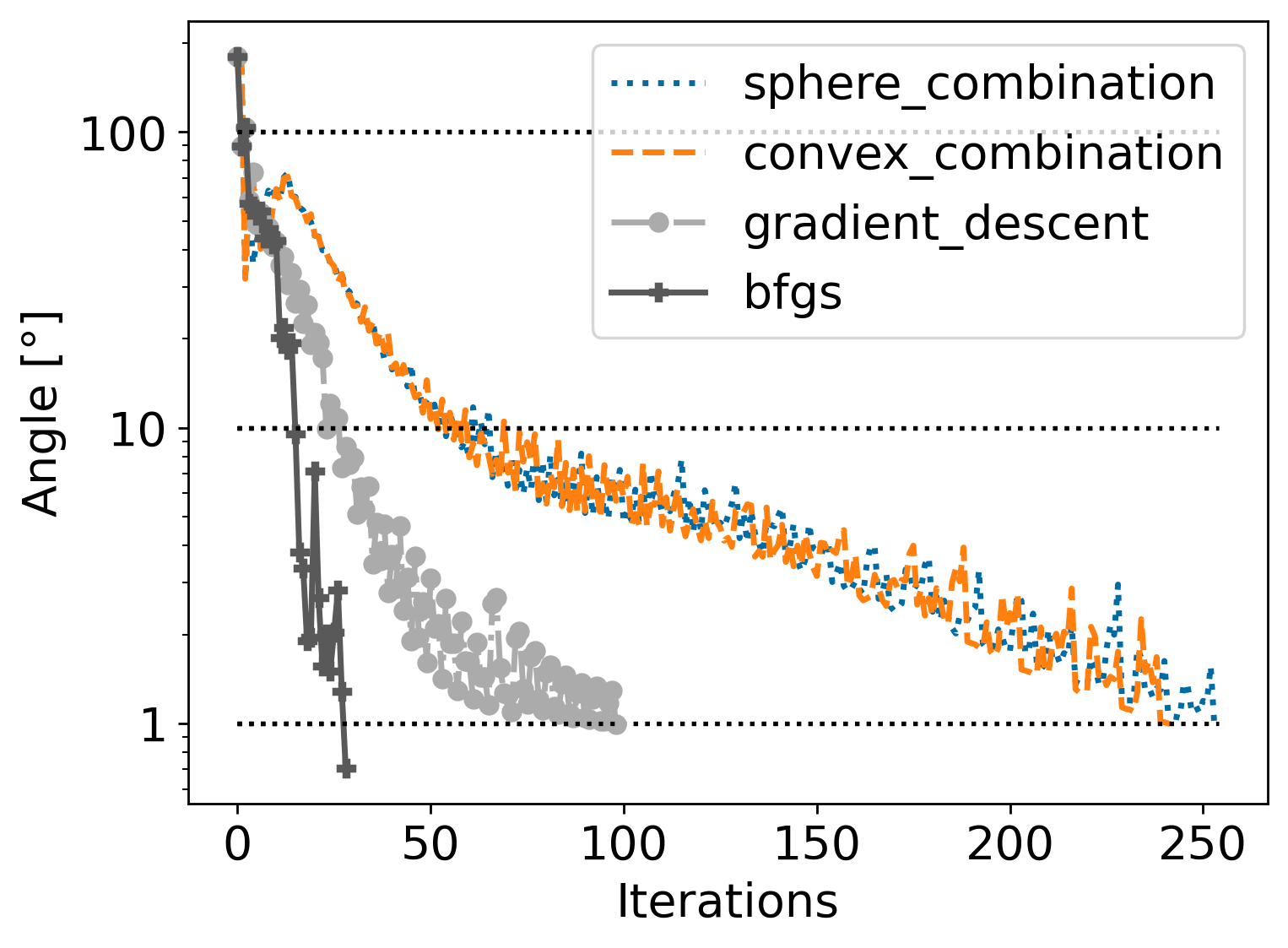}
		\caption{Angle.}
	\end{subfigure}%
	\begin{subfigure}{0.333\textwidth}
		\centering
		\includegraphics[width=\textwidth]{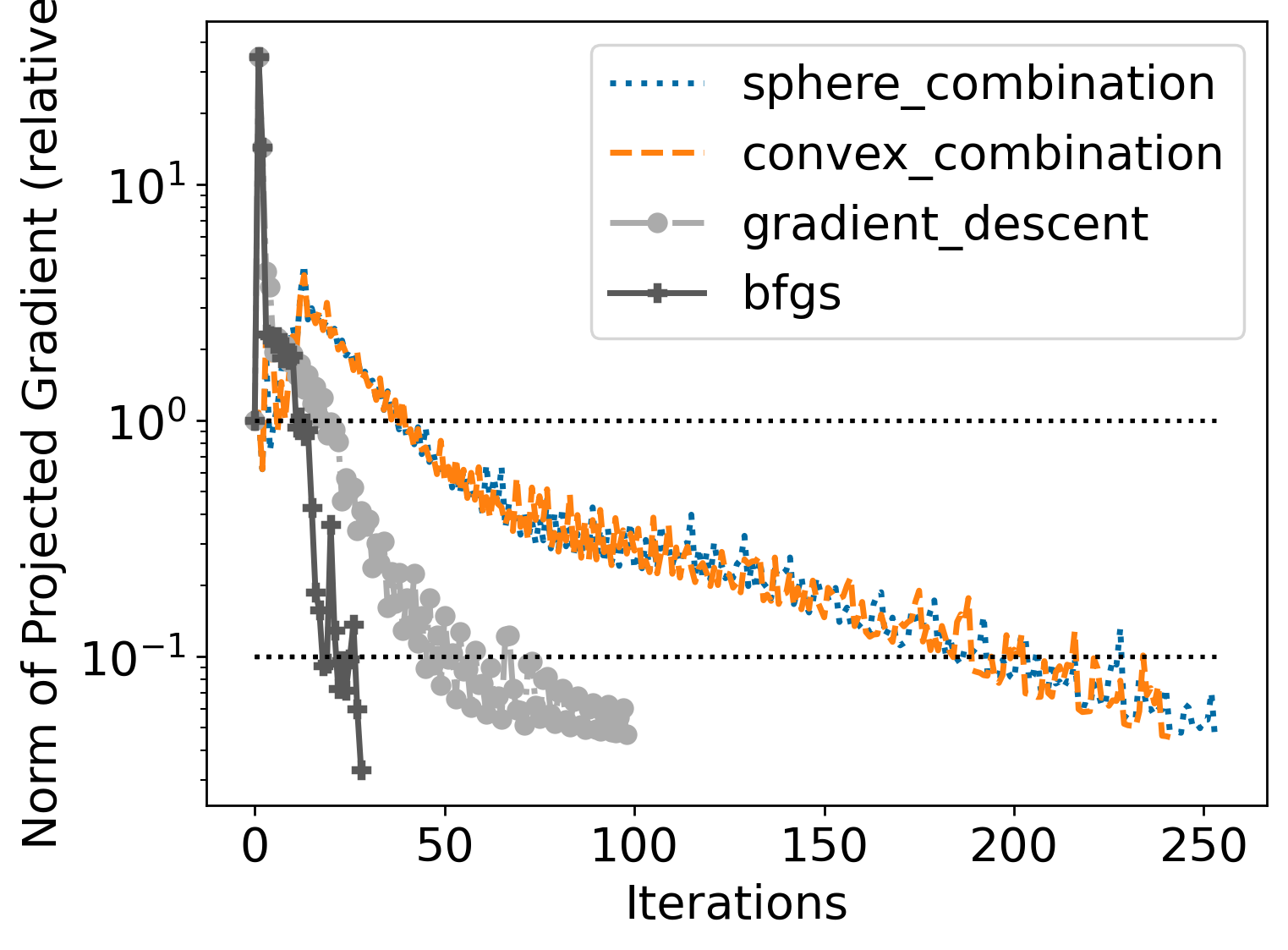}
		\caption{Norm of the projected topological derivative.}
	\end{subfigure}
	\caption{Evolution of the optimization for the pipe bend problem.}
	\label{fig:evolution_pipe_bend}
\end{figure}

\begin{figure}[!t]
	\centering
	\begin{subfigure}{0.25\textwidth}
		\centering
		\includegraphics[width=\textwidth]{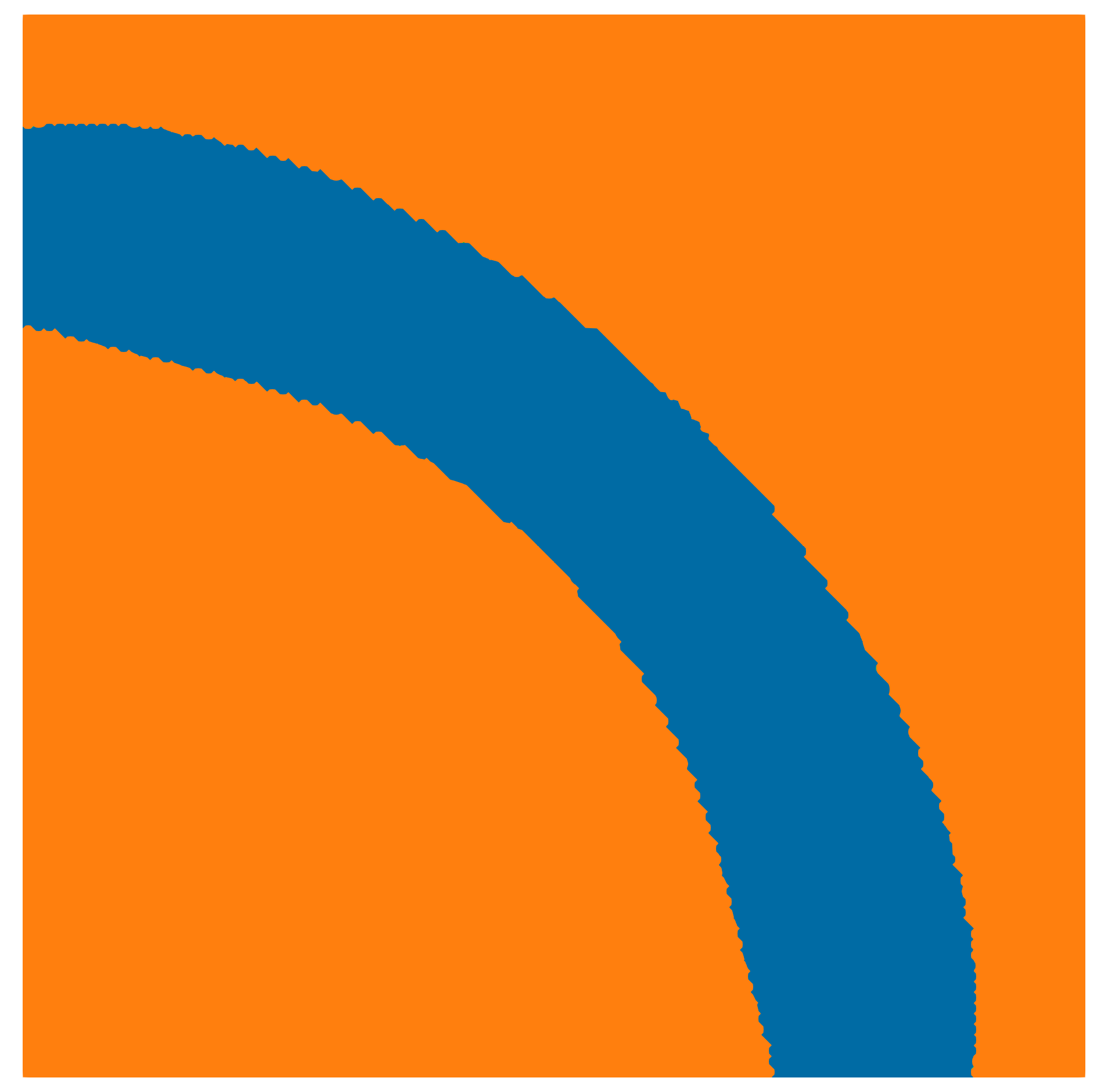}
		\caption{Sphere combination.}
	\end{subfigure}%
	\begin{subfigure}{0.25\textwidth}
		\centering
		\includegraphics[width=\textwidth]{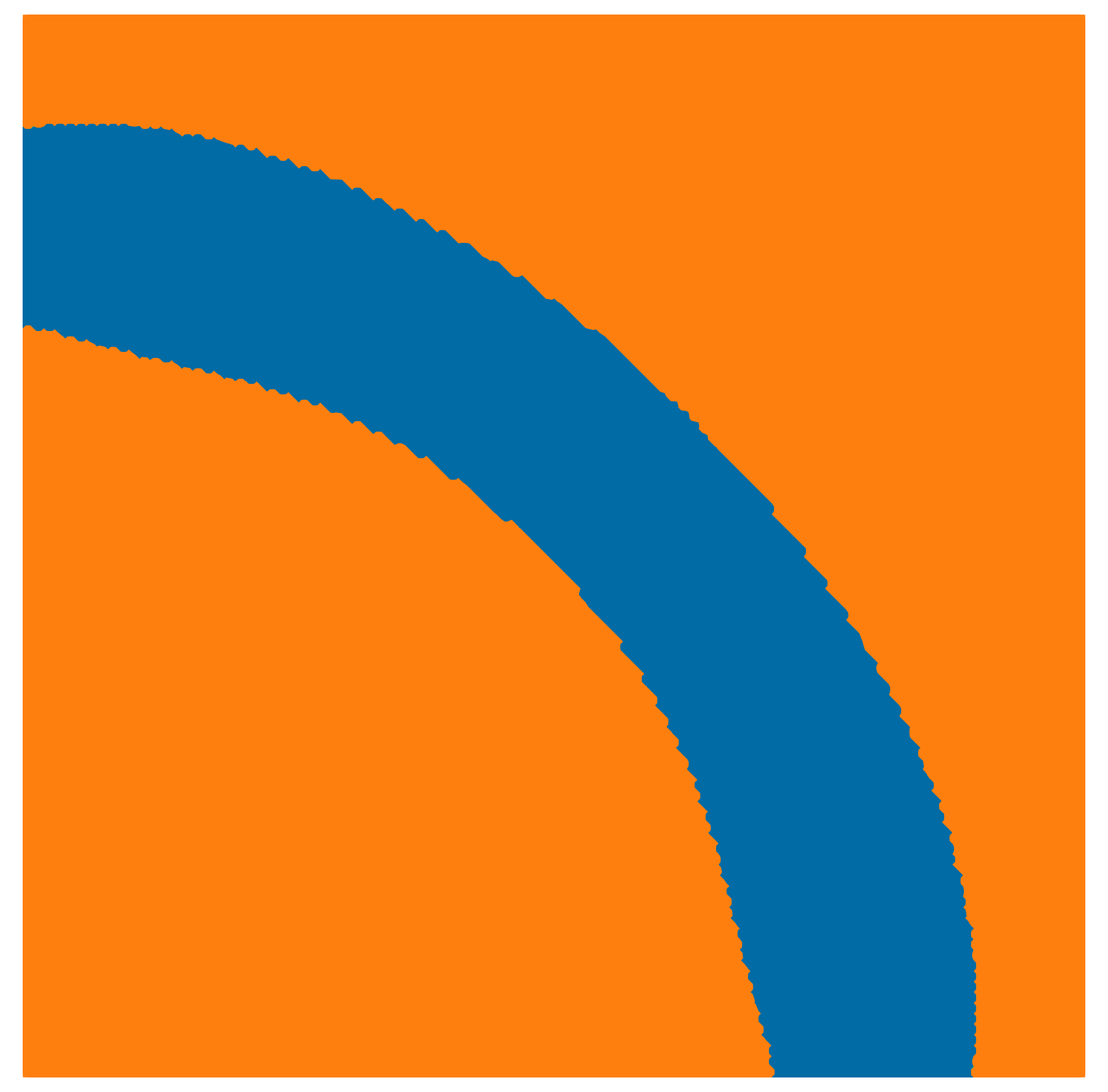}
		\caption{Convex combination.}
	\end{subfigure}%
	\begin{subfigure}{0.25\textwidth}
		\centering
		\includegraphics[width=\textwidth]{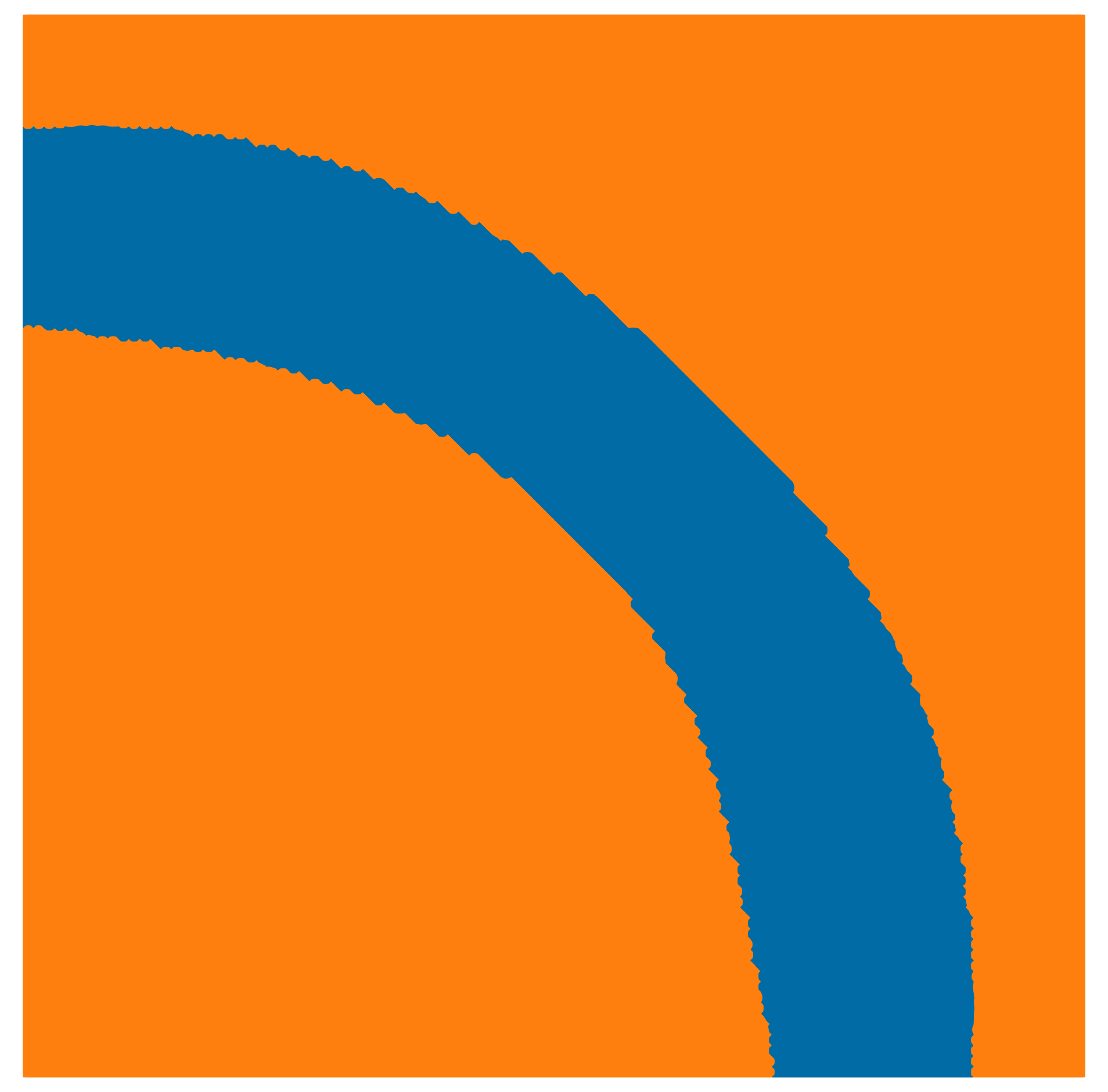}
		\caption{Gradient descent.}
	\end{subfigure}%
	\begin{subfigure}{0.25\textwidth}
		\centering
		\includegraphics[width=\textwidth]{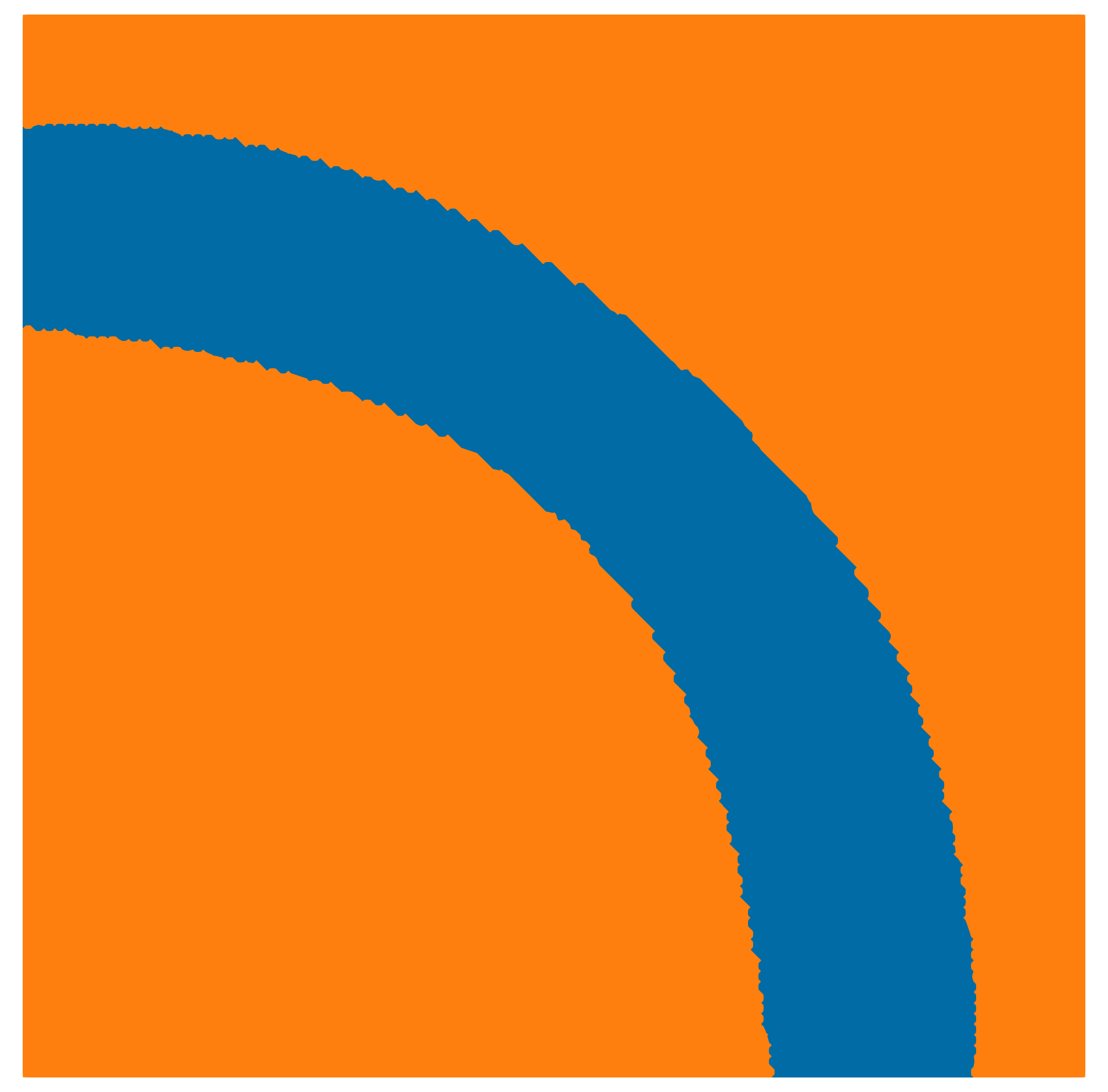}
		\caption{BFGS.}
	\end{subfigure}%
	\caption{Optimized geometries for the pipe bend problem.}
	\label{fig:geom_pipe_bend}
\end{figure}

The evolution of the cost functional, angle criterion, and norm of the projected topological derivative can be seen in Figure~\ref{fig:evolution_pipe_bend}. We observe that the BFGS method substantially outperforms the remaining methods as it only required about 30 iterations to reach the stopping criterion and find a local minimizer of the problem. The gradient descent method performed worse, but still way better than the already established methods: The former required around 100 iterations to satisfy the stopping criterion, whereas the sphere and convex combination methods performed very similarly, requiring about 250 iterations to reach a local minimizer. Considering the optimized geometries obtained by the method, which are depicted in Figure~\ref{fig:geom_pipe_bend}, we observe that all methods converge to a similar minimizer and that all optimized geometries show only the slightest visual differences.

\subsubsection{Rugby Ball}

\begin{figure}[!b]
	\centering
	\begin{tikzpicture}[scale=4]
	\draw [line width=0.4pt] (0,0) -- (1,0) -- (1,1) -- (0,1) -- cycle;
	\draw[pattern=north west lines, pattern color=black, draw=none] (0,0) rectangle (1,1);
	\draw [line width=0.4pt, fill=white] (0.05, 0.05) -- (0.95, 0.05) -- (0.95, 0.95) -- (0.05, 0.95) -- cycle;
	
	\draw [dimen] (0, 1.3) -- (1, 1.3) node {1};
	\draw [dimen] (1.1, 0) -- (1.1, 1) node {1};
	
	\draw [-latex] (0, -0.2) -- (0,0);
	\draw [-latex] (0.1, -0.2) -- (0.1,0);
	\draw [-latex] (0.2, -0.2) -- (0.2,0);
	\draw [-latex] (0.3, -0.2) -- (0.3,0);
	\draw [-latex] (0.4, -0.2) -- (0.4,0);
	\draw [-latex] (0.5, -0.2) -- (0.5,0);
	\draw [-latex] (0.6, -0.2) -- (0.6,0);
	\draw [-latex] (0.7, -0.2) -- (0.7,0);
	\draw [-latex] (0.8, -0.2) -- (0.8,0);
	\draw [-latex] (0.9, -0.2) -- (0.9,0);
	\draw [-latex] (1, -0.2) -- (1,0);
	
	\draw [-latex] (0, 1) -- (0,1.2);
	\draw [-latex] (0.1, 1) -- (0.1,1.2);
	\draw [-latex] (0.2, 1) -- (0.2,1.2);
	\draw [-latex] (0.3, 1) -- (0.3,1.2);
	\draw [-latex] (0.4, 1) -- (0.4,1.2);
	\draw [-latex] (0.5, 1) -- (0.5,1.2);
	\draw [-latex] (0.6, 1) -- (0.6,1.2);
	\draw [-latex] (0.7, 1) -- (0.7,1.2);
	\draw [-latex] (0.8, 1) -- (0.8,1.2);
	\draw [-latex] (0.9, 1) -- (0.9,1.2);
	\draw [-latex] (1, 1) -- (1,1.2);
	
	\draw [-latex] (0, 0) -- (0, 0.2);
	\draw [-latex] (0, 0.2) -- (0, 0.4);
	\draw [-latex] (0, 0.4) -- (0, 0.6);
	\draw [-latex] (0, 0.6) -- (0, 0.8);
	\draw [-latex] (0, 0.8) -- (0, 1);
	
	\draw [-latex] (1, 0) -- (1, 0.2);
	\draw [-latex] (1, 0.2) -- (1, 0.4);
	\draw [-latex] (1, 0.4) -- (1, 0.6);
	\draw [-latex] (1, 0.6) -- (1, 0.8);
	\draw [-latex] (1, 0.8) -- (1, 1);

	\node at (0.5, 0.5) {{\LARGE $D$}};
	\end{tikzpicture}
	\caption{Schematic of the rugby ball problem.}
	\label{fig:rugby}
\end{figure}
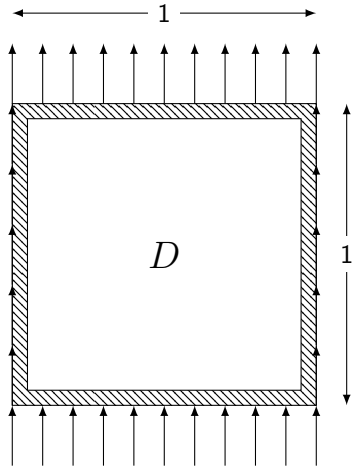

Let us now consider the rugby-ball problem from \cite{Borrvall2003Topology}, which is shown schematically in Figure~\ref{fig:rugby}. Here, the goal is to design an obstacle which minimizes the energy dissipation of the flow. For the volume constraint of the problem, we choose $\text{vol}_\mathrm{des} = \num{0.8}$ and use a regularization parameter of $l = \num{1e4}$ to enforce the volume constraint. For the boundary conditions, we prescribe a constant value of $u_\mathrm{D} = [0, 1]\transposed$ on the entire boundary $\partial \Dsf$.

\begin{figure}[!t]
	\centering
	\begin{subfigure}{0.333\textwidth}
		\centering
		\includegraphics[width=\textwidth]{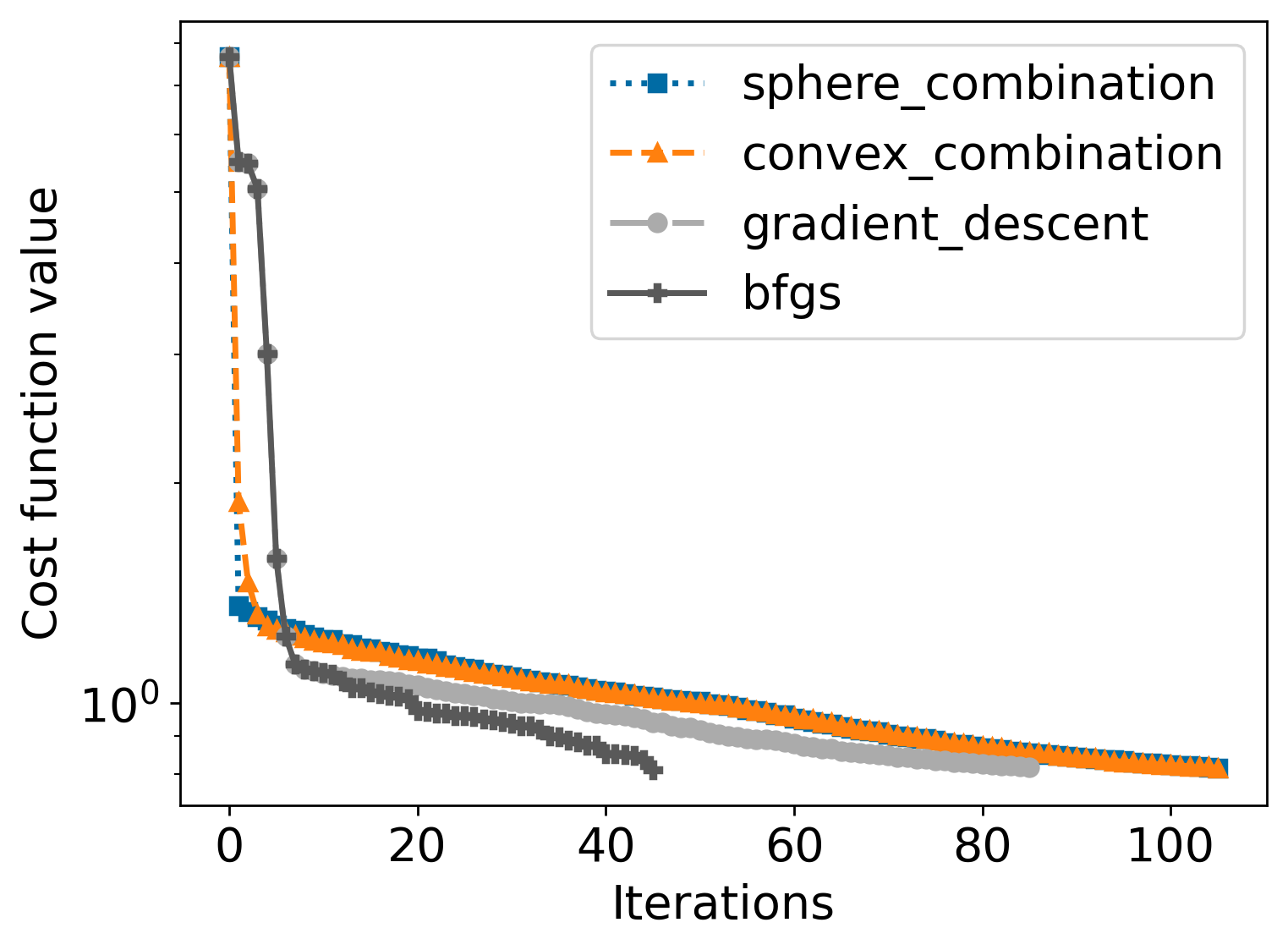}
		\caption{Cost functional.}
	\end{subfigure}%
	\begin{subfigure}{0.333\textwidth}
		\centering
		\includegraphics[width=\textwidth]{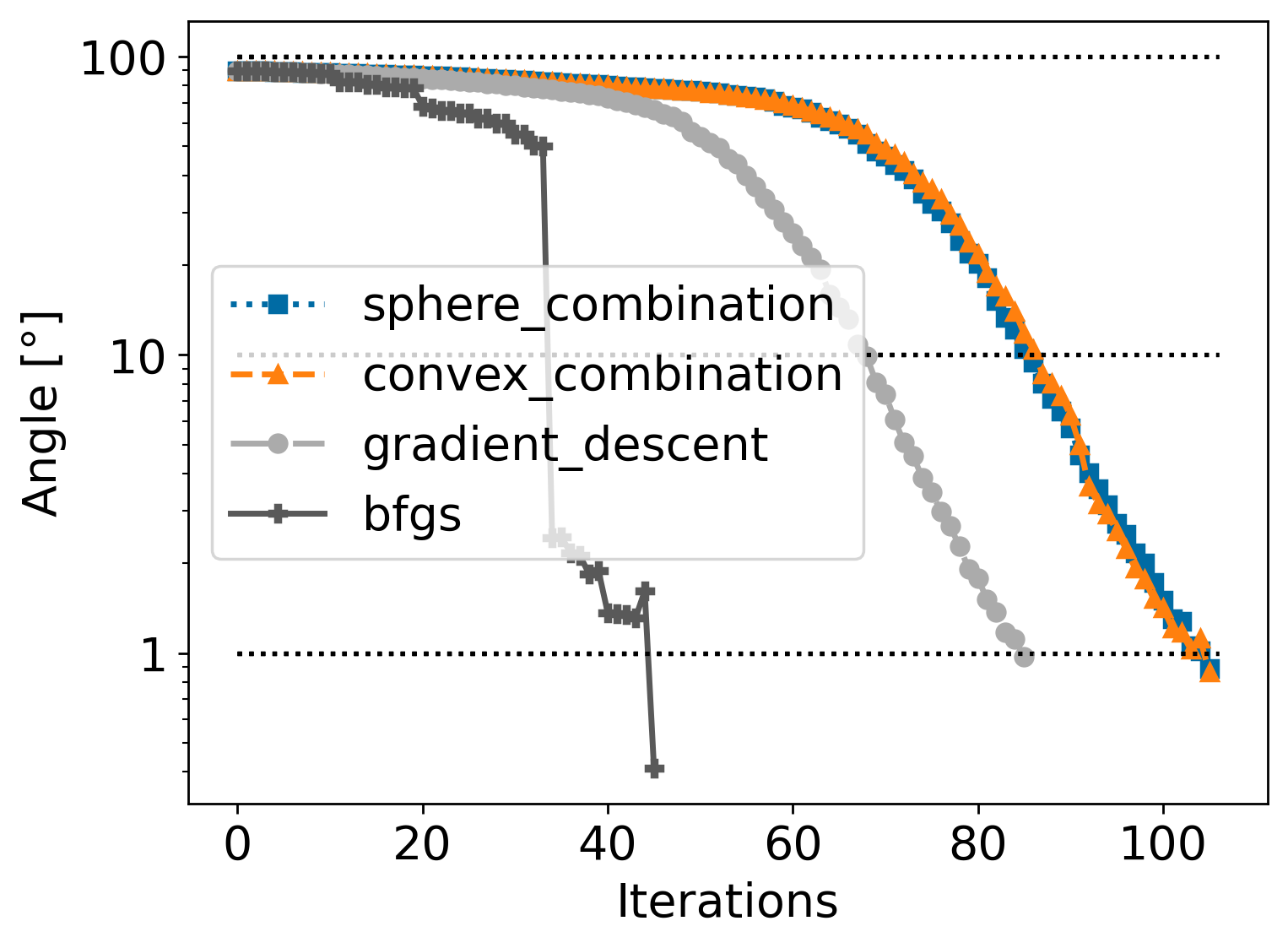}
		\caption{Angle.}
	\end{subfigure}%
	\begin{subfigure}{0.333\textwidth}
		\centering
		\includegraphics[width=\textwidth]{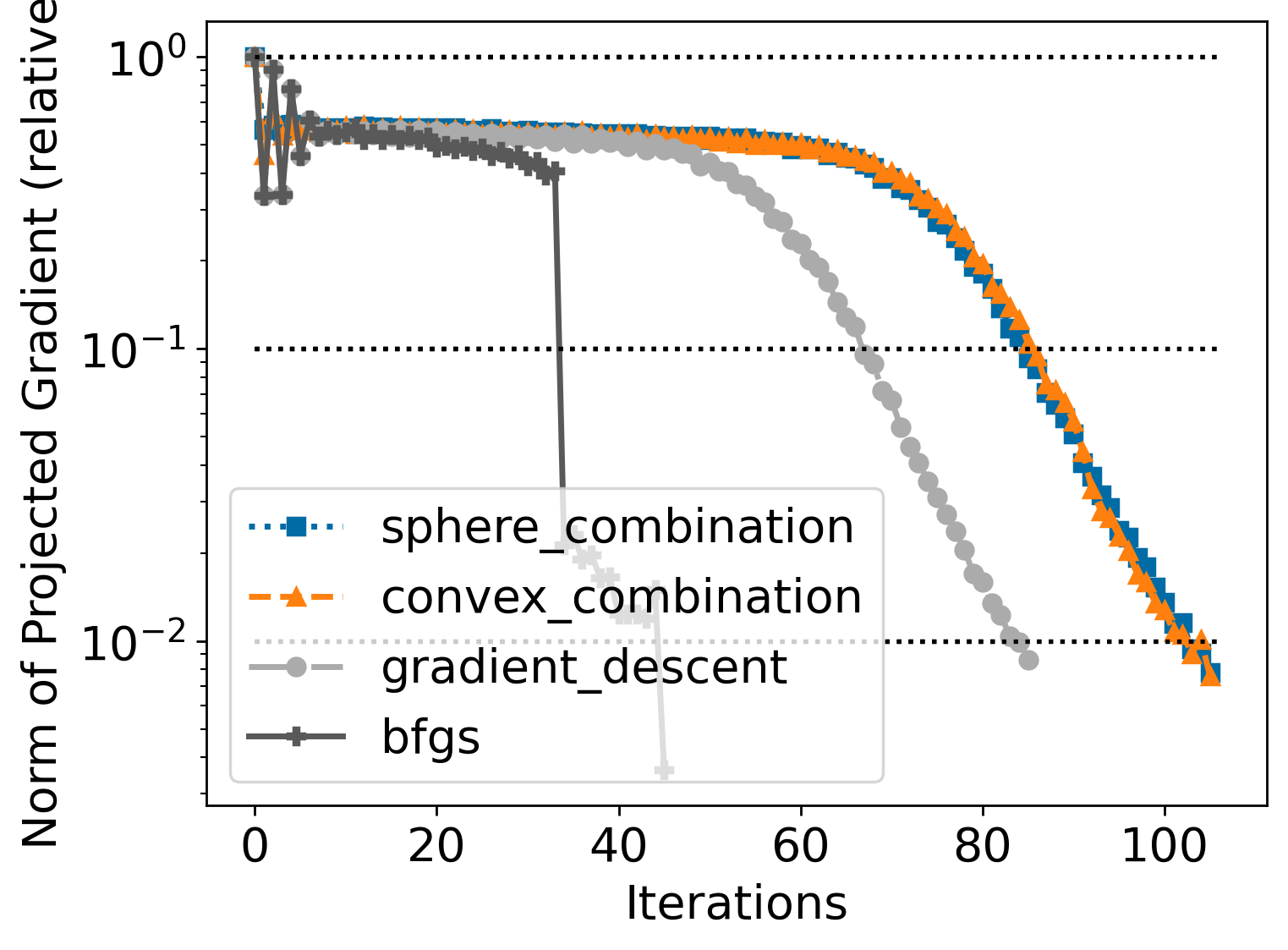}
		\caption{Norm of the projected topological derivative.}
	\end{subfigure}
	\caption{Evolution of the optimization for the rugby ball.}
	\label{fig:evolution_rugby_ball}
\end{figure}

The history of the cost functional, angle between topological derivative and level-set function, and the norm of the projected topological derivative are depicted in Figure~\ref{fig:evolution_rugby_ball}. Here, we can observe that the BFGS method again outperforms all other methods significantly, as it requires slightly less than 50 iterations to find a local minimizer. The gradient descent method shows the second best performance and is able to satisfy the stopping criterion after around 80 iterations. The performance of the sphere and convex combination methods is, again, very similar and both require slightly more than 100 iterations to converge. 

\begin{figure}[!t]
	\centering
	\begin{subfigure}{0.25\textwidth}
		\centering
		\includegraphics[width=\textwidth]{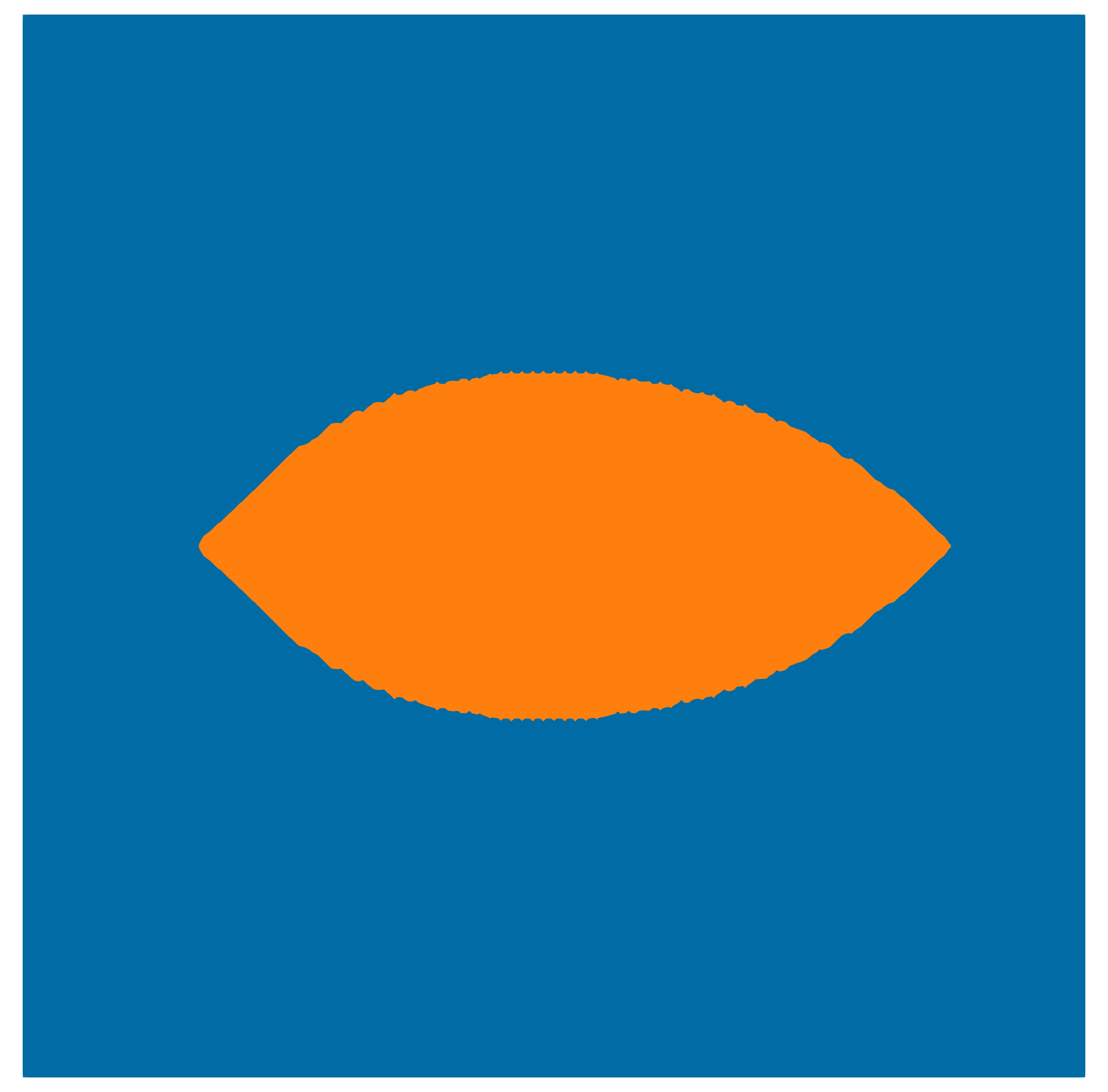}
		\caption{Sphere combination.}
	\end{subfigure}%
	\begin{subfigure}{0.25\textwidth}
		\centering
		\includegraphics[width=\textwidth]{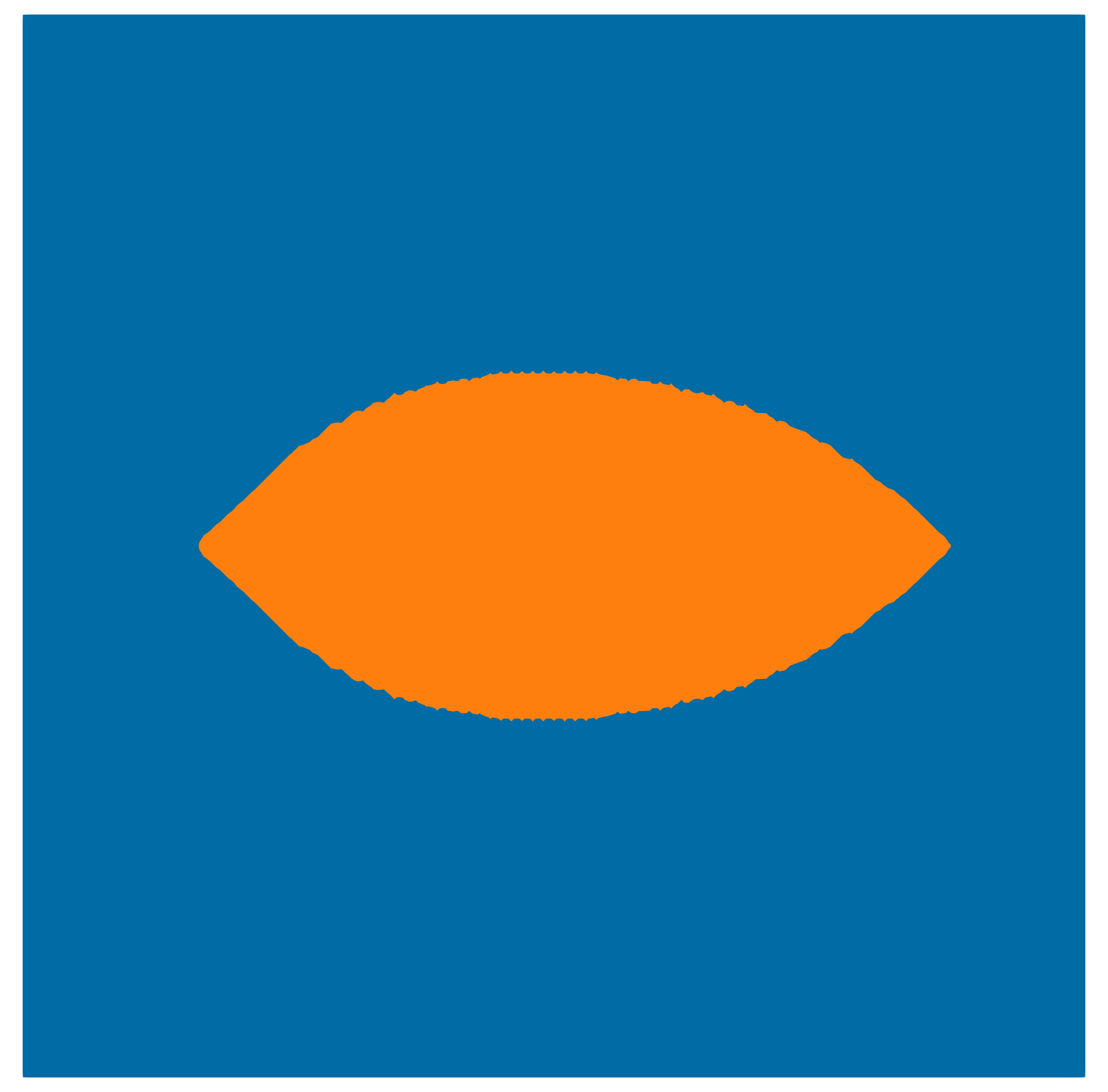}
		\caption{Convex combination.}
	\end{subfigure}%
	\begin{subfigure}{0.25\textwidth}
		\centering
		\includegraphics[width=\textwidth]{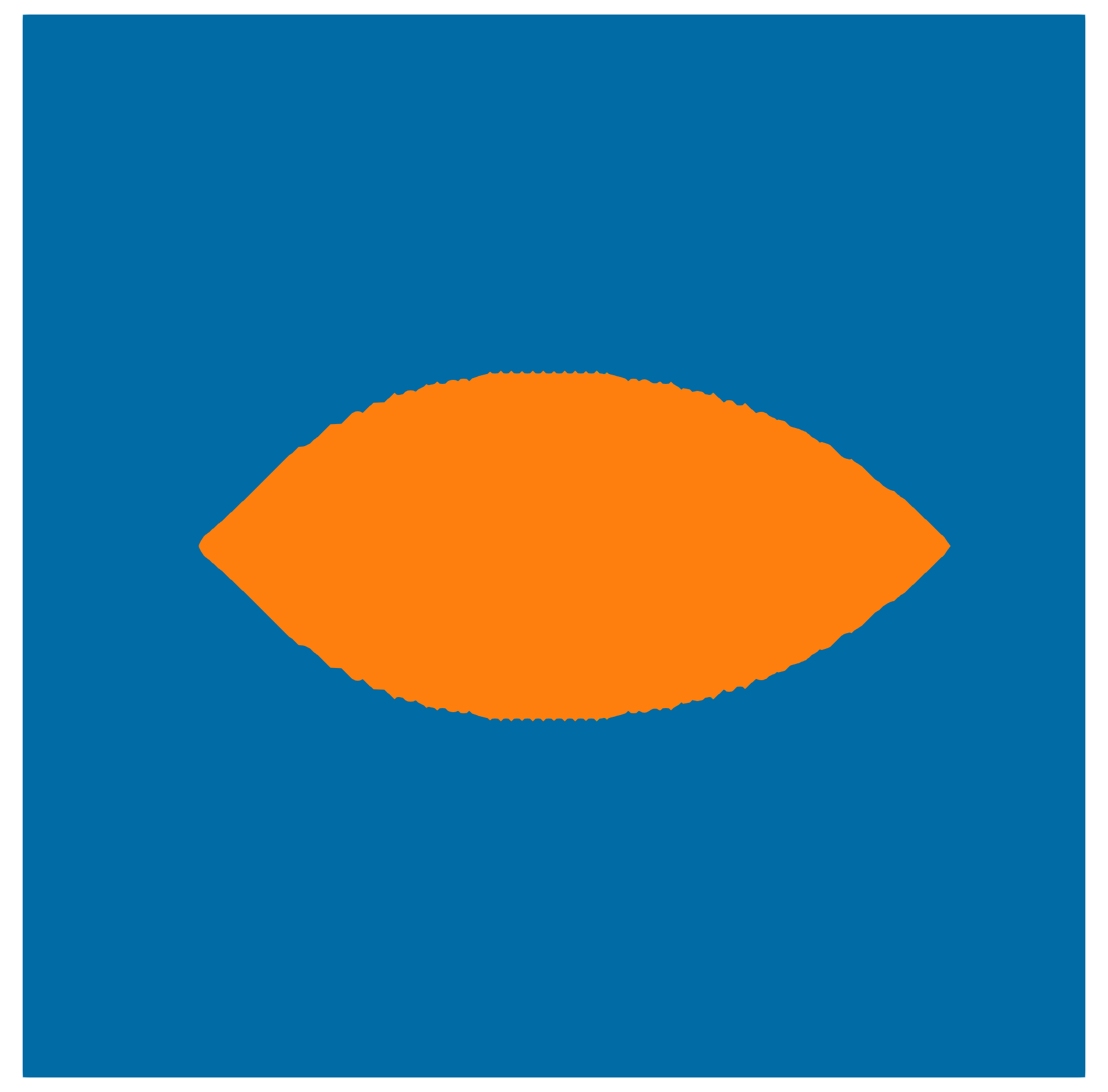}
		\caption{Gradient descent.}
	\end{subfigure}%
	\begin{subfigure}{0.25\textwidth}
		\centering
		\includegraphics[width=\textwidth]{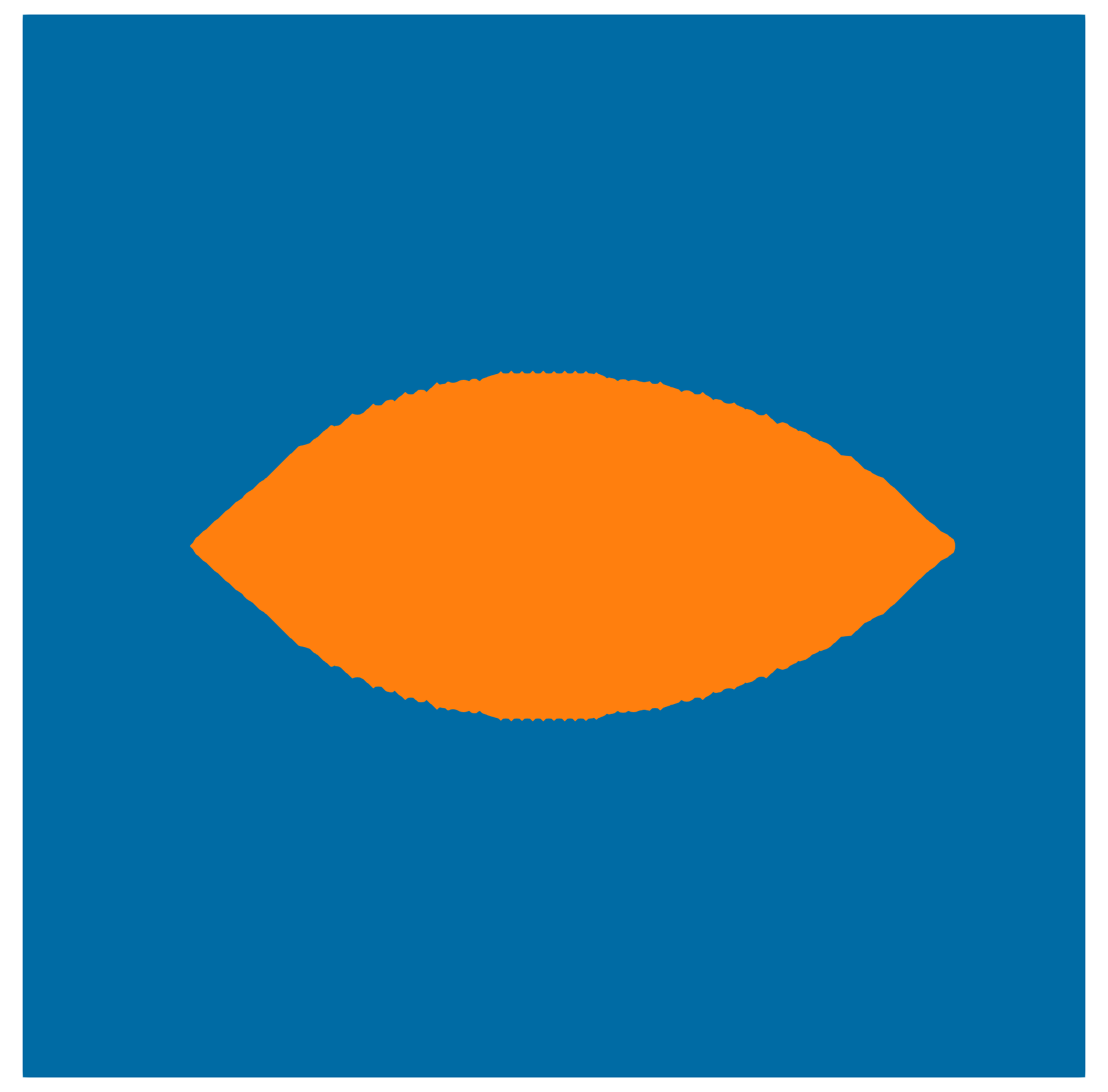}
		\caption{BFGS.}
	\end{subfigure}%
	\caption{Optimized geometries for the rugby ball problem.}
	\label{fig:geom_rugby_ball}
\end{figure}

The obtained optimized geometries are shown in Figure~\ref{fig:geom_rugby_ball}. We see that all methods produce the same desired shape, which is reminiscent of a rugby ball, so that all of them converge to the same local minimizer of the problem.

Altogether, in the context of fluid design optimization, we observe that the BFGS method performs substantially better than the other methods considered in this paper. It reaches the desired stopping criteria with significantly less iterations than the remaining methods. Our findings highlight the efficiency and potential of the proposed BFGS method for solving topology optimization problems.

\section{Conclusion and Outlook}
\label{sec:conclusion}

In this paper, we have presented novel quasi-Newton methods for topology optimization using a level-set method. We recalled the topological derivative, the level-set method for topology optimization, and the widely-used optimization algorithm proposed in \cite{Amstutz2006new}. Then, we presented a new perspective on the evolution equation for evolving the level-set function according to the topological derivative, which enables an interpretation as a classical gradient descent method. This method is the basis for our derivation of quasi-Newton methods for topology optimization and we present a limited-memory BFGS method in this paper. The derivation of the BFGS methods is analogous to the finite-dimensional case and is possible due to the change in perspective described above. We investigated the performance of the proposed gradient descent and BFGS methods on four problem classes: Inverse topology optimization problems constrained by linear and semilinear Poisson problems, compliance minimization in linear elasticity, and the optimization of fluids in Navier-Stokes flow. We compared the results to current state-of-the-art solution algorithms for topology optimization with level-set methods. Our results show that the novel BFGS methods often significantly outperform the other considered methods, requiring substantially less iterations to compute a (local) minimizer. The only exception was the problem of compliance minimization in linear elasticity, where all considered methods performed very similarly, so that the BFGS method performed at least as good as the other, already established methods. All in all, the proposed BFGS methods are efficient and attractive solution algorithms for topology optimization and show great potential for solving such problems.

For future research, there are several interesting directions. One could consider nonlinear conjugate gradient (NCG) methods for topology optimization in analogy to \cite{Blauth2021Nonlinear}, whose derivation is straightforward with the new perspective on the level-set evolution equation presented in this chapter. In fact, these methods are already implemented and available in our open-source software cashocs \cite{Blauth2021cashocs, Blauth2023Version}. However, a thorough numerical analysis of such NCG methods is still required for understanding their performance and behavior. Moreover, a theoretical analysis of the proposed BFGS methods is of great interest. For this, it would be useful to study the properties of the approximate Hessian operator, which the BFGS methods make use of, and its relation to higher order topological derivatives (see, e.g., \cite{Baumann2022Adjoint} and the references therein). Further, it remains an open question why the BFGS methods do not perform as well in the context of compliance minimization in linear elasticity as they do for the other problem classes considered in this paper. Finally, it is, of course, of particular interest to employ our proposed methods for solving practically relevant topology optimization problems, e.g., in the fields of fluid-dynamical optimization or the simulation of fracture evolution.

\section*{Conflict of Interest}

On behalf of all authors, the corresponding author states that there is no conflict of interest. 

\section*{Replication of Results}

The source code for our numerical examples is available publicly on GitHub \cite{Blauth2023Software} and the methods presented in this paper are implemented in our open-source software cashocs \cite{Blauth2021cashocs, Blauth2023Version}.


\bibliographystyle{siamplain}
\bibliography{literature_db.bib}

\end{document}